\documentstyle[amscd, amstex ]{amsart}

\def\NN        {{\Bbb N}}
\def\ZZ         {{\Bbb Z}}
\def\RR         {{\Bbb R}}
\def\CC         {{\Bbb C}}

\def\PP         {{\Bbb P}}

\def\AA         {{\Bbb A}}

\newtheorem{thm}{Theorem}[section]
\newtheorem{lem}[thm]{Lemma}
\newtheorem{cor}[thm]{Corollary}
\newtheorem{pr}[thm]{Proposition}

\theoremstyle{definition}
\newtheorem{rem}[thm]{Remark}
\newtheorem{ex}[thm]{Example}
\newtheorem{defn}[thm]{Definition}

\newcommand{\xsd}{{{X}_{ \Sigma_\Delta}}}

\newcommand{\xst}{{{X}_{\tilde\Sigma}}}

\newcommand{\Hom}{{\rm Hom}}

\newcommand{\spe}{{\rm Spec}}
\newcommand{\key}{\bibitem}

\newcommand{\xs}{{{ X}_{\Sigma}}}

\newcommand\hidot{{\raise1pt\hbox{$\scriptscriptstyle\bullet$}}}
\newcommand\lodot{{\raise.3pt\hbox{$\scriptscriptstyle\bullet$}}}

\newcommand{\co}{{\rm Conv}}
\newcommand{\inte}{{\rm int}}

\newcommand{\ds}{\displaystyle}
\newcommand{\ve}{{\scriptscriptstyle\vee}}
\begin{document}

\title[Deformations  of toric varieties]
{Deformations  of toric varieties \\ via   Minkowski sum
decompositions \\
 of Polyhedral Complexes}
\author{Anvar R. Mavlyutov}
\address {Department of Mathematics, Oklahoma State  University, Stillwater, OK 74078, USA.}
 \email{mavlyutov@@math.okstate.edu}

%Fax: +1(405)-744-8275, mavlyutov@@math.okstate.edu \vskip1cm

% Math Subject Classification

\maketitle

\begin{abstract}
We generalized the construction of deformations of affine toric
varieties of K. Altmann and our previous construction of
deformations of weak Fano toric varieties to the case of arbitrary
toric varieties by introducing the notion of Minkowski sum
decompositions of polyhedral complexes. Our construction embeds
the original toric variety into a higher dimensional toric variety
where the image is given by a prime binomial complete intersection
ideal in Cox homogeneous coordinates.  The deformations are
realized by families of complete intersections. For compact
simplicial toric varieties with at worst Gorenstein terminal
singularities, we show that our deformations span  the
infinitesimal space of deformations by Kodaira-Spencer map. For
Fano toric varieties, we show that their deformations can be
constructed in higher-dimensional Fano toric varieties related to
the Batyrev-Borisov mirror symmetry construction.
\end{abstract}
\tableofcontents
\setcounter{section}{-1}
\section{Introduction.}

In the 90's, Klaus
 Altmann studied   deformations of affine toric varieties
in a series of papers (see \cite{al1}-\cite{al5}).
  His main construction in \cite{al2} combinatorially
  described  the situation when one affine toric variety is embedded
  into
  a higher dimensional affine toric variety as a complete intersection
   given by a regular sequence of  binomial affine functions.
   The embedding of an affine
  toric variety associated with a convex cone $\sigma$
 corresponds to a Minkowski sum decomposition of a  polyhedral slice of
 $\sigma$, which is obtained as an intersection of the cone $\sigma$ with a
 hyperplane.
  There is a natural grading on the coordinate ring of
  of the ambient affine toric variety  coming from a projection of one lattice onto another one,
   and the equations of the complete intersection
  were assumed to have  the same
  degree with respect to this grading.
  This is why  the deformations of these complete
  intersections were called ``homogeneous'' deformations of the affine
  toric variety in \cite{al2}. Altmann's work with some insight from a question of
Bernd Sturmfels was ingenious to realize that Minkowski sums of
polyhedra are related to deformations of affine toric varieties.
However, homogeneous deformations of an affine toric variety in
\cite{al2} did not produce all the unobstructed directions of  the
infinitesimal space of
  deformations of the toric variety by Kodaira-Spencer map.
A later construction in \cite{al5} was supposed to produce
deformations so that infinitesimally they  span all the
unobstructed deformations, but that construction was
``non-toric''.

As it was conjectured in our previous work \cite{m2}, we expect
that all unobstructed deformations of Calabi-Yau complete
intersections in toric varieties can be realized as complete
intersections in higher dimensional toric varieties. Since for
Calabi-Yau hypersurfaces these deformations were induced by
deformations of the ambient toric varieties, which were also
realized by complete intersections in higher dimensional toric
varieties, this motivated our study of  deformations of  toric
varieties.

 Our first deformation construction  in \cite{m1} was obtained by   regluing of open
toric subvarieties of a weak Fano  toric variety by automorphisms
on their intersections. It is not easy if not possible to describe
deformations of abstract algebraic varieties via regluing of open
charts. Another approach is to embed a given variety into another
one as a complete intersection in Cox homogeneous coordinates as
in \cite{m2}  or, in the case of affine varieties, as a complete
intersection given by a regular sequence of affine coordinate
functions as in \cite{al2}. This approach is also not easy as many
such embeddings will only lead to trivial deformations, and a
priori we do not know which embedding will lead to nontrivial
deformations. However, the similarity of the construction of
\cite{m2}, which was obtained through working out some examples of
deformations of Calabi-Yau hypersurfaces in weighted projective
spaces similar to the construction of Sheldon Katz and David
Morrison for a Calabi-Yau in $\PP(1,1,2,2,2)$ in \cite{cofkm}, and
Altmann's construction in \cite{al2} has led us to finding some
compatible language, which turned out to be via the homogeneous
coordinates of toric varieties.

It is not common to use homogeneous coordinates on affine toric
varieties, but this approach produced more deformations of affine
toric varieties than Altmann's construction in \cite{al2}. Namely,
we found deformations of an affine toric variety realized as
complete intersections in another affine toric variety but the
defining equations are in terms of Cox homogeneous coordinates,
and, in general, the embedded variety is not a complete
intersection in terms of affine coordinate functions.  This
feature is very similar to the situation of a non-Cartier
hypersurface in a weighted projective space. Globally it is given
by one equation, and so it is a complete intersection in
homogeneous coordinates even if we restrict it to an affine open
chart of the ambient space. However, since the degree is
non-Cartier, the restriction of the hypersurface to some affine
chart is not a complete intersection in terms of the coordinate
functions of the chart.

In the present paper, we first generalized Altmann's construction
of  deformations of an  affine
  toric variety by embedding   the affine toric variety  into another affine toric
  variety, where the image is given by a  prime binomial complete intersection ideal in terms of
   homogeneous coordinates due to David Cox (see
  \cite{c}). Such ideals are called complete intersection toric
  ideals and have been characterized by K. Fischer and J. Shapiro in \cite{fs} in terms of
  mixed dominating matrices (see Appendix~\ref{s:ap}).
   Altmann's ``homogeneous'' toric
  deformations are described by complete intersections where all the defining equations have the
  same
  degree using the usual definition of degrees of homogeneous coordinates from
  \cite{c}.
After  the generalization of Altmann's construction for affine
toric varieties, we just needed to translate our fan construction
(which made the embedding) from \cite{m2} into the language of
Minkowski sums in order to obtain a more general construction by
introducing {\it Minkowski sum decompositions of polyhedral
complexes}.
 Our previous constructions of deformations of toric blow ups of Fano toric
  varieties in \cite{m1,m2} and Altmann's constructions in \cite{al2,al5} are isomorphic to  particular cases of the current
  one. We expect that by Kodaira-Spencer map our deformations of toric varieties should
  span all the unobstructed directions
 of the   infinitesimal space of deformations of a toric variety
 and verify this for compact simplicial    toric varieties with at worst Gorenstein terminal singularities.

   Once there is a sufficient
number of constructed deformations of an algebraic variety  so
that they span all the unobstructed deformations  by
Kodaira-Spencer map, the next step is to combine these
deformations in bigger families corresponding to the irreducible
components of the moduli of the algebraic variety. In general,
explicit solutions to the deformation problem for algebraic
varieties may not be possible, but for toric varieties we believe
that solution is possible  in terms of complete intersections in
toric varieties. Different irreducible components of the moduli
should be represented by families of complete intersections in
different toric varieties. Here, we show how to combine
deformation families  for Fano toric varieties, representing parts
of the irreducible components of the moduli.

Applying our general deformation construction to Gorenstein Fano
toric varieties $X_\Delta$ corresponding to reflexive polytopes
$\Delta$
 we discover that a large
number of deformations can be combined together in one family.
More precisely, we constructed  $(kl(\Delta)-k)$-parameter
deformation families corresponding to Minkowski sum decompositions
$\Delta^*=\Delta_0^*+\Delta_1^*+\cdots+\Delta_k^*$ of the
reflexive polytope dual to $\Delta$, where $l(\Delta)$ is the
number of lattice points in $\Delta$. The nice feature of this
construction is that the ambient toric variety of these
deformations is again a Gorenstein Fano toric variety, associated
with the dual of the Cayley cone corresponding to polytopes
$\Delta_0^*,\Delta_1^*,\cdots,\Delta_k^*$ used in the
Batyrev-Borisov Mirror Symmetry construction (see \cite{bb}). We
will explore this story in a parallel paper \cite{m3}.

{\it Acknowledgments.} Our key constructions in
Sections~\ref{s:1}--\ref{s:darb} via Cox homogeneous coordinates
and Minkowski sum decompositions of polyhedral complexes have
appeared as a result of inspiration to present a new result  at
the Valley Geometry Seminar at UMass, Amherst, in May 2007. We are
grateful to  David Cox and Evgeny Materov for inviting to the
seminar  and for their support. We also thank Paul Horja for
suggesting (during a local seminar) that condition $(**)$ on the
Minkowski sum decomposition is not necessary.

\section{Toric  varieties, subvarieties and morphisms in
homogeneous coordinates.}\label{s:toric}

In this section, we review the categorical quotient presentation
of a toric variety in terms of homogeneous coordinates from
\cite{c}. Then show that the quotient presentation extends to
subvarieties of toric varieties. We also define  complete
intersections  in arbitrary toric varieties in homogeneous
coordinates and describe toric morphisms between toric varieties
as monomial maps in homogeneous coordinates.

 Let $\xs$ be a $d$-dimensional toric variety associated
with a finite rational polyhedral fan $\Sigma$ in a
$d$-dimensional real vector space $N_\RR$ with lattice $N$ (see
\cite{d,f} for the definition). Denote by $\Sigma(1)$ the finite
set of the one-dimensional cones in  $\Sigma$, corresponding to
the torus invariant divisors $D_\rho$ in $\xs$. For a
1-dimensional cone $\rho\in\Sigma(1)$, its primitive lattice
generator will be denoted $v_\rho$. From the work of David Cox in
\cite{c}, every toric variety can be described as a categorical
quotient of a Zariski open subset of an affine space by a subgroup
of a torus. For simplicity we assume that the toric variety will
not have torus factors corresponding to the condition that
$\Sigma(1)$ span $N_\RR$. Consider the polynomial ring
$\CC[x_\rho\mid\rho\in\Sigma(1)]$, called the homogeneous
coordinate ring of the toric variety $\xs$, and the corresponding
affine space
$\CC^{\Sigma(1)}=\spe(\CC[x_\rho\mid\rho\in\Sigma(1)])$, where (by
abuse of notation) $\Sigma(1)$ also denotes the number of rays in
the fan $\Sigma(1)$. Let $B(\Sigma)=\langle \prod_{\rho
 \not\subseteq\sigma} x_\rho\mid\sigma\in\Sigma\rangle$ be the
ideal in $\CC[x_\rho\mid\rho\in\Sigma(1)]$. This ideal determines
a Zariski closed (exceptional) set $Z(\Sigma)$ in
$\CC^{\Sigma(1)}$.
 This set is
invariant under the diagonal group action of the subgroup
$$G(\Sigma)=\Biggl\{(\mu_\rho)_{\rho\in\Sigma(1)} \in(\CC^*)^{\Sigma(1)}\,\biggr|\,\prod_{\rho\in\Sigma(1)} \mu_\rho^{\langle
u,v_\rho \rangle}=1\, \forall\, u\in M:=\Hom(N,\ZZ)\Biggr\}$$ of
the torus $(\CC^*)^{\Sigma(1)}$ on the affine space
$\CC^{\Sigma(1)}$. Then by Theorem 2.1 in \cite{c}, the toric
variety  $\xs$ is naturally isomorphic to the categorical quotient
$$(\CC^{\Sigma(1)}\setminus Z(\Sigma))/G(\Sigma).$$
 As shown in Theorem 5.1.10 in \cite{torvar}, this isomorphism    arises from
the toric morphism $\pi:\CC^{\Sigma(1)}\setminus
Z(\Sigma)\rightarrow\xs $, constant on  $G(\Sigma)$-orbits. The
variety $\CC^{\Sigma(1)}\setminus Z(\Sigma)$ can be described as a
toric variety $X_{\bar\Sigma}$ corresponding to the fan
$$\bar\Sigma =\{\tau\mid\tau\preceq \bar\sigma \,\,{\rm for \,\,some}\,\,\sigma\in\Sigma\}$$
in $\RR^{\Sigma(1)}$, where $\bar\sigma =\sum_{\rho \subseteq
\sigma}\RR_{\ge0}\cdot e_\rho\subseteq\RR^{\Sigma(1)}
=\bigoplus_{\rho\in\Sigma(1)} \RR \cdot e_\rho$ corresponds to
$\sigma\in\Sigma$ and $\{e_\rho\mid\rho\in \Sigma(1)\}$ is the
standard basis of $\ZZ^{\Sigma(1)}$. Then the toric morphism
  $\pi$ is induced by the lattice morphism
$\bar\pi:\ZZ^{\Sigma(1)}\rightarrow N$, $\bar\pi(e_\rho)=v_\rho$,
sending cones of  $\bar\Sigma$ into cones of $\Sigma$.
  When $\Sigma$
is simplicial the above categorical quotient is a geometric
quotient, meaning that points in $\xs$ correspond to
$G(\Sigma)$-orbits in $\CC^{\Sigma(1)}\setminus Z(\Sigma)$.

 The quotient
presentation of a toric variety is important because it allows us
to work with closed subvarieties of the toric variety. The
homogeneous coordinate ring $S(\Sigma):=\CC[x_\rho\mid
\rho\in\Sigma(1)]$ of the toric variety $\xs$ is graded by the
class group
$${\rm Cl}(\Sigma)\simeq \Hom(G(\Sigma),\CC^*),$$
given by the short exact sequence $$0\rightarrow M\rightarrow
\ZZ^{\Sigma(1)}\rightarrow{\rm Cl}(\Sigma)\rightarrow0,$$ where
$m\mapsto(\langle m,v_\rho\rangle)_{\rho\in\Sigma(1)}$. The degree
of a monomial in $S(\Sigma)$ is defined by
$$\deg\Biggl(\prod_{\rho\in\Sigma(1) }
x_\rho^{a_\rho}\Biggr)= [(a_\rho)_{\rho\in\Sigma(1)}]\in {\rm
Cl}(\Sigma).$$

 By Proposition 5.2.8 in \cite{torvar}, all closed
subvarieties of $\xs$ correspond to homogeneous ideals $I\subseteq
S(\Sigma)$: $${\bf V}(I)= \{p\in\xs\mid \exists\,
x\in\pi^{-1}(p)\,{\rm with }\,\,f(x)=0\,\forall\, f\in I\}.$$ For
compact simplicial toric varieties $\xs$, the closed subvarieties
${\bf V}(I)$ bijectively correspond to radical homogeneous ideals
$I\subseteq B(\Sigma)$ by Proposition 2.4 in \cite{c}.

 As in the classical case of complete intersections in an affine or  projective space,
  we   define an (ideal-theoretic) {\it complete intersection} in a toric variety $\xs$ (in homogeneous coordinates) as
  a
  subvariety ${\bf V}(I)$ corresponding to a radical homogeneous ideal $I\subseteq
S(\Sigma)$ generated
  by a  regular sequence of homogeneous polynomials $f_1,\dots,f_k\in S(\Sigma)$ such that $k=\dim
  \xs-\dim {\bf V}(I)$. We will say that the complete intersection is given by the
  equations $f_i(x)=0$, for $i=1,\dots,k$. Quasismooth complete
  intersections in simplicial toric varieties have been studied in
  \cite{m0}.

\begin{rem} If the toric variety $\xs$ is simplicial, then the
group $G(\Sigma)$ has zero dimensional stabilizers on
$\CC^{\Sigma(1)}\setminus Z(\Sigma)$ with all orbits closed (see
\cite{bc}), and, in particular, if the radical ideal $I$ of the
affine subvariety ${\bf V}_a(I)\subset\CC^{\Sigma(1)}$ is
generated  by a regular sequence $f_1,\dots,f_k\in S(\Sigma)$ of
Cox homogeneous polynomials, then the closed subvariety ${\bf
V}(I)\subset \xs$ is of pure codimension $k$. In the nonsimplicial
case, the varieties ${\bf V}_a(I)$ and ${\bf V}(I)$ may have
different codimension as Example 5.2.9 in \cite{torvar} shows.
\end{rem}

 The above
presentation ${\bf V}(I)$ of closed subvarieties  in $\xs$ is good
for set-theoretic  purposes, but it is missing the
scheme-theoretic  structure, essential for studying of their
deformations. We now show that a closed subvariety in a toric
variety $\xs$ admits a categorical quotient presentation, which
equips the closed subvariety with the scheme structure in terms of
its homogeneous ideal.

\begin{thm}\label{t:quotient}  Let $\xs$ be a toric variety without torus factors. Let $Y= {\bf V}(I)\subseteq \xs$ be a
closed subvariety corresponding to a radical homogeneous ideal
$I\subseteq S(\Sigma)$.  Then $Y$ is
 naturally isomorphic
to the categorical quotient $${\bf
V}_a(I)\cap(\CC^{\Sigma(1)}\setminus Z(\Sigma))/G(\Sigma),$$ where
${\bf V}_a(I)\subseteq\CC^{\Sigma(1)}$  denotes the affine
subvariety defined by $I$. If the fan $\Sigma$ is simplicial, this
quotient is geometric.
\end{thm}

\begin{pf}
The toric variety $X_{\bar\Sigma} =\CC^{\Sigma(1)}\setminus
Z(\Sigma)$ has an affine open cover by
$$X_{\bar\sigma}:=\spe\bigl(\CC\bigl[\bar\sigma^\ve\cap\ZZ^{\Sigma(1)}\bigr]\bigr)$$
which is $G(\Sigma)$-invariant and the toric variety $\xs$ has an
affine open cover by $X_\sigma={\rm Spec}(\CC[{ {\sigma}^\ve}\cap
M])$, for $\sigma\in\Sigma$. Moreover, by the proof of
Theorem~5.1.10 in \cite{torvar}, we know that
$$\pi|_{X_{\bar\sigma}}:X_{\bar\sigma}=\pi^{-1}(X_\sigma)\rightarrow X_\sigma $$ is a good categorical quotient, corresponding to the fact
that on the level of the coordinate rings
$\pi_\sigma:=\pi|_{X_{\bar\sigma}}$ induces an isomorphism $$
\CC[{ {\sigma}^\ve}\cap M]\simeq
\CC\bigl[\bar\sigma^\ve\cap\ZZ^{\Sigma(1)}\bigr]^{G(\Sigma)}\simeq
(S(\Sigma)_{x^{\hat\sigma}})^{G(\Sigma)},$$ where
$x^{\hat\sigma}:= \prod_{\rho
 \not\subseteq\sigma} x_\rho$.

The morphism $\pi$ induces a surjective morphism
$$\pi_I:{\bf V}_a(I)\cap X_{\bar\Sigma}\rightarrow
Y=\textbf{V}(I). $$ The closed subvariety $Y\subset\xs$ has an
affine open cover by $Y\cap X_\sigma$, for $\sigma\in\Sigma$, and
$\pi_I^{-1}(Y\cap X_\sigma)={\bf V}_a(I)\cap X_{\bar\sigma}$.
 Denote $\pi_{I,\sigma}:=\pi_I|_{\pi_I^{-1}(Y\cap X_\sigma)}$. Then we have
commutative diagrams $$ \begin{CD}
{\bf V}_a(I)\cap X_{\bar\sigma}&\,\hookrightarrow  &X_{\bar\sigma}\\
@VV{ \pi}_{I, \sigma}V  @VV\pi_\sigma V \\
Y\cap X_\sigma & \hookrightarrow &\,\,X_\sigma,
\end{CD} $$
where the horizontal morphisms are embeddings and vertical
morphisms ${ \pi}_{I,\sigma}$, $ \pi_\sigma$  are surjections. The
corresponding diagram  of coordinate rings is $$ \begin{CD}
S(\Sigma)_{x^{\hat\sigma}}/I_{x^{\hat\sigma}}  @<<<   S(\Sigma)_{x^{\hat\sigma}}\\
@A\pi^*_{I,\sigma} AA  @A\pi^*_\sigma AA \\
\CC[Y\cap X_\sigma] @<<<     \CC[{ {\sigma}^\ve}\cap M],
\end{CD}   $$
where the horizontal arrows are surjections and the vertical ones
are injections. Since the image of $\pi^*_\sigma$ is
$(S(\Sigma)_{x^{\hat\sigma}})^{G(\Sigma)}$, by commutativity of
the diagram we get  isomorphism $$\CC[Y\cap X_\sigma]\simeq
(S(\Sigma)_{x^{\hat\sigma}}/I_{x^{\hat\sigma}})^{G(\Sigma)}$$
induced by $\pi^*_{I,\sigma}$. Hence, $\pi_{I,\sigma}$ and $\pi_I$
are good categorical quotients by Propositions~5.0.9 and 5.0.12 in
\cite{torvar}.
\end{pf}

The proof of the above theorem tells us that   the closed
subvariety $Y={\bf V}(I)\subseteq\xs$    locally corresponds to
the rings of invariants
$(S(\Sigma)_{x^{\hat\sigma}}/I_{x^{\hat\sigma}})^{G(\Sigma)}$. In
an affine toric chart $X_\sigma\subset \xs$, it is easier to work
with the following equivalent description assuming that the
dimension of $\sigma$ is maximal.

\begin{pr}\label{p:local} Let $\sigma\in\Sigma$ be such that $\dim\sigma=\dim
N_\RR$ and let $\sigma(1)$ be the set of edges of $\sigma$. Let
$\psi: S(\Sigma)=\CC[x_\rho\mid\rho\in\Sigma(1)]\rightarrow
S(\sigma):=\CC[x_\rho\mid\rho\in\sigma(1)]$ be the ring
homomorphism given by evaluation at $x_\rho=1$ for
$\rho\not\subseteq\sigma$.
 Then
$\psi$ induces natural ring isomorphisms:
$$(S(\Sigma)_{x^{\hat\sigma}})^{G(\Sigma)}\simeq
S(\sigma)^{G(\sigma)},\qquad
 (S(\Sigma)_{x^{\hat\sigma}}/I_{x^{\hat\sigma}})^{G(\Sigma)}\simeq(S(\sigma)/\psi(I))^{G(\sigma)},$$
where  $G(\sigma):=\bigl\{(\mu_\rho)
\in(\CC^*)^{\sigma(1)}\mid\prod_{\rho\in\sigma(1)}
\mu_\rho^{\langle u,v_\rho \rangle}=1\, \forall\, u\in M\bigr\}$.
\end{pr}

\begin{pf} The affine toric variety $X_\sigma$ can be viewed as a
toric variety corresponding to the fan consisting of the faces of
the cone $\sigma$. Hence, we can represent  $X_\sigma$ by the
categorical quotient. Note that the corresponding exceptional set
in $\CC^{\sigma(1)}$ is empty since there are no rays of the fan
of $X_\sigma$ not lying in  $\sigma$. By the construction, we have
the surjective toric morphism
$$\varphi:\CC^{\sigma(1)}=\spe\bigl(\CC\bigl[\ZZ_{\ge0}^{\sigma(1)}\bigr]\bigr)\rightarrow
X_\sigma=\spe(\CC[{ {\sigma}^\ve}\cap M])$$ induced by the
morphism  of lattices:
$$\bar\varphi:\ZZ^{\sigma(1)}=\bigoplus_{\rho\in\sigma(1)} \ZZ \cdot
e_\rho \rightarrow N,\quad e_\rho\mapsto v_\rho,$$  where
$\{e_\rho\mid\rho\in \sigma(1)\}$ is the standard basis of
$\ZZ^{\sigma(1)}$. This fits into a commutative diagram of toric
morphisms
$$
\begin{array}{cccc}
  \CC^{\sigma(1)} &  \hspace{-0.3cm} \longrightarrow & \hspace{-0.4cm} X_{\bar\sigma} &\\
    \downarrow & \swarrow\quad  &   \\
  \,\, X_\sigma, &&&
\end{array}$$
induced by the composition $\bar\varphi:\ZZ^{\sigma(1)}
\hookrightarrow\ZZ^{\Sigma(1)} \rightarrow N$, where the first
morphism is  the natural inclusion and the second one is the same
one that induces the toric morphism
$\pi:X_{\bar\Sigma}\rightarrow\xs$. On the level of coordinate
rings
  this translates into the diagram
$$
\begin{array}{cccc}
  S(\sigma)&  \hspace{-1cm} \longleftarrow & \hspace{-0.4cm} S(\Sigma)_{x^{\hat\sigma}}&\\
    \uparrow & \nearrow\quad  &   \\
  \,\, \CC[{ {\sigma}^\ve}\cap M]. &&&
\end{array}$$
Since $\varphi:\CC^{\sigma(1)} \rightarrow X_\sigma$ and
$\pi_\sigma:X_{\bar\sigma}\rightarrow X_{ \sigma}$ are good
categorical quotients, the last diagram induces isomorphisms:
$$\CC[{ {\sigma}^\ve}\cap M]\simeq (S(\Sigma)_{x^{\hat\sigma}})^{G(\Sigma)}\simeq S(\sigma)^{G(\sigma)},$$  where the last one is   given by the evaluation
at $x_\rho=1$ for $\rho\not\subseteq\sigma$ since it is the
restriction of the homomorphism
$$\CC[x_\rho\mid\rho\in\Sigma(1)]_{x^{\hat\sigma}}
\simeq\CC\bigl[\bar\sigma^\ve\cap\ZZ^{\Sigma(1)}\bigr]\rightarrow
\CC\bigl[\ZZ_{\ge0}^{\sigma(1)}\bigr]\simeq\CC[x_\rho\mid\rho\in\sigma(1)]
$$
arising from the projection $\ZZ^{\Sigma(1)}\rightarrow
\ZZ^{\sigma(1)}$ dual to the inclusion $\ZZ^{\sigma(1)}
\hookrightarrow\ZZ^{\Sigma(1)}$.

The restriction of the above diagram  to closed subvarieties gives
$$
\begin{array}{cccc}
    & {\bf V}_a(\psi(I))  \longrightarrow &  {\bf V}_a(  I ) \cap X_{\bar\sigma} &\\
   & \downarrow &   \hspace{-2.5cm}\swarrow  &   \\
  &  {\bf V}(I)\cap X_\sigma. &&
\end{array}$$
Equivalently, we get the commutative diagram
$$
\begin{array}{cccc}
  S(\sigma)/ \psi(I) &  \hspace{-0.3cm} \longleftarrow & \hspace{-0.4cm} S(\Sigma)_{x^{\hat\sigma}}/I_{x^{\hat\sigma}}&\\
    \uparrow & \nearrow\quad  &   \\
  \,\, \CC[Y\cap X_\sigma], &&&
\end{array}$$
where $Y={\bf V}(I)\subseteq \xs$, and,   by the proof of
Theorem~\ref{t:quotient}, the vertical morphisms  induce
isomorphisms from $\CC[Y\cap X_\sigma]$ onto the respective rings
of invariants.
\end{pf}

Everything said above can be naturally generalized for
    a toric variety $\xs$ with a torus factor, in which case the
  rays of $\Sigma(1)$ do not span $N_\RR$ and   we can
write $N=N'\oplus N''$ so that $\Sigma(1)$ spans $N'_\RR$. If
$\Sigma'$ denotes the fan consisting of the cones of $\Sigma$ in
$N'_\RR$, then by (5.1.10) in \cite{torvar}, $$\xs\simeq
X_{\Sigma'}\times (\CC^*)^r\simeq (\CC^{\Sigma(1)}\setminus
Z(\Sigma'))/G(\Sigma')\times (\CC^*)^r\simeq
(\CC^{\Sigma(1)+r}\setminus Z'(\Sigma))/G(\Sigma'),$$ where
$r=\dim N''_\RR$ and $Z'(\Sigma)=Z(\Sigma')\times\CC^r\cup
\CC^{\Sigma(1)}\times{\bf V}_a(x_1 \cdots x_r)$. Here,
$G(\Sigma')$ and $Z(\Sigma')$ coincide  with the previously
defined $G(\Sigma)$ and $Z(\Sigma)$ which did not depend on the
condition that $\Sigma(1)$ spans $N_\RR$. We note that
$Z'(\Sigma)={\bf V}_a(B'(\Sigma))\subseteq\CC^{\Sigma(1)+r}$ where
$B'(\Sigma)=\langle \prod_{\rho
 \not\subseteq\sigma} x_\rho\mid\sigma\in\Sigma\rangle$ assuming
 that $x_1,\dots,x_r$ are among $x_\rho$ corresponding to
   ``nonexistent'' rays of $\Sigma$. The variables $x_1,\dots,x_r$ can be thought as coordinates on the torus factor
$(\CC^*)^r$.   The variety $\CC^{\Sigma(1)+r}\setminus
 Z'(\Sigma) $ can be viewed as a toric variety with a surjective toric
 morphism $$\CC^{\Sigma(1)+r}\setminus
 Z'(\Sigma)\rightarrow\xs$$ similar to the construction of $\pi$.
 Define the homogeneous coordinate ring of the toric variety $\xs$
 as
 $$S'(\Sigma)=\CC[x_\rho\mid
 \rho\in\Sigma(1)]\otimes\CC[x_1,\dots,x_r],$$
which is graded by the class group ${\rm Cl}({\Sigma'})$ given by
the short exact sequence $$0\rightarrow M'\rightarrow
\ZZ^{\Sigma(1)}\rightarrow{\rm Cl}(\Sigma')\rightarrow0,$$ where
$m\mapsto(\langle m,v_\rho\rangle)_{\rho\in\Sigma(1)}$ for $m\in
M'=\Hom(N',\ZZ)$. The grading of   $S'(\Sigma)$ is induced by the
grading of the first tensor factor  given by
$$\deg\Biggl(\prod_{\rho\in\Sigma(1) }
x_\rho^{a_\rho}\Biggr)= [(a_\rho)_{\rho\in\Sigma(1)}]\in {\rm
Cl}(\Sigma').$$ Then closed subvarieties of $\xs$ correspond to
the homogeneous ideals of $S'(\Sigma)$ and can be  described by
categorical quotients the same way as in the proof of
Theorem~\ref{t:quotient}.

We will finish this section by working out   toric morphisms
  in terms of
homogeneous coordinates. First, suppose that $X_{\Sigma_2}$ is a
smooth toric variety and    $\phi:X_{\Sigma_1}\rightarrow
X_{\Sigma_2}$ is a toric morphism corresponding to a lattice
homomorphism $\bar\phi:N_1\rightarrow N_2$ compatible with  the
fans $\Sigma_1$ and $\Sigma_2$.   From the categorical quotient
presentations we have additional morphisms
$$\pi_1: X_{\bar\Sigma_1}\rightarrow X_{\Sigma_1},\qquad
\pi_2:X_{\bar\Sigma_2}\rightarrow X_{\Sigma_2}$$ induced by the
natural lattice morphisms $\bar\pi_1:\ZZ^{\Sigma_1(1)}\rightarrow
N_1$ and $\bar\pi_2:\ZZ^{\Sigma_2(1)}\rightarrow N_2$,
respectively.

For every $\rho\in\Sigma_1(1)$, denote by $\sigma_\rho$   the
minimal cone in $\Sigma_2$ containing $\bar\phi(\rho)$. Since
$\Sigma_2$ is smooth, we can write $\bar\phi
(v_{\rho})=\sum_{\xi\in\sigma_{\rho}(1)} a_{\rho,\xi} v_\xi$ with
unique $a_{\rho,\xi}\in\ZZ_{\ge0}$. Define the lattice morphism
$$ {\bar\phi}':\ZZ^{\Sigma_1(1)}\rightarrow\ZZ^{\Sigma_2(1)} \quad
\text{ by }\quad \bar\phi'(e_{\rho})=\sum_{\xi\in\sigma_{\rho}(1)}
a_{\rho,\xi} e_\xi,$$ where $\{e_{\rho}\mid\rho\in\Sigma_1(1)\}$
and $\{e_{\xi}\mid\xi\in\Sigma_2(1)\}$ are the standard bases of
$\ZZ^{\Sigma_1(1)}$ and  $\ZZ^{\Sigma_2(1)}$, respectively. By the
construction,   the diagram
$$ \begin{CD}
  \ZZ^{\Sigma_1(1)} @>\bar\pi_1>> N_1  \\
  @VV{\bar\phi}'V @VV\bar\phi V \\
\ZZ^{\Sigma_2(1)} @>\bar\pi_2>>   N_2
\end{CD}   $$
commutes and induces commutative diagram of toric morphisms
$$ \begin{CD}
\CC^{\Sigma_1(1)} &\,\,\hookleftarrow\,& X_{\bar\Sigma_1} @> \pi_1>> X_{ \Sigma_1}  \\
 @VV\tilde\phi V  @VV{ \phi}'V  @VV \phi V \\
\CC^{\Sigma_2(1)}&\,\,\hookleftarrow\,\,& X_{\bar\Sigma_2} @>
\pi_2>> X_{ \Sigma_2}
\end{CD}   $$
where $\tilde\phi$ is a monomial morphism corresponding to the
ring homomorphism
\begin{equation}\label{e:morphism}{\tilde\phi}^*:\CC[x_\xi\mid
 \xi\in\Sigma_2(1)]\rightarrow\CC[x_\rho\mid
 \rho\in\Sigma_1(1)],\quad
  x_\xi\mapsto \prod_{\rho\in\Sigma_1(1)}x_\rho^{a_{\rho,\xi}}.
\end{equation}
(Here, we assume that $a_{\rho,\xi}=0$ if
$\xi\not\in\sigma_\rho(1)$.) By this formula we see that  if a
point $p=\pi_1(x)\in X_{ \Sigma_1}$ is represented by
$x=(x_\rho)_{\rho\in\Sigma_1(1)}\in\CC^{\Sigma_1(1)}$, then
$\phi(p)$  is represented by
$\tilde\phi(x)=(x_\xi)_{\xi\in\Sigma_2(1)}\in\CC^{\Sigma_2(1)}$,
where $x_\xi=\prod_{\rho\in\Sigma_1(1)}x_\rho^{a_{\rho,\xi}}$ in
homogeneous coordinates.

We can recover the toric morphism $\phi$ from (\ref{e:morphism})
by
$$\CC[{ \sigma_2^\ve}\cap M_2]\simeq (S(\Sigma_2)_{x^{\hat\sigma_2}})^{G(\Sigma_2)} \rightarrow
 (S(\Sigma_1 )_{x^{\hat\sigma_1}})^{G(\Sigma_1)}\simeq \CC[{ \sigma_1^\ve}\cap M_1],$$
where $\sigma_1\in\Sigma_1$, $\sigma_2\in\Sigma_2$, such that
$\bar \phi(\sigma_1)\subseteq \sigma_2$, and the middle
homomorphism is induced by ${\tilde\phi}^*$: $$\prod_{\xi\in
\Sigma_2(1)}x_\xi^{ \langle u,v_\xi\rangle}\mapsto \prod_{\xi\in
\Sigma_2(1)}\prod_{\rho\in\Sigma_1(1)}{x_\rho^{a_{\rho,\xi}\cdot
\langle u,v_\xi\rangle}}=\prod_{\rho\in\Sigma_1(1)}x_\rho^{
\langle
u,\bar\phi(v_\rho)\rangle}=\prod_{\rho\in\Sigma_1(1)}x_\rho^{
\langle \bar\phi^\#( u) ,v_\rho\rangle}$$ for
$u\in\sigma_2^\ve\cap M_2$, where
$\bar\phi^\#:M_2=\Hom(N_2,\ZZ)\rightarrow M_1=\Hom(N_1,\ZZ)$ is
the lattice homomorphism dual to $\bar\phi$.

Note that our formulas for a toric morphism
$\phi:X_{\Sigma_1}\rightarrow X_{\Sigma_2}$ in homogeneous
coordinates work   for any fan $\Sigma_2$ such that the minimal
cones $\sigma_\rho\in\Sigma_2$, containing $\bar\phi(\rho)$ for
$\rho\in\Sigma_1(1)$, give $\bar\phi
(v_{\rho})=\sum_{\xi\in\sigma_{\rho}(1)} a_{\rho,\xi} v_\xi$ with
integer coefficients  $a_{\rho,\xi}\in\ZZ_{\ge0}$. If this
condition fails, the morphisms $\phi'$ and $\tilde\phi$ to make
the commutative diagram do not exist. However, one can still
consider a multi-valued map $\CC^{\Sigma_1(1)}\rightarrow
\CC^{\Sigma_2(1)}$ given by the monomials
$\prod_{\xi\in\Sigma_1(1)}x_\xi^{a_{\xi,\rho}}$ with rational
exponents and show that it descends to a well defined map
$X_{\Sigma_1}\rightarrow X_{\Sigma_2}$ between toric varieties
coinciding with the morphism $\phi$.

\begin{ex} Let $N_1=\ZZ$,   $N_2=\ZZ^2$, $X_{\Sigma_1}=\PP^1$ and
  $X_{\Sigma_2}=\PP(1,1,2)$ be  the weighted projective space
corresponding to the fan $\Sigma_2$ in  $(N_2)_\RR=\RR^2$ with the
primitive lattice generators $ v_1=(1,0)$, $v_2=(-1,-2)$,
$v_3=(0,1)$:
$$ \setlength{\unitlength}{1cm}
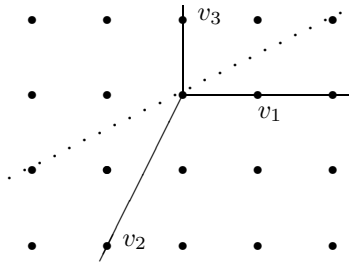
\begin{picture}(8,5)
\put( 3,2.7){$v_{1}$} \put( 1.2,1){$v_{2}$} \put( 2.2,4){$v_{3}$}
\put(0,1){\circle*{0.1}}
\put(0,2){\circle*{0.1}}\put(0,3){\circle*{0.1}}\put(0,4){\circle*{0.1}}
\put(4,1){\circle*{0.1}}
\put(4,2){\circle*{0.1}}\put(4,3){\circle*{0.1}}\put(4,4){\circle*{0.1}}
\put(1,1){\circle*{0.1}} \put(3,3){\circle*{0.1}}
\put(2,2){\circle*{0.1}} \put(2,4){\circle*{0.1}}
\put(3,4){\circle*{0.1}} \put(1,2){\circle*{0.1}}
\put(1,2){\circle*{0.1}}
 \put(1,3){\circle*{0.1}}
\put(1,4){\circle*{0.1}} \put(2,1){\circle*{0.1}}
\put(3,1){\circle*{0.1}} \put(3,2){\circle*{0.1}}
\put(2,3){\circle*{0.1}} \put(2,3){\line(0,1){1.2}}
\put(2,3){\line(1,0){2.2}} \put(2,3){\line(-1,-2){1.1}}
 \multiput(4.1,4.1 )(-0.2,-0.1){23}{\circle*{0.05}}
\end{picture}\label{e:pic}
$$
Consider the lattice homomorphism  $\bar\phi: \ZZ\rightarrow
 \ZZ^2$ given by  $\bar\phi(z)=(2z,z)$, which induces the map of the fan $\Sigma_1$ of $\PP^1$ to the fan
  $\Sigma_2$ (the dotted line in the
picture is the image of $(N_1)_\RR=\RR$) and the corresponding
toric morphism $\phi:\PP^1\rightarrow\PP(1,1,2)$. If $y_1,y_2,y_3$
are the homogeneous coordinates on $\PP(1,1,2)$ corresponding to
$v_1,v_2,v_3$, and $x_1,x_2$ are the homogeneous coordinates on
$\PP^1$ corresponding to the lattice generators $1,-1$ of
$N_1=\ZZ$, then our formulas say that the morphism $\phi$ is given
by $ (x_1,x_2) \mapsto(x_1^2,x_2^2,x_1x_2^3)$, since
$\bar\phi(1)=2v_1+v_3$  and $\bar\phi(-1)= 2v_2+3v_3 $.
\end{ex}

\section{Deformations of affine  toric varieties.}\label{s:1}

In this section, we first  describe Altmann's construction of
deformations of an affine toric variety from \cite{al2} after some
modification and then generalize it using the categorical quotient
presentation of an affine toric variety. While we use an embedding
of the affine toric variety into a higher dimensional affine toric
variety using Minkowski sum decompositions of a polyhedra similar
to \cite{al2}, the essential difference in our construction is
that the image of the affine toric variety is   given by a prime
binomial complete intersection ideal in Cox homogeneous
coordinates, which is a {\it toric ideal}.   Our construction
allows more ways to deform an affine toric variety and to combine
these deformations together than \cite{al2}, and we expect that
 our deformations will span the unobstructed
deformations of an affine toric variety by Kodaira-Spencer map. We
also show that our deformations into certain directions   are
isomorphic to the ``non-toric'' Altmann's deformations in
\cite{al5} and compute the Kodaira-Spencer map for our families.

Let $N$ be a lattice and $M$ be its dual lattice.
 An affine
toric variety  $X_\sigma={\rm Spec}(\CC[{ {\sigma}^\ve}\cap M])$
corresponds to a rational strongly
 convex polyhedral cone  $\sigma\subset N_\RR$ with apex at $0$,
 where $ {\sigma^\ve}$ is the dual cone. We will view elements of the semigroup ring $\CC[{ {\sigma}^\ve}\cap M]$ as finite sums
 $$\sum_{m\in{\sigma}^\ve\cap M}a_m\chi^m,$$ where $a_m\in\CC$ are complex coefficients and  $\chi^m$ are the elements of the ring
 corresponding to lattice points $m$. For $u\in M$ and the cone
 $\sigma$ we will use notation:
$$\sigma(u):=\sigma\cap\{n\in N_\RR\mid\langle u,n\rangle=-1\},$$   the slice of $\sigma$ in the hyperplane where $u$ takes value $-1$.
Take $R\in{\sigma^\ve}\cap M$ and decompose
$\sigma(-R)=Q_0+Q_1+\cdots+Q_k$ into an (ordered) Minkowski sum of
polyhedra with the same recession cone $\sigma\cap R^\perp$ such
that the following two conditions are satisfied:

 $(*)$   the induced decomposition of a vertex of $\sigma(-R)$
 $${\rm vert}(\sigma(-R))={\rm vert}(Q_0)+{\rm vert}(Q_1)+\cdots+{\rm vert}(Q_k)$$  into the sum
 of vertices of $Q_0,Q_1,\dots,Q_k$ has all but possibly one of the summands
  being lattice
 points,

$(**)$  $Q_0\subset\{n\in N_\RR|\langle
 R,n\rangle=1\}$, $Q_i \subset\{n\in N_\RR|\langle
 R,n\rangle=0\}$ for $i=1,\dots,k$.

\begin{rem}
An induced decomposition of a vertex of $\sigma(-R)$ is unique. If
$\sigma(-R) $ is a lattice polyhedron, then  condition $(*)$ means
 that all $Q_0,Q_1,\dots,Q_k$ are lattice polyhedrons as well.
The recession cone of a polyhedron $Q$ is  defined as ${\rm
 rec}(Q)=\{n\in N_\RR|\,q+tn\in Q \,\,\forall\,\, q\in Q,\,t\ge0\}$.
\end{rem}

 Consider the lattice $\tilde{N}=N\oplus\ZZ^k$ denoting by  $\{e_1,\dots,e_k\}$ the standard basis for the second summand
 $\ZZ^k$. The dual lattice is then $\tilde{M}=M\oplus\ZZ^k$. Let
 $\{e_1^*,\dots,e_k^*\}$ be the dual to $\{e_1,\dots,e_k\}$ basis for the second
 component.
Let the cone \begin{equation}\label{e:cons}\tilde{\sigma}=\langle
\sigma,
 Q_0-e_1-\dots-e_k,Q_1+e_1,\dots,Q_k+e_k\rangle\end{equation} in
 $\tilde{N}_\RR$ be generated by the indicated sets.

The inclusion of cones $\sigma\subset\tilde{\sigma}$  induces a
ring homomorphism $$\CC[{ {\tilde{\sigma}}^{\ve}}\cap
\tilde{M}]\rightarrow\CC[{\sigma}^{\ve}\cap M]$$ corresponding to
a morphism of toric varieties $X_\sigma\rightarrow
X_{\tilde{\sigma}}=\spe (\CC[{ {\tilde{\sigma}}^{\ve}}\cap
\tilde{M}] )$.

\begin{lem}\label{l:1} The ring
homomorphism $\CC[{ {\tilde{\sigma}}^{\ve}}\cap
\tilde{M}]\rightarrow\CC[{\sigma}^{\ve}\cap M]$ is surjective.
\end{lem}

\begin{pf}   This ring
homomorphism  is induced by the natural projection of lattices
$\tilde{M}=M\oplus\ZZ^k\rightarrow M$ and it is sufficient to show
that the map of the semigroups ${\tilde{\sigma}}^\ve\cap
\tilde{M}\rightarrow{\sigma}^{\ve}\cap M$ is surjective. Let $s\in
{\sigma}^{\ve}\cap M$ and denote by $\min\langle s,Q\rangle$ the
minimal value of the vector $s$ on the polyhedron $Q=\sigma(-R)$,
which occurs at one of the vertices of $Q$ since the value of $s$
on the points  of $Q\subset\sigma$ is bounded from below by $0$.
Since $Q_i$ is a Minkowski summand of $Q=\sigma(-R)$, the value
$\min\langle s,Q_i\rangle$ is also attained at a vertex of $Q_i$.
Denote by $[[\min\langle s,Q_i\rangle]]$
 the greatest integer value of the corresponding real number.
 Then we claim that the lattice point $$\tilde{s}=s-\sum_{i=1}^k[[\min\langle
 s,Q_i\rangle]]e_i^*$$
 is in ${\tilde{\sigma}}^\ve\cap
\tilde{M}$, whence the statement follows. We can check this claim
on the generators of $\tilde\sigma$. If $n\in\sigma$, then
$\langle \tilde s,n\rangle=\langle   s,n\rangle\ge0$ since
$s\in\sigma^{\ve}$.  If $q_i\in Q_i$ for $i\ge1$, then $\langle
\tilde s, q_i+e_i\rangle=\langle   s,q_i\rangle-[[\min\langle
 s,Q_i\rangle]]\ge0$. If $ q_0\in Q_0$, then $$\langle
\tilde s, q_0-\sum_{i=1}^k e_i\rangle=\langle
s,q_0\rangle+\sum_{i=1}^k[[\min\langle
 s,Q_i\rangle]]\ge \sum_{i=0}^k[[\min\langle
 s,Q_i\rangle]]=[[\min\langle
 s,Q\rangle]]\ge0. $$
The last equality is true because all but possibly one of
$\min\langle
 s,Q_i\rangle$  are integers by the condition $(*)$ on the Minkowski sum
 decomposition.
\end{pf}

\begin{rem}  An analog of this lemma in a complicated setting was proved in (4.1.1) in \cite[p.~162]{al2} but
    we found  this elegant  proof of it. Moreover, our proof shows
    that a more general result holds: $\sigma(-R)$ can be
    replaced by any rational polyhedron $Q$ inside $\sigma$ with a Minkowski sum decomposition $Q=Q_0+Q_1+\cdots+Q_k$
    satisfying only condition
    $(*)$.
\end{rem}

  The above lemma tells us that we have an
embedding of the  toric varieties $X_\sigma\hookrightarrow
X_{\tilde{\sigma}}$. The next one shows the defining equations.

\begin{pr}\label{p:comp} Let $R\in {\sigma}^{\ve}\cap M$ and $\tilde{\sigma}$ be as in $(\ref{e:cons})$.
Then the ideal  $\ker (\CC[{ {\tilde{\sigma}}^{\ve}}\cap
\tilde{M}]\rightarrow\CC[{\sigma}^{\ve}\cap M])$  is a complete
intersection generated by the regular sequence
$\chi^{R+e_1^*}-\chi^R,\dots,\chi^{R+e_k^*}-\chi^R$.
\end{pr}

\begin{pf} While our setting is slightly different from that of
Altmann (up to a linear transformation our cone $\tilde\sigma$
coincides with the one in (3.2) of \cite{al3}), the proof of this
result is the same as in \cite[p.~162-163]{al2} and we will
generalize this result in Proposition~\ref{p:iso}. The regularity
of the sequence follows from the fact that this is a complete
intersection given by ``homogeneous'' elements, where the
$\NN$-grading on $\CC[{ {\tilde{\sigma}}^{\ve}}\cap \tilde{M}]$
can be induced by the paring
$\langle{\underline{\,\,\,}},q\rangle:\tilde M\rightarrow\ZZ$ for
some $q\in \RR_{>0}\cdot Q\cap N$.
\end{pf}

\begin{rem}\label{r:2ndcon} Condition $(**)$ on the Minkowski sum decomposition was
 only necessary  to fix the defining equations of $X_\sigma$ inside
 $X_{\tilde{\sigma}}$. It was used to present the construction in \cite{al5}. If   this condition is dropped,
 one can show
 that $R$ must be replaced by $\tilde R=  R-\sum_{i=1}^k \langle
 R,Q_i\rangle e_i^*$ in the equations. Here,  $\langle
 R,Q_i\rangle$ denotes the constant value of $R$ on $Q_i$, which must be an integer from the condition $(*)$ and
 the equation   $\sum_{i=0}^k  \langle R,Q_i\rangle=\langle
 R,\sigma(-R)\rangle=1$.  Note that in this case, for $i\ne0$,  $\langle
 \tilde R   ,Q_i+e_i\rangle=0$, and
   $\langle\tilde R   ,Q_0-e_1-\dots-e_k\rangle=\sum_{i=0}^k  \langle
   R,Q_i\rangle=1$.
 These are the same equalities if we left condition $(**)$ and  used $R$
instead of $\tilde R   $.
 \end{rem}

Thus, we have our toric variety $X_\sigma$ embedded into another
toric variety  $X_{\tilde{\sigma}}$ as a complete intersection
given by $\chi^{R+e_i^*}-\chi^R$ for $i=1,\dots,k$. A
$k$-parameter embedded deformation of $X_\sigma$ corresponds to
the following natural diagram:
\begin{equation}\label{diagr}
\begin{array}{ccccc}
X_\sigma & \subset & {  \mathcal   X} & \subset & X_{\tilde{\sigma}}\times\CC^k\\
\downarrow & & \downarrow& \swarrow  &\\
\{0\}&\subset &\CC^k,&&
\end{array}
\end{equation}
where the deformation family $ {  \mathcal   X}$ in
$X_{\tilde{\sigma}}\times\CC^k$ is given by
$\chi^{R+e_i^*}-\chi^R-\lambda_i=0$, for $i=1,\dots,k$, with
$\lambda_i$ being coordinates on $\CC^k$.

\begin{rem} In \cite{al2,al5},    the deformation   of  $X_{ {\sigma}}$
was described   differently   by the flat morphism
$X_{\tilde{\sigma}}\rightarrow\CC^k$ given by the values of the
functions $(\chi^{R+e_1^*}-\chi^R,\dots,\chi^{R+e_k^*}-\chi^R)$
with the special fiber over the origin isomorphic to $X_\sigma$.
But the total space $X_{\tilde{\sigma}}$ of this family is
naturally isomorphic to the above $ {  \mathcal   X}$.
\end{rem}

By using homogeneous coordinates on the affine toric variety we
found the following generalization of  Altmann's deformations of
$X_\sigma$.  As we described the above construction, notice that
there is no problem to proceed with the construction of
$\tilde{\sigma}$ even if $R\not\in{\sigma^\ve}\cap M$.  We still
get the embedding $X_\sigma\hookrightarrow X_{\tilde{\sigma}}$.
But the image may no longer be a complete intersection in terms of
affine functions on $X_{\tilde{\sigma}}$  if    $X_\sigma$ is
singular. We also note that
$\chi^R,\chi^{R+e_1^*},\dots,\chi^{R+e_k^*}$ will not be functions
on $X_{\tilde{\sigma}}$ if $R\not\in{\sigma^\ve}\cap M$.

 From now on, we will use a more general
construction  of $\tilde\sigma =\langle \sigma,
 Q_0-\sum_{i=1}^k e_i ,Q_1+e_1,\dots,Q_k+e_k\rangle$ by replacing
$\sigma(-R) $ with any rational polyhedron $Q$ in $\sigma$  such
that $0\notin Q$ and the Minkowski sum decomposition
 $Q=Q_0+\cdots +Q_k$ into polyhedra satisfies condition $(*)$. We have the
following property for $\tilde\sigma$:

\begin{lem}\label{l:coneprop}  Let $\sigma$ be a strongly convex rational polyhedral cone in $N_\RR$
 and $Q$ be a rational polyhedra contained in $\sigma$. Let
  $Q=Q_0+\cdots +Q_k$ be a Minkowski sum decomposition into
  polyhedra and $$\tilde\sigma =\langle \sigma,
 Q_0-\sum_{i=1}^k e_i ,Q_1+e_1,\dots,Q_k+e_k\rangle.$$ Then
$\tilde\sigma\cap N_\RR=\sigma$ and, if $0\notin Q$, then
$\tilde\sigma$ is strongly convex and its one-dimensional faces
are among the rays through the vertices of $Q_0-\sum_{i=1}^k e_i
,Q_1+e_1,\dots,Q_k+e_k$ and the edges of $\sigma$ not lying in the
cone  $\langle Q_0 -\sum_{i=1}^k e_i ,Q_1 +e_1,\dots,Q_k
+e_k\rangle$.
\end{lem}

\begin{pf} To show that $\tilde\sigma$ is strongly convex take $u$ in the interior of $\sigma^\ve$ such that
$\min\langle u,Q\rangle>0$ and check that  $ u-\sum_{i=1}^k
\min\langle u,Q_i\rangle e_i^* +\frac{\min\langle
u,Q\rangle}{k+1}\sum_{i=1}^k
 e_i^*$ is positive on all generators of $\tilde\sigma$. The rest
 of the
statements are elementary to prove and left as an exercise to the
reader.\end{pf}

The crucial idea to generalize Altmann's construction is to use
the categorical quotient presentation of a toric variety from
\cite{c}. For simplicity we will assume that $X_\sigma$ does not
have torus factors (equivalently, $\dim \sigma=\dim N_\RR$), but
this assumption will  not be essential as described in
Section~\ref{s:toric}. If $\dim \sigma$ is maximal, then $\dim
\tilde\sigma=\dim \tilde N_\RR$ and $X_{\tilde{\sigma}}$ will be
without torus factors as well and for the affine toric variety
$X_{\tilde\sigma}$ we have presentation
$${\rm
Spec}(\CC[\tilde\sigma^\ve\cap \tilde{M}]) \simeq
\CC^{\tilde\sigma(1)}/ G(\tilde\sigma)={\rm
Spec}\bigl(\CC[x_\xi\mid\xi\in\tilde\sigma(1)]^{G(\tilde\sigma)}\bigr),$$
where $x_\xi$ are the    homogeneous coordinates of
$X_{\tilde{\sigma}}$.
 On
the level of rings this corresponds to the isomorphism
\begin{equation}\label{e:isom}\CC[{ {\tilde{\sigma}}^{\ve}}\cap
\tilde{M}]\simeq
\CC[x_\xi\mid\xi\in\tilde\sigma(1)]^{G(\tilde\sigma)},\end{equation}
which sends $\chi^u$ to $\prod_{\xi\in\tilde\sigma(1)}x_\xi^{
\langle u,v_\xi\rangle}$, where $v_\xi$ are  the primitive lattice
generators of the edges $\xi$ of $\tilde{\sigma}$. We will use
notation $x^u:=\prod_{\xi\in\tilde\sigma(1)}x_\xi^{ \langle
u,v_\xi\rangle}$ for any $u\in \tilde{M}$.

\begin{pr}\label{p:iso} Let $\tilde\sigma$ be as in Lemma~\ref{l:coneprop}, corresponding to $0\not\in Q\subset\sigma$
and decomposition $Q=Q_0+\cdots +Q_k$ satisfying condition $(*)$.
Then there is a natural ring isomorphism
$$\CC[{\sigma}^{\ve}\cap M]\simeq
(\CC[x_\xi\mid\xi\in\tilde\sigma(1)]/I)^{G(\tilde\sigma)},$$ where
the ideal $I$ is generated by  binomials
$$ \prod_{ \langle e_i^*,v_\xi\rangle>0 } x_\xi^{ \langle
e_i^*,v_\xi\rangle}-\prod_{ \langle e_i^*,v_\xi\rangle<0 } x_\xi^{
-\langle e_i^*,v_\xi\rangle},$$ for $i=1,\dots,k$, forming a
regular sequence in $\CC[x_\xi\mid\xi\in\tilde\sigma(1)]$.
Moreover, the ideal $I$ is   prime.
\end{pr}

\begin{pf} Let $I_1$ denote  the $\ker (\CC[{ {\tilde{\sigma}}^{\ve}}\cap
\tilde{M}]\rightarrow\CC[{\sigma}^{\ve}\cap M])$. Then
$$\CC[{\sigma}^{\ve}\cap M]\simeq \CC[{ {\tilde{\sigma}}^{\ve}}\cap
\tilde{M}]/I_1.$$ It is easy to find that the ideal $I_1$ is
generated by binomials:
$$I_1=\langle\chi^r-\chi^s\mid r,s\in { {\tilde{\sigma}}^{\ve}}\cap\tilde{M},r-s\in {\rm ker}(\tilde{M}\rightarrow M)\rangle.$$
Consider the ring homomorphism $\phi:
\CC[x_\xi\mid\xi\in\tilde\sigma(1)]^{G(\tilde{\sigma})}\rightarrow(\CC[x_\xi\mid\xi\in\tilde\sigma(1)]/I)^{G(\tilde{\sigma})},$
assigning $f\mapsto f+I$. To prove the ring isomorphism in the
first part of the proposition, it suffices to show that
$$ {\rm ker}(\phi)=\langle x^r-x^s\mid r,s\in {
{\tilde{\sigma}}^{\ve}}\cap\tilde{M},r-s\in {\rm
ker}(\tilde{M}\rightarrow M)\rangle.$$

 Now we use some ideas from \cite[p.~162-163]{al2}.
 Suppose $r,s\in {
{\tilde{\sigma}}^{\ve}}\cap\tilde{M}$ and $r-s\in {\rm
ker}(\tilde{M}\rightarrow M)$, then we can write $\ds
r-s=\sum_{i=1}^k\alpha_i^+ e_i^*-\sum_{i=1}^k\alpha_i^- e_i^*$
with $\alpha_i^+\cdot\alpha_i^-=0$ and $\alpha_i^+,\alpha_i^-\in
\ZZ_{\ge0}$ for all $i$.   Denote
$$q=r-\sum_{i=1}^k\alpha_i^+ e_i^*=s-\sum_{i=1}^k\alpha_i^-
e_i^*,$$ which is in ${ {\tilde{\sigma}}^{\ve}}\cap \tilde{M}$ by
an easy check restricting one of the presentations of $q$ to the
generators of the cone $\tilde{\sigma}$. We have
 $$x^r-x^s= x^{q+\sum_{i=1}^k\alpha_i^+ e_i^*}-x^{q+\sum_{i=1}^k\alpha_i^-e_i^*}.$$
Since $q+\sum_{i=1}^k\alpha_i^+ e_i^* =r\in {
{\tilde{\sigma}}^{\ve}}\cap \tilde{M}$ and $q\in{
{\tilde{\sigma}}^{\ve}}\cap \tilde{M}$, the monomial $
x^{q+\sum_{i=1}^k\alpha_i^+ e_i^*}$ is divisible by
$$\prod_{i=1}^k\Biggl( \prod_{ \langle e_i^*,v_\xi\rangle>0 } x_\xi^{ \langle
e_i^*,v_\xi\rangle}\Biggr)^{\alpha_i^+}.$$

Using relations
$$\ds\prod_{ \langle e_i^*,v_\xi\rangle>0 } x_\xi^{ \langle
e_i^*,v_\xi\rangle}-\prod_{ \langle e_i^*,v_\xi\rangle<0 } x_\xi^{
-\langle e_i^*,v_\xi\rangle}=(1-x^{-e_i^*})\prod_{ \langle
e_i^*,v_\xi\rangle>0 } x_\xi^{ \langle e_i^*,v_\xi\rangle}=0 $$
modulo $I$, we get $x^{q+\sum_{i=1}^k\alpha_i^+ e_i^*}-x^q\in I$.
Similarly, $x^{q+\sum_{i=1}^k\alpha_i^- e_i^*}-x^q\in I$. Hence
$x^r-x^s\in I$, showing that these binomials lie in ${\rm
ker}(\phi)$.
  The
other inclusion can be shown straightforward.

To show the second part of the proposition, we will present  the
ideal $I$ as a lattice basis ideal (see  Appendix~\ref{s:ap}) and
show that it is a complete intersection toric ideal. Consider the
 lattice $L=\bigoplus_{i=1}^k \ZZ e_i^*$ and embed it into
 $\ZZ^{\tilde\sigma(1)}$ by the homomorphism
 $$L=\bigoplus_{i=1}^k \ZZ e_i^*\longrightarrow\ZZ^{\tilde\sigma(1)},
  \quad l\mapsto (\langle l,v_\xi\rangle)_{\xi\in\tilde\sigma(1)},$$
where the pairing $\langle\underline{\,\,}\,,\underline{\,\,}\,
\rangle$ is induced from the pairing of $\tilde M=M\bigoplus
\bigl(\bigoplus_{i=1}^k \ZZ e_i^*\bigr)$ and $\tilde N$. Note that
for a vertex of $q$ of $Q$, the Minkowski sum decomposition
$Q=Q_0+Q_1+\cdots+Q_k$ induces the decomposition
$q=q_0+q_1+\cdots+q_k$, where $k$ of the summands are lattice
points. In particular, the primitive lattice generators $v_\xi$ of
the edges of $\tilde\sigma$ will contain $k$ of the points
$q_0-\sum_{i=1}^k e_i, q_1+e_1,\dots, q_k+e_k$. Since  any $k$ of
the  points $ -\sum_{i=1}^k e_i,  e_1,\dots,  e_k$ form a lattice
basis of the lattice $\bigoplus_{i=1}^k \ZZ e_i$ dual to
$\bigoplus_{i=1}^k \ZZ e_i^*$, the lattice point $l\in
 L=\bigoplus_{i=1}^k \ZZ e_i^*$ is completely determined by the
 integral
 values $\langle l,v_\xi\rangle$ for $\xi\in\tilde\sigma(1)$.
 Therefore,  the lattice homomorphism  is injective.

Now, our ideal $I$ equals the lattice basis ideal $I_{\bf L}$ (see
 Appendix~\ref{s:ap}) corresponding to the basis $(\langle
e_i^*,v_\xi\rangle)_{\xi\in\tilde\sigma(1)}$, for $i=1,\dots,k$,
of the homomorphic image of the lattice $L$ in
$\ZZ^{\tilde\sigma(1)}$.
 Then
note that the lattice $L$ is saturated in $\ZZ^{\tilde\sigma(1)}$
(equivalently, the quotient $\ZZ^{\tilde\sigma(1)}/L$ is torsion
free), since $(\langle
l,v_\xi\rangle/s)_{\xi\in\tilde\sigma(1)}\in \ZZ^{\tilde\sigma(1)}
$ for a positive integer $s$ and $l\in \bigoplus_{i=1}^k \ZZ
e_i^*$ implies that $l/s\in \bigoplus_{i=1}^k \ZZ e_i^*$ by
arguing  as above  that the integral values of $l/s$ on $k$ of the
points $ -\sum_{i=1}^k e_i,  e_1,\dots,  e_k$ completely determine
the lattice point in $\bigoplus_{i=1}^k \ZZ e_i^*$. Hence, the
lattice ideal $$I_L=\langle {\bf x}^u-{\bf x}^v\mid u,v\in
\NN^{\tilde\sigma(1)} \text{ and } u-v\in L\rangle$$ is a toric
ideal (i.e., a prime binomial ideal) by Theorem~7.4 in \cite{ms}
(also, see Appendix~\ref{s:ap}).

To show that the ideal $I=I_{\bf L}$ coincides with the  prime
ideal $I_L$ and that the binomial sequence is regular, it suffices
to show that the $k\times |\tilde\sigma(1)|$ matrix $M=(\langle
e_i^*,v_\xi\rangle)_{1\le i\le k,\xi\in\tilde\sigma(1)}$ is mixed
dominating by Corollary~\ref{c:toricci} and Remark~\ref{r:satur}.
But this sign pattern of the matrix entries (every row has
positive and negative entries, and every square submatrix does not
have the same property) we can easily see if we order the edges
$\xi$ of $\tilde\sigma$ according to the sets $Q_0-\sum_{i=1}^k
e_i ,Q_1+e_1,\dots,Q_k+e_k$ and the edges of $\sigma$ not lying in
the cone  $\langle Q_0 -\sum_{i=1}^k e_i ,Q_1 +e_1,\dots,Q_k
+e_k\rangle$ (see Lemma~\ref{l:coneprop}):

$$M=\left[\begin{array}{rrrrrrrrrrr}- & +& 0&0&\cdots&0&0&\cdots&0\\ -& 0&
+&0& \cdots&0&0&\cdots&0\\-& 0& 0&+ &\cdots&0 &0&\cdots&0\\
\vdots&\vdots&\vdots&\vdots&&\vdots&\vdots&&\vdots\\-& 0&
0&0&\cdots&+&0&
  \cdots&0
\end{array}\right],$$ where $-$ and $+$  actually represent
several consecutive negative or positive entries, respectively.
\end{pf}

The cone $\tilde\sigma$ is strongly convex, whence $\dim
X_{\tilde{\sigma}}=\dim X_\sigma+k$. By our definition of the
complete intersection in a toric variety and its categorical
quotient presentation in Section~\ref{s:toric}, the last
proposition means

\begin{thm}\label{c:eq}  Affine toric variety $X_\sigma$  is embedded into
affine toric variety  $X_{\tilde{\sigma}}$ as a complete
intersection given by a regular sequence $$\ds\prod_{ \langle
e_i^*,v_\xi\rangle>0 } x_\xi^{ \langle e_i^*,v_\xi\rangle}-\prod_{
\langle e_i^*,v_\xi\rangle<0 } x_\xi^{ -\langle
e_i^*,v_\xi\rangle}=0$$ for $i=1,\dots,k$.
\end{thm}

\begin{rem} \label{r:gener}  To see that Proposisition~\ref{p:iso} generalizes
Proposition~\ref{p:comp}, assume that
$Q=\sigma(-R)=Q_0+\cdots+Q_k$ satisfies $(*)$ and $(**)$ and
notice the equality
$$
\prod_{ \langle e_i^*,v_\xi\rangle>0 } x_\xi^{ \langle
e_i^*,v_\xi\rangle}-\prod_{ \langle e_i^*,v_\xi\rangle<0 } x_\xi^{
-\langle e_i^*,v_\xi\rangle}=\ds(x^{R+e_i^*}-x^R) \prod_{\langle
R,v_\xi\rangle<0 } x_\xi^{-\langle R,v_\xi\rangle}$$ since
$-e_i^*$ and $R$ have the same value 1 on $Q_0 -\sum_{i=1}^k e_i$.
 If $R\in{\sigma}^{\ve}\cap
M$, then the above binomial reduces to $x^{R+e_i^*}-x^R$ which
corresponds to $\chi^{R+e_i^*}-\chi^R$ by the isomorphism
(\ref{e:isom}). The exponents $\langle e_i^*,v_\xi\rangle$ can be
viewed as the ``distances'' of the lattice points $v_\xi$  along
$e_i$ from the linear subspace $N_\RR$ inside $\tilde N_\RR$.
  If $R$ is not in ${\sigma}^{\ve}\cap M$, the
ideal $I_1$ defining $X_\sigma$ in $X_{\tilde{\sigma}}$ in the
proof of Proposisition~\ref{p:iso} need not be a complete
intersection in terms of affine coordinate funcitons, but using
homogeneous coordinates we can still represent our toric variety
$X_\sigma$ as a complete intersection in $X_{\tilde{\sigma}}$,
providing  us more opportunities to deform the affine toric
variety.
\end{rem}

Denote by $\sigma(u)^c$ the  convex hull of vertices of a possibly
unbounded polyhedra $\sigma(u)=\sigma\cap\{n\in N_\RR\mid\langle
u,n\rangle=-1\}$. For any $u\in M$,
  a polyhedron $Q\subset\sigma$ such
that $$0\notin Q, \quad \min\langle u,Q\rangle\ge-1,\quad
\sigma(u)^c\subset \RR_{>0} \cdot Q$$ and a Minkowski sum
decomposition $Q=Q_0+ \cdots+Q_k$ with the condition $(*)$ we get
the $k$-parameter embedded deformation of $X_\sigma$ in
$X_{\tilde{\sigma}}$ corresponding to the same diagram as
 in (\ref{diagr}) but given by the equations in homogeneous
 coordinates:
\begin{equation}\label{e:deformh} \prod_{ \langle e_i^*,v_\xi\rangle>0 }
x_\xi^{ \langle e_i^*,v_\xi\rangle}-\prod_{ \langle
e_i^*,v_\xi\rangle<0 } x_\xi^{ -\langle
e_i^*,v_\xi\rangle}-\lambda_i x^{\tilde u}\prod_{ \langle
e_i^*,v_\xi\rangle<0 } x_\xi^{ -\langle
e_i^*,v_\xi\rangle}=0\end{equation} for $i=1,\dots,k$, where
$\tilde u=u-\sum_{i=1}^k[[\min\langle u,Q_i\rangle]] e_i^*$. We
will denote the total space of this deformation family by  $ {
\mathcal X}_{u Q Q_0\dots Q_k}$  because it depends on the lattice
point $u\in M$ and the Minkowski sum decomposition of $Q$.

\begin{rem}\label{r:mon}  By Lemma~\ref{l:coneprop},
notice   $\langle e_i^*,v_\xi\rangle<0 $ only for $v_\xi$, which
are multiples of vertices of $Q_0-\sum_{i=1}^k e_i$. Moreover,
$$\min\Bigl\langle\tilde u-e_i^*,Q_0-\sum_{i=1}^k e_i\Bigr\rangle\ge1+\sum_{i=0}^k[[\min\langle
 u,Q_i\rangle]]=1+ [[\min\langle
 u,\sum_{i=0}^k Q_i \rangle]]\ge0,$$
  $\min\langle\tilde u, Q_i + e_i\rangle=\min\langle
u,Q_i\rangle
 -[[\min\langle u,Q_i\rangle]]\ge0$, as in the proof of
of Lemma~\ref{l:1}, and $  \langle \tilde u,v_\xi\rangle=\langle
  u,v_\xi\rangle\ge0$ for every edge generator $v_\xi$ of $\tilde\sigma$, which is in $\sigma=\tilde \sigma\cap N_\RR$,
  since $\sigma(u)^c\subset \RR_{\ge0} \cdot Q$. These conditions
  guarantee that $x^{\tilde u}\prod_{
\langle e_i^*,v_\xi\rangle<0 } x_\xi^{ -\langle
e_i^*,v_\xi\rangle}$ is a  monomial in
$\CC[x_\xi\mid\xi\in\tilde\sigma(1)]$ and we get a well defined
 deformation family. We also note that in most cases there is an infinite number
 of such $u$'s.
\end{rem}

\begin{rem} When  $Q=\sigma(-R)=Q_0+\cdots+Q_k$ for some $R\in M$, the deformation family
$ { \mathcal X}_{-R\sigma(-R)Q_0\dots Q_k}$ is given  by
polynomials
\begin{multline*} \prod_{ \langle
e_i^*,v_\xi\rangle>0 } x_\xi^{ \langle e_i^*,v_\xi\rangle}-\prod_{
\langle e_i^*,v_\xi\rangle<0 } x_\xi^{ -\langle
e_i^*,v_\xi\rangle}-\lambda_i x^{-\tilde R}\prod_{ \langle
e_i^*,v_\xi\rangle<0 } x_\xi^{ -\langle e_i^*,v_\xi\rangle}=\\=
(x^{e_i^*}-1-\lambda_i x^{-\tilde R})\prod_{ \langle
e_i^*,v_\xi\rangle<0 } x_\xi^{ -\langle e_i^*,v_\xi\rangle}=
(x^{\tilde R+e_i^*}-x^{\tilde R}-\lambda_i) \prod_{\langle \tilde
R,v_\xi\rangle<0 } x_\xi^{-\langle \tilde
R,v_\xi\rangle}.\end{multline*} If $R\in{\sigma}^{\ve}\cap M$ and
$\sigma(-R)=Q_0+\cdots+Q_k$ satisfies condition $(**)$, these
polynomials reduce to $x^{  R+e_i^*}-x^{  R}-\lambda_i$ and by the
isomorphism (\ref{e:isom}) correspond to $\chi^{ R+e_i^*}-\chi^{
R}-\lambda_i$ as in the deformation family $(\ref{diagr})$
isomorphic to Atltmann's construction in \cite{al2}.
\end{rem}

  We  now  relate our construction of deformations of
an affine toric variety $X_\sigma$ to  the ``non-toric''
construction in (3.5) of \cite{al5} for an arbitrary $R\in M$.
After some modification Altmann's construction can be described as
follows. Consider cone $\tau=\sigma\cap\{n\in N_\RR\mid\langle
R,n\rangle\ge0\}$. Corresponding to a Minkowski sum decomposition
$\sigma(-R)=Q_0+Q_1+\dots+Q_k$ satisfying conditions $(*)$ and
$(**)$, we have the cone $\tilde\sigma$ from $(\ref{e:cons})$ and
$\tilde\tau=\langle\tau,
Q_0-e_1-\dots-e_k,Q_1+e_1,\dots,Q_k+e_k\rangle$. Let
 $$B:=\CC[{ {\tilde{\sigma}}^{\ve}}\cap
\tilde{M}][\chi^{R+e_1^*}-\chi^R,\dots,\chi^{R+e_k^*}-\chi^R]\subseteq
\CC[{ {\tilde{\tau}}^{\ve}}\cap \tilde{M}].$$ Then projection
$\tilde M\rightarrow M$ induces a surjective ring homomorphism
$B\rightarrow \CC[{ {{\sigma}}^{\ve}}\cap {M}]$ giving an
embedding $X_{\sigma}\hookrightarrow \spe\, B$. Moreover, Altmann
showed that the sequence of functions
$\chi^{R+e_1^*}-\chi^R,\dots,\chi^{R+e_k^*}-\chi^R$ is regular in
the ring $B$ giving a flat morphism $\spe\, B\rightarrow \CC^k$
with the fiber over the origin isomorphic to $X_{\sigma}$.

\begin{pr} Let $\tilde \sigma$ be as in $(\ref{e:cons})$. There is a natural ring isomorphism
$$(\CC[x_\xi\mid\xi\in\tilde\sigma(1)]\otimes\CC[\lambda_1,\dots,\lambda_k]/I_{\lambda})^{G(\tilde\sigma)}\simeq B,$$
 where the ideal $I_{\lambda}$ is generated by
$$(x^{
R+e_i^*}-x^{  R}-\lambda_i) \prod_{\langle  R,v_\xi\rangle<0 }
x_\xi^{-\langle   R,v_\xi\rangle}$$ for $i=1,\dots,k$.
\end{pr}

\begin{pf} Consider the ring homomorphism
$$
(\CC[x_\xi\mid\xi\in\tilde\sigma(1)]\otimes\CC[\lambda_1,\dots,\lambda_k]/I_{\lambda})^{G(\tilde\sigma)}
\rightarrow B$$ assigning the generators $
x^u+I_{\lambda}\mapsto\chi^u$ for $u\in\tilde\sigma^\ve$ and $
\lambda_i+I_{\lambda} \mapsto \chi^{R+e_i^*}-\chi^R$ for
$i=1,\dots,k$. It is well defined since
${G(\tilde\sigma)}$-invariant elements of $I_{\lambda}$ are the
linear combinations of
 $\ds(x^{R+e_i^*}-x^R-\lambda_i)x^u$
such that $u,u+R\in { {\tilde{\sigma}}^{\ve}}\cap\tilde{M}$ by
definition of $I_\lambda$ and Remark~\ref{r:gener}.
 The ring homomorphism  is clearly surjective. To show
injectivity we can use the fact from (3.6) in \cite{al5} that
elements of $B$ can be  uniquely written in the form
$$\sum_{(n_1,\dots,n_k)\in\NN^k}c_{n_1,\dots,n_k}\cdot(\chi^{R+e_1^*}-\chi^R)^{n_1}\cdots(\chi^{R+e_k^*}-\chi^R)^{n_k}$$
with $c_{n_1,\dots,n_k}\in\CC[{
{\tilde{\sigma}}^{\ve}}\cap\tilde{M}]$ such that  its every
monomial term $\chi^s$ has $s-(R+e_i^*)\notin\tilde\sigma^\ve$ for
$i=1,\dots,k$. Similarly, every element of
$(\CC[x_\xi\mid\xi\in\tilde\sigma(1)]\otimes\CC[\lambda_1,\dots,\lambda_k])^{G(\tilde\sigma)}$
modulo the ideal $I_\lambda$ can be written in the form
$$\sum_{(n_1,\dots,n_k)\in\NN^k}c'_{n_1,\dots,n_k}\cdot  \lambda_1 ^{n_1}\cdots \lambda_k^{n_k}$$
with
$c'_{n_1,\dots,n_k}\in\CC[x_\xi\mid\xi\in\tilde\sigma(1)]^{G(\tilde\sigma)}$
such that its every monomial $x^s$ has
$s-(R+e_i^*)\notin\tilde\sigma^\ve$ for $i=1,\dots,k$. The last
claim can be shown as in (3.6) of \cite{al5}: if for $x^s\in\CC[{
{\tilde{\sigma}}^{\ve}}\cap\tilde{M}]$ we have
$s-(R+e_i^*)\in\tilde\sigma^\ve$ for some $i$, then write
$$x^s+I_\lambda=x^{s-e_i^*}+x^{s-(R+e_i^*)}(x^{R+e_i^*}-x^R)+I_\lambda=x^{s-e_i^*}+x^{s-(R+e_i^*)}\lambda_i+I_\lambda,$$
where $s-e_i^*\in\tilde\sigma^\ve$ and  repeat this step for
$x^{s-e_i^*}$ and $x^{s-(R+e_i^*)}$.
 From the uniqueness of the first sum representation of the elements
 of $B$ we get the injectivity.
\end{pf}

The last proposition means that we have a natural isomorphism
between our deformation   and Altmann's ``non-toric'' construction
in \cite{al5}:
$$
\begin{array}{ccccc}
 { \mathcal X}_{-R\sigma(-R)Q_0\dots Q_k} & \simeq & \spe\, B \\
\downarrow & & \downarrow&  \\
\CC^k&=&\CC^k.&&
\end{array}
$$

 Next we compute the
Kodaira-Spencer map for  our deformation families. The  space
$T^1_{X_\sigma}$ of infinitesimal deformations of $X_\sigma$ was
calculated by K. Altmann in \cite{al1}. It is $M$-graded:
$$T^1_{X_\sigma}=\bigoplus_{R\in M}T^1_{X_\sigma}(-R),$$
where the graded pieces can be described as follows. Let
$a_1,\dots,a_N$ be the primitive lattice generators of the edges
of $\sigma$ and $S=\{s_1,\dots,s_w\}$ be a set of generators of
${\sigma}^{\ve}\cap M$. Denote $S_l^R=\{s\in S|\langle
s,a_l\rangle<\langle R,a_l\rangle\}$ for $l=1,\dots,N$. Then by
Theorem in \cite[p.~243]{al1}
$$T^1_{X_\sigma}(-R)=\Biggl(L\Biggl(\bigcup_{l=1}^N S_l^R\Biggr)\Biggl/\sum_{l=1}^N L(S_l^R)\Biggr)^*\otimes_\RR
\CC,$$ where $L(\dots)$ denotes the $\RR$-vector space of linear
dependencies of the corresponding subset. All of the vector spaces
$L(\dots)$  in the formula are naturally embedded into $L(S)$,
whose elements can be represented by the $w$-tuples of real
coefficients of the relations.

\begin{thm}\label{t:ks} For the family $  { \mathcal X}_{u Q Q_0 \dots Q_k}\rightarrow\CC^k$, the
Kodaira-Spencer map  $\kappa :T_{\CC^k,0}\rightarrow
T^1_{X_\sigma}$ sends the basis vector $\partial/\partial
\lambda_i$ to the element in $T^1_{X_\sigma}(u)$ induced by the
map
$$L(S)\rightarrow\RR,\quad
 (c_1,\dots,c_w)\mapsto-\sum_{j=1}^w[[\min\langle
s_j,Q_i\rangle]] c_j.$$
\end{thm}

\begin{pf} Classification of infinitesimal deformations of an
affine algebraic variety is well described in \cite{e}. We will
follow its procedure, however computation of the Kodaira-Spencer
map for Altmann's toric construction was done in
\cite[pp.~167--171]{al2} with a slightly different description.

Without loss of generality, assume $i=1$. By \cite[p.~80]{h} the
tangent vector $\partial/\partial \lambda_1$ to $\CC^k$ at $0$
corresponds to the ring homomorphism
$\CC[\lambda_1,\dots,\lambda_k]\rightarrow\CC[\varepsilon]/(\varepsilon^2)$,
assigning $\lambda_1\mapsto\varepsilon$, $\lambda_i\mapsto0$, for
$i\ne1$. By the base change the infinitesimal deformation of
$X_\sigma$ corresponding to  $\partial/\partial \lambda_1$ is the
family $${\cal X}'= { \mathcal
X}_{u,Q,Q_0,\dots,Q_k}\times_{{\CC}^k}{\rm
Spec}(\CC[\varepsilon]/(\varepsilon^2))\rightarrow {\rm
Spec}(\CC[\varepsilon]/(\varepsilon^2)).$$  The first step to
classify this infinitesimal deformation is to present it as an
embedded infinitesimal deformation in an affine space. For this
reason  embed $X_\sigma$ in the affine space $\CC^w$ by the
surjective ring homomorphism
$\CC[z_1,\dots,z_w]\rightarrow\CC[{\sigma}^{\ve}\cap M]$, sending
$z_j\mapsto \chi^{s_j}$, where $s_1,\dots,s_w\ $ are the
generators of ${\sigma}^{\ve}\cap M$. Denote by $J$ the kernel of
this ring homomorphism, then $\CC[{\sigma}^{\ve}\cap M]\simeq
\CC[z_1,\dots,z_w]/J$. The infinitesimal deformation ${\cal X}'$
of $X_\sigma$  corresponds to the ring $$\tilde
A=\Bigl(\bigl(\CC[x_\xi\mid\xi\in\tilde\sigma(1)]\otimes\CC[\varepsilon]/(\varepsilon^2)\bigr)\bigl/I_\varepsilon\Bigr)^{G(\tilde\sigma)},$$
where the ideal $I_\varepsilon$ is generated by
\begin{equation}\label{e:reldef}(1-x^{-e_1^*}-\varepsilon x^{\tilde u-e_1^*})\prod_{ \langle
e_1^*,v_\xi\rangle>0 } x_\xi^{ \langle e_1^*,v_\xi\rangle},\quad
 (1-x^{-e_i^*})\prod_{ \langle e_i^*,v_\xi\rangle>0
} x_\xi^{ \langle e_i^*,v_\xi\rangle},\quad
i=1,\dots,k,\end{equation} where $\tilde
u=u-\sum_{i=1}^k[[\min\langle u,Q_i\rangle]] e_i^*$. Associated to
the embedding $X_\sigma\hookrightarrow{\cal X}'$ we have a
surjective ring homomorphism $\tilde
A\rightarrow\CC[{\sigma}^{\ve}\cap M]$. We lift generator
$\chi^{s_j}$ of $\CC[{\sigma}^{\ve}\cap M]$ by this surjective map
to the coset of $x^{\tilde{s}_j}$, where
$\tilde{s}_j=s_j-\sum_{i=1}^k[[\min\langle
 s_j,Q_i\rangle]]e_i^*$. From the proof of Lemma~\ref{l:1},
 $\tilde{s}_j\in{\tilde{\sigma}}^\ve\cap
\tilde{M}$. Then by \cite[p.~414]{e}, we get a surjective
homomorphism
\begin{equation}\label{e:1}\CC[z_1,\dots,z_w]\otimes
\CC[\varepsilon]/(\varepsilon^2)\rightarrow \tilde
A,\end{equation}
 sending $z_j$ to the coset of $x^{\tilde{s}_j}$,
whose kernel we denote by $\tilde J$. From \cite[p.~178]{e} we
know that $\tilde J/J\varepsilon$ corresponds to the graph of a
homomorphism
$$\psi\in{\rm Hom}(J,\CC[{\sigma}^{\ve}\cap M])\simeq{\rm
Hom}(J/J^2,\CC[{\sigma}^{\ve}\cap M])$$ representing an element of
$$T^1_{X_\sigma}={\rm Hom}(J/J^2,\CC[{\sigma}^{\ve}\cap M])/{\rm
Hom}(\CC[{\sigma}^{\ve}\cap M]^w,\CC[{\sigma}^{\ve}\cap M]).$$ In
practice,  take a binomial generator $z^a-z^b$ of the ideal $J$,
where $a=( a_1,\dots, a_w)$, $b=(b_1,\dots,b_w)\in \ZZ_{\ge0}^w$
with $\ds\sum_{j=1}^w( a_j-b_j)s_j=0$, lift it to
$\CC[z_1,\dots,z_w]\otimes \CC[\varepsilon]/(\varepsilon^2)$ and
then map it by $(\ref{e:1})$ to the coset of $x^{\sum_{j=1}^w a_j
\tilde s_j}-x^{\sum_{j=1}^w b_j \tilde s_j}$ in $\tilde A$. Since
$\sum_{j=1}^w a_j \tilde s_j -\sum_{j=1}^w b_j \tilde s_j\in
\ker(\tilde M\rightarrow M)$, we can write $\sum_{j=1}^w a_j
\tilde s_j- \sum_{j=1}^w b_j\tilde s_j= \sum_{i=1}^k\alpha_i^+
e_i^*-\sum_{i=1}^k\alpha_i^- e_i^*$ with
$\alpha_i^+\cdot\alpha_i^-=0$ and $\alpha_i^+,\alpha_i^-\in
\ZZ_{\ge0}$ for all $i$. Denote
$$q=\sum_{j=1}^w a_j
\tilde s_j-\sum_{i=1}^k\alpha_i^+ e_i^*=\sum_{j=1}^w b_j\tilde
s_j-\sum_{i=1}^k\alpha_i^- e_i^*,$$   where without losing
generality we assume $\alpha_1^-=0$. Note
$q\in{\tilde{\sigma}}^\ve\cap \tilde{M}$.
    Then in the ring $\tilde A$:
\begin{multline*}
x^{\sum_{j=1}^w \alpha_j \tilde s_j}-x^{\sum_{j=1}^w \beta_j
\tilde s_j}=x^q(x^{\sum\alpha_i^+ e_i^*}-x^{\sum\alpha_i^-
e_i^*})= \\
=x^q(  (1+\varepsilon x^{\tilde u })^{ \alpha_1^+}-
1)=\varepsilon\alpha_1^+ x^{q+\tilde u }
\end{multline*}
modulo $\varepsilon^2$ and modulo the relations (\ref{e:reldef}).
Dividing this element by $\varepsilon$ and mapping by $\tilde
A\rightarrow\CC[{\sigma}^{\ve}\cap M]$ we get the image of
$z^a-z^b$ by $\psi$ is
$$\alpha_1^+ \chi^{\sum_{j=1}^w a_j
 s_j+u} = \varphi(a-b)\chi^{\sum_{j=1}^w a_j
 s_j+u},$$
 where $\varphi(c_1,\dots,c_w)=- \sum_{j=1}^w[[\min\langle
s_j,Q_1\rangle]] c_j$.

By the correspondence in Theorem~(3.4) in \cite{al3} we are left
to show that $\varphi$ vanishes on $L(S_l^{-u})$. If
$(c_1,\dots,c_w)\in L(S_l^{-u})$ then $c_j\ne0$ only when $\langle
s_j,a_l\rangle<\langle -u,a_l\rangle$. In this case $a_l/\langle
-u,a_l\rangle$ is a vertex of $\sigma(u)^c$ as noticed in
\cite{al5}. By the condition on $u$ we have $\max\langle
-u,Q\rangle\le1$ and $\sigma(u)^c\subset \RR_{>0} \cdot Q$, whence
there is a vertex $q$ of $Q$ such that $q=ta_l/\langle
-u,a_l\rangle$ for some $0<t\le1$. Consequently, $0\le\langle
s_j,q\rangle<1$ for all $s_j\in S_l^{-u}$. Let $q=q_0+\cdots+q_k$
be the vertex decomposition induced by $Q=Q_0+\cdots+Q_k$. By the
condition $(*)$ on the Minkowski sum decomposition,   all but
possibly one of $q_i$ are lattice points.  Then
$$\langle s_j,q\rangle=\sum_{i=0}^k \langle
 s_j,q_i\rangle \ge \sum_{i=0}^k \min\langle
 s_j, Q_i\rangle\ge\sum_{i=0}^k [[\min\langle
 s_j, Q_i\rangle]]=[[\min\langle
 s_j, Q \rangle]]=0$$
 as in the
proof of Lemma~\ref{l:1}. Since $[[\langle s_j,q\rangle]]=0$ and
with at most one exception $\langle
 s_j,q_i\rangle$, for $i=1,\dots,k$, are integers,  we conclude from the above inequalities that
$\langle
 s_j,q_i\rangle=  \min\langle
 s_j, Q_i\rangle$ for all $s_j\in S_l^{-u}$ and for all but possibly one of $i$. If the
 exceptional $i\ne1$, then $\sum_{j=1}^w[[\min\langle
s_j,Q_1\rangle]] c_j=\sum_{j=1}^w  \langle s_j,q_1 \rangle
c_j=\sum_{j=1}^w  \langle c_j s_j,q_1 \rangle=0$ for
$(c_1,\dots,c_w)\in L(S_l^{-u})$. If   $i=1$ is exceptional
($q_1\notin N$), then $[[\min\langle
s_j,Q_1\rangle]]=-\sum_{i\ne1}  [[\min\langle
 s_j, Q_i\rangle]]=-\sum_{i\ne1}\langle
 s_j, q_i\rangle= \langle
 s_j, -\sum_{i\ne1}q_i\rangle$ for all $s_j\in S_l^{-u}$ and the vanishing of $\varphi$ again follows.
\end{pf}

\begin{rem}\label{r:ks} The
Kodaira-Spencer map $\kappa  $ sends $-\sum_{i=1}^k
{\partial}/{\partial \lambda_i}$ to the element in
$T^1_{X_\sigma}(u)$ induced by the map
$$L (S)\rightarrow\RR,\quad
 (c_1,\dots,c_w)\mapsto-\sum_{j=1}^w[[\min\langle
s_j,Q_0\rangle]] c_j.$$ This follows from the equation $
\sum_{i=0}^k[[\min\langle
 s_j,Q_i\rangle]]=[[\min\langle
 s_j,Q\rangle]]=0$ for all $s_j\in\bigcup_{l=1}^N S_l^{-u}$ as in the
proof of the above theorem.
\end{rem}

While comparison of  deformations ${ \mathcal X}_{u Q Q_0 \dots
Q_k}$ for different $u$ and Minkowski sum decompositions
$Q=Q_0+\cdots+Q_k$ requires a careful study and different $Q$ with
the same $u$ can lead to   isomorphic deformations, we expect that
the families ${ \mathcal X}_{u\sigma(u)Q_0\dots Q_k}$, for $u\in
M$, which are isomorphic to the  non-toric Altmann's construction,
are sufficient to generate all unobstructed infinitesimal
deformations of $X_\sigma$ by Kodaira-Spencer map.
 Similar to Proposition in \cite[p.~172]{al2} we get a criterion
 when this deformation family   has a
 non-trivial Kodaira-Spencer map, which implies a non-trivial
 deformation.

\begin{pr} For the family $f: { \mathcal X}_{u\sigma(u)Q_0\dots Q_k}\rightarrow\CC^k$, the
Kodaira-Spencer map $$\kappa :T_{\CC^k,0}\rightarrow
T^1_{X_\sigma}$$ is trivial if and only if either $k$ of the
summands $Q_0,Q_1,\dots,Q_k$ are cones with a vertex at a lattice
point or $\sigma(u)$ is a lattice polyhedron and  the polyhedra
$Q_0,Q_1,\dots,Q_k$ are either homothetic to $\sigma(u)$  or equal
to cones with a vertex at a lattice point.
\end{pr}

\begin{pf} As in \cite{s}, polyhedra $P$ and $Q$ in $N_\RR$ are {\it
homothetic} if $Q=\delta\cdot P+n$ for some $\delta>0$ and $n\in
N_\RR $. Suppose $Q_i=\delta_i\cdot\sigma(u)+n_i$ for some
$\delta_i>0$ and $n_i\in N_\RR $.  If $\sigma(u) $ is a lattice
polyhedron, then its every Minkowski summand is so by the
condition $(*)$ on the  Minkowski sum decomposition.  Hence,
$$ \sum_{j=1}^w[[\min\langle
s_j,Q_i\rangle]] c_j=\sum_{j=1}^w \min\langle s_j,Q_i\rangle
c_j=\sum_{j=1}^w \delta_i c_j\min\langle
s_j,\sigma(u)\rangle+\sum_{j=1}^w \langle c_j s_j,n_i\rangle =0$$
 for every $(c_1,\dots,c_w)\in L\bigl(\bigcup_{l=1}^N S_l^{R}\bigr)$
 because
  $\min\langle
s_j,\sigma(u) \rangle=0$ as in Remark~\ref{r:ks}. Thus, by
Theorem~\ref{t:ks}, the Kodaira-Spencer map $\kappa $ is trivial.
If $Q_i $ is a cone with a vertex at a lattice point $q_i$, we
again get
 $ \sum_{j=1}^w[[\min\langle
s_j,Q_i\rangle]] c_j=  \langle \sum_{j=1}^w c_j s_j,q_i\rangle =0$
for  $(c_1,\dots,c_w)\in L\bigl(\bigcup_{l=1}^N S_l^{R}\bigr)$.
Using Remark~\ref{r:ks}, we conclude that $k$ of the vectors $
{\partial}/{\partial \lambda_1} ,\dots,
 {\partial}/{\partial \lambda_k}$,
$ -\sum_{i=1}^k {\partial}/{\partial \lambda_i}$ are sent to $0$
by the Kodaira-Spencer map $\kappa $  if $k$ of the summands
$Q_0,Q_1,\dots,Q_k$ are  cones with a lattice vertex. Therefore,
$\kappa $ is trivial in this case as well.

Conversely, suppose that the Kodaira-Spencer map  $\kappa $ is
trivial, then  $$  \sum_{j=1}^w[[\min\langle s_j,Q_i\rangle]]
c_j=0$$ for every $ \sum_{j=1}^w c_j s_j=0$ where $c_j=0$ if
$s_j\notin\bigcup_{l=1}^N S_l^{-u}$. Hence,  we get $
[[\min\langle s_j,Q_i\rangle]]=\langle s_j,n_i\rangle$ for all
$s_j\in\bigcup_{l=1}^N S_l^{-u}$ and some $n_i\in N_\RR$
independent of $j$. Then  $[[\min\langle s_j,Q_i-n_i\rangle]]=0$.
The facets of polyhedron $\sigma(u)$ in the hyperplane where $u$
has value $-1$ are determined by the conditions $\min\langle
s_j,\sigma(u)\rangle=0$ for a subset of $s_j$ in $\bigcup_{l=1}^N
S_l^{-u}$. But for such $j$,
$$0=\min\langle s_j,\sigma(u)\rangle= \sum_{i=0}^k \min\langle
s_j,Q_i\rangle,$$ whence  $\min\langle s_j,Q_i\rangle$ are
integers by the condition $(*)$ on the Minkowski sum
decomposition. Therefore, we get stricter condition $\min\langle
s_j,Q_i-n_i\rangle =0$ for $s_j$ that were normals to the facets
of $\sigma(u)$. Since $Q_i+n_i$ is a Minkowski summand of
$\sigma(u)$ in a parallel hyperplane, it follows that $Q_i-n_i$ is
either   a multiple of $\sigma(u)$ or   the recession cone of
$\sigma(u)$. Thus, $Q_i=\delta_i\sigma(u)+n_i$ for
$0<\delta_i\le1$ or $Q_i$ is a cone. If $\sigma(u)$ is a lattice
polyhedron, then we are done. Suppose $a_l/\langle -u,a_l\rangle$
is a non-lattice vertex of $\sigma(u)$, or equivalently $\langle
-u,a_l\rangle>1$. Since $S$ is generating ${\sigma}^{\ve}\cap M$,
there is always $s_{j}$ such that $\langle s_{j},a_l\rangle=1$.
Then $s_j\in S_l^{-u}$ and $\langle s_j,n_i\rangle$ is an integer
by the above. Now the Minkowski sum decomposition
$\sigma(u)=Q_0+Q_1+\cdots+Q_k$ induces decomposition of the vertex
$$\frac{a_l}{\langle -u,a_l\rangle}=\Bigl(\delta_0\frac{a_l}{\langle -u,a_l\rangle}+n_0\Bigr)+\cdots
+\Bigl(\delta_k\frac{a_l}{\langle -u,a_l\rangle}+n_k\Bigr)$$ into
the sum of vertices of $Q_0,Q_1,\dots,Q_k$, where we assume that
$\delta_i=0$ if $Q_i$ is a cone. By the condition $(*)$ on this
decomposition only one of the summands can be non-lattice. But if
$\delta_i\frac{a_l}{\langle -u,a_l\rangle}+n_i$ is a lattice point
then $$\Bigl\langle s_j, \delta_i\frac{a_l}{\langle
-u,a_l\rangle}+n_i\Bigr\rangle-\langle s_j,n_i\rangle=
\delta_i\frac{\langle s_j,a_l\rangle}{\langle
-u,a_l\rangle}=\frac{\delta_i}{\langle -u,a_l\rangle}$$ is an
integer which can happen only if $\delta_i=0$ implying $Q_i$ is a
cone.
\end{pf}

\section{Minkowski sum decompositions of polyhedral
complexes.}\label{s:mindec}

In order to generalize the construction of deformations of affine
toric varieties and the construction of deformations of  blow ups
of Fano toric varieties in \cite{m2} we have to introduce the
notion of a Minkowski sum decomposition of a polyhedral complex.
While elements of this  notion have appeared in   the mixed
subdivisions of a Minkowski sum decomposition of a polytope
introduced in \cite{ps} and \cite{st} (also, see
\cite[p.~360]{clo}, \cite{hrs}, \cite{sa}), this combinatorial
object has not been formally introduced and studied in a general
case.

A   polyhedral complex in a real vector space is  a collection of
convex polyhedra with their faces such that intersection of
 any two polyhedra in the complex must be a face of each.
 We define a {\it Minkowski sum decomposition  of a polyhedral
 complex} as
  a Minkowski sum decomposition of each polyhedra   in the
 complex such that for any two different intersecting polyhedra $Q$ and $Q'$ of the
 complex the decompositions $Q=Q_1+\cdots+Q_n$ and
 $Q'=Q'_1+\cdots+Q'_n$ induce compatible decompositions on their
 common face $Q\cap Q'=Q_1\cap Q'_1+\cdots+Q_n\cap Q'_n$. As in
 \cite{s}, we use the convention that if a polyhedra   is
 unbounded, then all of its summands have the same recession cone.

\begin{ex}\label{ex:line} Consider the polyhedral complex obtained by subdividing
a line segment   in $N_\RR$ by points $n_1,\dots, n_k$:
$$ \setlength{\unitlength}{1cm}
\begin{picture}(8,2)
\put(-1,1.2){$n_{1}$} \put(-1,1){\circle*{0.1}}
\put(0,1){\circle*{0.1}} \put(1,1){\circle*{0.1}}
\put(2,1){\circle*{0.1}} \put(3,1){\circle*{0.1}}
\put(4,1){\circle*{0.1}}
 \put(0,1.2){$n_{2}$}
\multiput(1.5,1.2)(.4,0){3}{\circle*{0.03}} \put(3,1.2){$n_{k-1}$}
  \put(4,1.2){$n_k$}
\put(-1,1){\line(1,0){5}}
\end{picture}\label{e:pic}
$$
Polyhedra in this complex are line segments $[n_i,n_{i+1}]$ with
their  vertices. We decompose them as follows:
$$
\begin{array}{ccccc} [n_1,n_2]&[n_2,n_3]&\cdots&[n_{k-2},n_{k-1}]&[n_{k-1},n_k] \\
\parallel  &\parallel  &\cdots  &\parallel  &\parallel  \\ {[n_1,n_2]} & n_2 & \cdots & n_2 & n_2\\
 +&+&\cdots&+&+\\
  0 & [0,n_3-n_2] &   \cdots& n_3-n_2& n_3-n_2\\
  +&+&\cdots&+&+\\ \vdots&\vdots&&\vdots  &\vdots\\
   +&+&\cdots&+&+ \\
  0&0&\cdots&[0,n_{k-1}-n_{k-2}]&n_{k-1}-n_{k-2}\\
  +&+&\cdots&+&+ \\
  0&0&\cdots&0&[0,n_k-n_{k-1}],
\end{array}
$$
where by $0$ we denoted the origin point in $N_\RR$. Note that
each vertex $n_i$ has the same induced decomposition from the
nearby line segments. Thus, we have a Minkowski sum decomposition
of the polyhedral complex.
\end{ex}

\begin{ex}\label{ex:2nd} Another, less trivial example, of a Minkowski sum
decomposition of a polyhedral complex comes from the skeleton of a
reflexive polytope. Consider the convex hull of points $(-1,-1)$,
$(2,-1)$, $(-1,2)$ in $\RR^2$. This is the dual of of the
reflexive polytope coming from the fan of $\PP^2$. Subdividing the
skeleton of this reflexive polytope by the lattice points
$n_1,\dots,n_7$ inside   its edges we obtain the following
polyhedral complex.

$$ \setlength{\unitlength}{1cm}
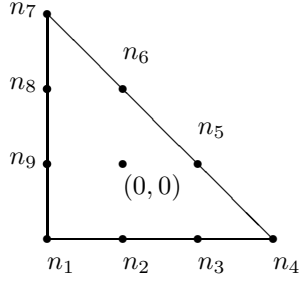
\begin{picture}(8,4)
\put(1,0.6){$n_{1}$} \put(2,0.6){$n_{2}$} \put(3,0.6){$n_{3}$}
\put(2 ,1.6){$(0,0)$} \put(4,0.6){$n_{4}$}  \put(3,2.4){$n_{5}$}
\put(2,3.4){$n_{6}$} \put(0.5,4){$n_{7}$} \put( 0.5,3){$n_{8}$}
\put( 0.5,2){$n_{9}$} \put(1,1){\circle*{0.1}}
\put(2,2){\circle*{0.1}} \put(1,2){\circle*{0.1}}
\put(1,2){\circle*{0.1}} \put(1,3){\circle*{0.1}}
\put(1,3){\circle*{0.1}} \put(1,4){\circle*{0.1}}
\put(2,1){\circle*{0.1}} \put(3,1){\circle*{0.1}}
\put(4,1){\circle*{0.1}} \put(3,2){\circle*{0.1}}
\put(2,3){\circle*{0.1}} \put(1,1){\line(0,1){3}}
\put(1,1){\line(1,0){3}} \put(1,4){\line(1,-1){3}}
\end{picture}\label{e:pic}
$$
Here is a nontrivial Minkowski sum decompostion of this complex:
$$
\begin{array}{ccccccccc} [n_1,n_2]=&[n_1,n_2]&+&0&+&0\\
{[n_2,n_3]}=&n_2&+&[0,n_3-n_2]&+&0\\
{[n_3,n_4]}=&n_2&+&n_3-n_2&+&[0,n_4-n_3]\\
{[n_4,n_5]}=&[n_2,n_9]&+&n_3-n_2&+& n_4-n_3 \\{[n_5,n_6]}=&n_9&+&[n_3-n_2,n_6]&+&n_4-n_3\\
{[n_6,n_7]}=&n_9&+&n_6&+&[n_4-n_3,n_6]  \\
   {[n_7,n_8]}= & [n_9,n_1] &+&n_6&+&n_6\\
  {[n_8,n_9]}= &  n_1 &+&[n_6,0]&+&n_6\\   {[n_9,n_1]}= &  n_1 &+& 0
  &+&[n_6,0]
\end{array}
$$
\end{ex}

The above two examples are related to mixed subdivisions as
follows. Similar to \cite{ps} and  \cite{sa},  given polytopes
$P_1,\dots,P_n$ in $\RR^d$, a {\it mixed subdivision} of
$P=P_1+\cdots+P_n$ is a polyhedral subdivision of $P$  into
polytopes  (called cells) of the form $Q=Q_1+\cdots+Q_n$, where
$Q_i$ is a polytope contained in $P_i$, such that for any two
cells $Q=Q_1+\cdots+Q_n$ and $Q'= Q'_1+\cdots+Q'_n$ the
intersection $Q_i\cap Q_i'$ is a face of $Q_i$ and $Q_i'$, for
$i=1,\dots,n$. (The intersection can be
 an empty face.)

\begin{rem}
Definition of the mixed subdivision in \cite{sa} assumes that the
vertices of $Q_i$ are among the vertices of $P_i$, while \cite{ps}
requires that $Q_i$ is a face of $P_i$.
\end{rem}

The polyhedral Cayley Trick in \cite{st} and  in Theorem~3.1 of
\cite{hrs}   extends to the general case as follows. Given
polytopes $P_1,\dots,P_n$ in $\RR^d$, define the Cayley embedding
polytope:
$${\cal C}(P_1,\dots,P_n)={\rm Conv}(P_1+r_1,\dots,P_n+r_n)$$ in
$\RR^d\oplus \RR^{n-1}$, where
  $r_1,\dots,r_n$ is an affine basis of $\RR^{n-1}$.

\begin{thm}\label{t:mix} There is a  one-to-one correspondence
between mixed subdivisions   of $P_1+\cdots+P_n$  and polyhedral
subdivisions    of ${\cal C}(P_1,\dots,P_n)$, whose
$0$-dimensional cells are contained in the faces
$P_1+r_1,\dots,P_n+r_n$:
\\ (i)  If ${\cal Q}=\{ Q_1+\cdots+Q_n\mid Q_i\subseteq
P_i\text{ for } i=1,\dots n\}$ is a mixed subdivision of
$P_1+\cdots+P_n$, then the cells ${\cal C}(Q_1,\dots,Q_n)$
form a subdivision of ${\cal C}(P_1,\dots,P_n)$;\\
(ii) If ${\cal B}$ is a polyhedral subdivision   of ${\cal
C}(P_1,\dots,P_n)$, whose $0$-dimensional cells are contained in
the faces $P_1+r_1,\dots,P_n+r_n$, then the intersection of ${\cal
B}$ with the affine space $\RR^d\oplus \{\frac{1}{n}\sum_{i=1}^n
r_i\}$ gives a mixed subdivision of $\frac{1}{n}(P_1+\cdots+P_n)$.
\end{thm}

\begin{pf} The proof is the same as in Propositions~3.4 and 3.5 in
\cite{hrs}. It also shows that  for the definition of the mixed
subdivision it is sufficient that for any two adjacent cells
$Q=Q_1+\cdots+Q_n$ and $Q'= Q'_1+\cdots+Q'_n$ their common face
$Q\cap Q'$ has a compatible induced Minkowski sum decomposition.
\end{pf}

\begin{ex} Polyhedral subdivision of the line segment $[n_1,n_k]$
by the points $n_2,\dots,n_{k-1}$ with the Minkowski sum
decompositions of the line segments in Example~\ref{ex:line} is a
mixed subdivision, but it is easy to see that a different
Minkowski sum decomposition of the polyhedral subdivision gives a
different mixed subdivision structure of the same polyhedral
complex. These different structures are encoded by different
polyhedral subdivisions of  the Cayley embedding polytope.
\end{ex}

\begin{lem} Let $P$ be a polytope and ${\cal Q}$ be a polyhedral
complex  whose underlying set is $P$. If
 $Q =Q_1+\cdots+Q_n$,  $ Q\in{\cal  Q}$,  is a Minkowski
sum decomposition of  ${\cal  Q}$, then   ${\cal  Q}$ is a mixed
subdivision of $P=P_1+\cdots+P_n$, where $P_i=\bigcup_{Q\in{\cal
Q}}Q_i$ is convex.
\end{lem}

By this lemma,  a    mixed subdivision  of a Minkowski sum of
polytopes $P=P_1+\cdots+P_n$ is a polyhedral complex ${\cal Q}$
with underlying set   $|{\cal Q}|=P$ and a Minkowski sum
decomposition of the   complex $Q=Q_1+\cdots+Q_n$, for $Q\in{\cal
Q}$, such that $Q_i\subseteq P_i$. Therefore, Minkowski sum
decompositions of a polyhedral subdivision of a polytope $P$
correspond  to polyhedral subdivisions of a  Cayley embedding
polytope for a Minkowski sum $P=P_1+\cdots+P_n$.

\begin{ex}  Consider the polyhedral
subdivision of the line segment $[0,3]$ in $\RR$ by 1 and 2 and a
Minkowski sum decomposition of this complex $[0,1]=[0,1]+0$,
$[1,2]=1+[0,1]$, $[2,3]=[1,2]+1$, which is a mixed subdivision of
$[0,2]+[0,1]=[0,3]$.
 For the affine basis $\{-1,1\}$   of $\RR^1$, we get the Cayley
 embedding polytope   ${\cal
 C}([0,2],[0,1])=\co((0,-1),(2,-1),(0,1),(1,1))$ in $\RR^2\cong
 \RR\oplus \RR$, subdivided by  ${\cal
 C}([0,1],\{0\})$, ${\cal
 C}(\{1\},[0,1])$ and ${\cal
 C}([1,2],\{1\})$. Intersection of this subdivision with $\RR\oplus\{0\}$ gives a
 mixed subdivision of $\frac{1}{2}[0,3]$, which is a scaled version of the mixed subdivision of
 $[0,3]$.
$$ \setlength{\unitlength}{1cm}
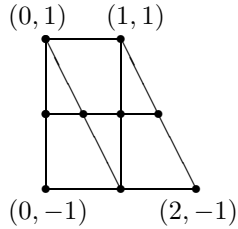
\begin{picture}(4,4)
\put(0.5,0.6){$(0,-1)$} \put(2.5,0.6){$(2,-1)$}\put(1.8,3.2){$(1,
1)$}   \put( 0.5,3.2){$(0,1)$} \put(1,1){\circle*{0.1}}
\put(2,2){\circle*{0.1}} \put(1,2){\circle*{0.1}}
\put(1,2){\circle*{0.1}} \put(1,3){\circle*{0.1}}
\put(1,3){\circle*{0.1}} \put(2,1){\circle*{0.1}}
\put(3,1){\circle*{0.1}}
\put(2.5,2){\circle*{0.1}}\put(1.5,2){\circle*{0.1}}
\put(2,3){\circle*{0.1}} \put(2,1){\line(0,1){2}}
\put(1,3){\line(1,0){1}} \put(1,1){\line(0,1){2}}
\put(1,1){\line(1,0){2}}\put(1,2){\line(1,0){1.5}}
\put(3,1){\line(-1,2){1}} \put(2,1){\line(-1,2){1}}
\end{picture}
$$
\end{ex}

Example~\ref{ex:2nd} provides a Minkowski sum decomposition of a
polyhedral subdivision of the boundary of a polytope, and, in
particular, induces mixed subdivisions on the faces of the
polytope. However, one can show that there is no Minkowski sum
decomposition of a polyhedral subdivision of the polytope
compatible with the Minkowski sums in this example. In other
words, compatible mixed subdivisions on the boundary of a polytope
may not have an extension  to  a  mixed subdivision of the whole
polytope.

As a consequence of Theorem~\ref{t:mix}, we can form  Minkowski
sum decompositions of a polyhedral subdivision  of the boundary of
the Minkowski sum polytope $P=P_1+\cdots+P_n$ by subdividing the
boundary of the associated Cayley polytope $ {\cal
C}(P_1,\dots,P_n)$.

\begin{cor}\label{c:subd} Let $P_1,\dots,P_n$ be polytopes in
$\RR^d$. Let ${\cal B}$ be a polyhedral subdivision   of the union
of faces of $ {\cal C}(P_1,\dots,P_n)$ which intersect the affine
space $\RR^d\oplus \{\frac{1}{n}\sum_{i=1}^n r_i\}$, such that the
$0$-dimensional cells of ${\cal B}$ are contained in
$P_1+r_1,\dots,P_n+r_n$. Then   Minkowski sums $B_1+\cdots+B_n$,
for $B\in {\cal B}$, such that $B={\cal C}(B_1,\dots,B_n)$ with
$B_i\subseteq P_i$, form a Minkowski sum decomposition of a
polyhedral subdivision of the boundary of $P=P_1+\cdots+P_n$.
\end{cor}

\section{Embedding of toric varieties via Minkowski sum decompositions of polyhedral
complexes.}\label{s:embed}

Using Minkowski sum decompositions of polyhedral complexes we will
naturally generalize  the embedding of an affine toric variety
into a higher dimensional affine toric variety from
Section~\ref{s:1}.

\begin{defn}\label{d:compat}  Let $\Sigma$ be a fan in $N_\RR$.  We say that   a
  polyhedral complex ${\cal Q}$  in $N_\RR$   is {\it
$\Sigma$-compatible} if  $Q^\sigma:=\sigma\cap|{\cal Q}|\in{\cal
Q}$ for all $\sigma\in\Sigma$,  and every polyhedra in ${\cal Q}$
is a face of $Q^\sigma$  for some $\sigma\in \Sigma$.
\end{defn}

 It is clear that the fan $\Sigma$ is a
  $\Sigma$-compatible polyhedral complex. Another example: proper faces of a   polytope with an origin in its interior will form
  a $\Sigma$-compatible polyhedral complex, if the fan $\Sigma$ is formed by the cones generated by the proper faces
  of the polytope.

\begin{defn}\label{d:comp}
Suppose that ${\cal  Q} $ is  a $\Sigma$-compatible polyhedral
complex in $N_\RR$ and let
$$Q^\sigma=Q^\sigma_0+Q^\sigma_1+\cdots+Q^\sigma_k,  \quad \sigma\in\Sigma,$$ be  a Minkowski
sum decomposition of polyhedral complex ${\cal  Q}$. We will say
that this decomposition is {\it $\Sigma$-separated} if for every
pair of cones $\sigma_1,\sigma_2\in\Sigma$, there exists $m\in
M_\RR$ such that $m\in \sigma_1^\ve$, $-m\in\sigma_2^\ve$,
$\sigma_1\cap m^\perp=\sigma_2\cap m^\perp$, and, for all $i$,
either $\min\langle m,Q_i^{\sigma_1}\rangle>\max\langle
m,Q_i^{\sigma_2}\rangle$, or $$\{q\in Q_i^{\sigma_1}\mid\langle m,
q\rangle=\min\langle m,Q_i^{\sigma_1}\rangle\}= \{q\in
Q_i^{\sigma_2}\mid\langle m, q\rangle=\max\langle
m,Q_i^{\sigma_2}\rangle\}.$$  \end{defn}

\begin{rem} The first three conditions say that $m$ determines a
hyperplane which intersects both cones in their common face
$\sigma_1\cap\sigma_2$. The other conditions imply that the
Minkowski summands $Q_i^{\sigma_1}$ and $Q_i^{\sigma_2}$ are
separated by  parallel hyperplanes with the normal $m$. These
conditions are sufficient to form the following fan.
\end{rem}

As in Section~\ref{s:1},  we have extended lattice
$\tilde{N}=N\oplus\ZZ^k$ and the cones
$$\tilde{\sigma}=\langle \sigma,
 Q^\sigma_0-e_1-\dots-e_k,Q^\sigma_1+e_1,\dots,Q^\sigma_k+e_k\rangle$$
in $\tilde{N}_\RR$
 for each $\sigma\in\Sigma$.

\begin{lem}\label{l:tilfan} Let $\Sigma$ be a  fan in $N_\RR$ and ${\cal  Q} $ be  a $\Sigma$-compatible   polyhedral
complex with $0\notin |{\cal  Q}|$, and let
 $Q^\sigma=Q^\sigma_0+Q^\sigma_1+\cdots+Q^\sigma_k$, for $\sigma\in\Sigma$,   be  a $\Sigma$-separated  Minkowski sum
decomposition of polyhedral complex ${\cal  Q}$. Then  the
collection of cones $\tilde{\Sigma}_{\cal
Q}=\{\tilde{\sigma}\mid\sigma\in\Sigma\}$ is a fan in
$\tilde{N}_\RR$. Additionally, if $Q^\sigma_i $ are rational
polyhedra for all $i$ and $\sigma$, then the fan
$\tilde{\Sigma}_{\cal Q}$ is rational polyhedral.

\end{lem}

\begin{pf} By Lemma~\ref{l:coneprop}, the cones $\tilde{\sigma}$
are strongly convex.
 By
Definition~\ref{d:comp}, for a pair of cones
$\sigma_1,\sigma_2\in\Sigma$, there exists $m\in M_\RR$ satisfying
the special properties. We need to show existence of   $\tilde
m=m+\sum_{i=1}^k a_i e_i^*$ with $a_i\in\RR$, such that
\begin{equation}\label{e:conditions}\tilde m\in \tilde\sigma_1^\ve,\quad -\tilde
m\in\tilde\sigma_2^\ve, \quad\tilde\sigma_1\cap\tilde
m^\perp=\tilde\sigma_2\cap\tilde m^\perp.\end{equation}
 The first two conditions are equivalent to $\langle
\tilde m,
 q_i^1+e_i\rangle\ge \min\langle
m,Q_i^{\sigma_1}\rangle+a_i\ge0$, $\langle \tilde m,
 q_0^1-\sum_{i=1}^ke_i\rangle\ge \min\langle
m,Q_0^{\sigma_1}\rangle-\sum_{i=1}^k a_i\ge0$,
 for $q_i^1\in Q_i^{\sigma_1}$,  and $\langle - \tilde m,
 q_i^2+e_i\rangle\ge \min\langle
-m,Q_i^{\sigma_2}\rangle-a_i=-\max\langle
m,Q_i^{\sigma_2}\rangle-a_i\ge0$, $\langle -\tilde m,
 q_0^2-\sum_{i=1}^ke_i\rangle\ge \min\langle
-m,Q_0^{\sigma_2}\rangle+\sum_{i=1}^k a_i=-\max\langle
m,Q_0^{\sigma_2}\rangle+\sum_{i=1}^k a_i\ge0$,
 for $q_i^2\in Q_i^{\sigma_1}$. Then  (\ref{e:conditions}) reduces to
$$-\max\langle m,Q_i^{\sigma_2}\rangle\ge a_i\ge-\min\langle
m,Q_i^{\sigma_1}\rangle, \text{ for } i=1,\dots,k,$$ $$
\min\langle m,Q_0^{\sigma_1}\rangle \ge \sum_{i=1}^k a_i\ge
\max\langle m,Q_0^{\sigma_2}\rangle$$ with strict inequalities $>$
replacing $\ge$ in the case $\min\langle
m,Q_i^{\sigma_1}\rangle\ne\max\langle m,Q_i^{\sigma_2}\rangle$ for
$i=0,\dots,k$.
 It follows from
Definitions~\ref{d:compat} and \ref{d:comp} that either
$\min\langle m,Q^{\sigma_1}\rangle=0=\max\langle
m,Q^{\sigma_2}\rangle$ or $\min\langle
m,Q^{\sigma_1}\rangle>0>\max\langle m,Q^{\sigma_2}\rangle$. The
first case implies that $\min\langle m,Q_i^{\sigma_1}\rangle=
\max\langle m,Q_i^{\sigma_2}\rangle$ for $i=0,\dots,k$ and we can
take $a_i=-\min\langle m,Q_i^{\sigma_1}\rangle= -\max\langle
m,Q_i^{\sigma_2}\rangle$ to get $\tilde m$ satisfying
(\ref{e:conditions}). In the second case we take
$$a_i= \frac{\max\langle
m,Q^{\sigma_2}\rangle\cdot\min\langle
m,Q_i^{\sigma_1}\rangle-\min\langle
m,Q^{\sigma_1}\rangle\cdot\max\langle m,Q_i^{\sigma_2}\rangle
}{\min\langle m,Q^{\sigma_1}\rangle-\max\langle
m,Q^{\sigma_2}\rangle},$$ for which
\begin{multline*}\sum_{i=1}^k a_i= \frac{\max\langle
m,Q^{\sigma_2}\rangle\cdot\sum_{i=1}^k\min\langle
m,Q_i^{\sigma_1}\rangle-\min\langle
m,Q^{\sigma_1}\rangle\cdot\sum_{i=1}^k\max\langle
m,Q_i^{\sigma_2}\rangle }{\min\langle
m,Q^{\sigma_1}\rangle-\max\langle
m,Q^{\sigma_2}\rangle}=\\=\frac{-\max\langle
m,Q^{\sigma_2}\rangle\cdot \min\langle
m,Q_0^{\sigma_1}\rangle+\min\langle m,Q^{\sigma_1}\rangle\cdot
\max\langle m,Q_0^{\sigma_2}\rangle }{\min\langle
m,Q^{\sigma_1}\rangle-\max\langle
m,Q^{\sigma_2}\rangle}.\end{multline*} Hence, the required
inequalities  to get $\tilde m$ satisfying (\ref{e:conditions})
follow.
\end{pf}

\begin{ex}\label{e:subdiv}
Let $P$ be a polytope $N_\RR$ with the origin in its interior and
let $P=P_0+\cdots +P_k$ be a Minkowski sum decomposition into
rational polytopes.  Let ${\cal B}$ be a polyhedral subdivision of
the union of faces of the Cayley polytope $\co(P_0-\sum_{i=1}^k
e_i,P_1+e_1,\dots,P_k+e_k)$ not intersecting $N_\RR\subset \tilde
N_\RR$. Then the intersection of ${\cal B}$ with $N_\RR$ forms a
polyhedral subdivision  ${\cal Q}$ of the boundary of $P$, which,
by Corollary~\ref{c:subd}, admits a Minkowski sum decomposition
$B_0+\cdots+B_k$, for $B\in {\cal B}$, such that
$B=\co(B_0-\sum_{i=1}^k e_i,B_1+e_1,\dots,B_k+e_k)$ with
$B_i\subseteq P_i$. If $\Sigma$ is the fan    formed by the cones
generated by polyhedra in $\cal Q$, then $\tilde{\Sigma}_{\cal
Q}$, corresponding to the Minkowski sum decomposition of ${\cal
Q}$, is the fan over the polyhedral subdivision ${\cal B}$ of the
Cayley polytope. In this case, the Minkowski sum decomposition of
${\cal Q}$ is automatically $\Sigma$-separated.
\end{ex}

Another example is provided by a Minkowski sum decomposition of a
slice of a complete fan:
\begin{lem}\label{l:slice}
Given a complete fan $\Sigma$ in $N_\RR$ and    $u\in M$, then
$\Sigma(u):=\{ \sigma(u) \mid\sigma\in\Sigma\}$, where
$\sigma(u)=\sigma\cap \{n\in N_\RR\mid\langle u,n\rangle=-1\}$, is
a $\Sigma$-compatible   polyhedral complex.  Let  $
\sigma(u)=Q_0+Q_1+\cdots+Q_k$, for $\sigma\in\Sigma$,
  be a Minkowski sum decomposition of $\Sigma(u)$.
  If $|\Sigma|=N_\RR$, then,   the polyhedra $ {\cal C}(Q_0,\dots,Q_k)={\rm
Conv}(Q_0-\sum_{i=1}^ke_i,Q_1+e_1,\dots,Q_k+e_k)$ form a
polyhedral subdivision of the hyperplane in $N_\RR\oplus\RR^k$
where $(k+1)(u-\sum_{i=1}^k\langle u,Q_i\rangle
e_i^*)-\sum_{i=1}^ke_i^*$ takes value $-1$, whence
 a Minkowski sum decomposition of $\Sigma(u)$ is
 $\Sigma$-separated and $\tilde{\Sigma}_{\cal Q}=\{\tilde{\sigma}\mid\sigma\in\Sigma\}$ is a
fan in $\tilde{N}_\RR$.
\end{lem}

\begin{pf}
The arguments to show   that $ {\cal C}(Q_0,\dots,Q_k)$ form a
polyhedral subdivision are similar to the proof of Proposition~3.5
in \cite{hrs}. The key condition is that the polyhedral complex
$\Sigma(u)$ fills up the whole hyperplane $\{n\in N_\RR\mid\langle
u,n\rangle=-1\}$.
\end{pf}

\begin{ex}
 The following picture illustrates a
 slice of a fan in $\RR^3$ and its Minkowski sum decomposition.
$$ \setlength{\unitlength}{0.7cm}
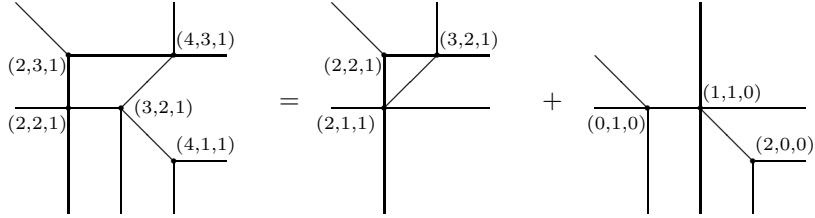
\begin{picture}(20,4)
 \put(2,2){\line(-1,0){1}} \put(2,2){\line(0,-1 ){2}}\put(3,2){\line(0,-1 ){2}}
\put(2,2){\line(0,1){1}}\put(2,2){\line(1,0){1}}
\put(4,3){\line(0,1){1}}\put(4,1){\line(0,-1){1}}
\put(4,3){\line(1,0){1}}\put(4,1){\line(1,0){1}}
\put(3,2){\line(1,-1){1}}\put(3,2){\line(1, 1){1}}
 \put(2,3){\line(1,0){2}}\put(1,4){\line(1,-1){1}}

\put(2,3){\circle*{0.1}} \put(0.7,2.7){\SMALL{ (2,3,1)}}
\put(4,3){\circle*{0.1}}\put(3.9,3.2){\SMALL{ (4,3,1)}}
 \put(2,2){\circle*{0.1}}\put(0.7,1.6){\SMALL{ (2,2,1)}}
   \put(3,2){\circle*{0.1}}\put(3.1,1.9){\SMALL{ (3,2,1)}}
\put(4,1){\circle*{0.1}}\put(3.9,1.2){\SMALL{ (4,1,1)}}

\put( 6,2 ){$=$}

\put(8,2){\line(-1,0){1}} \put(8,2){\line(0,-1 ){2}}
\put(8,2){\line(0,1){1}}\put(8,2){\line(1,0){2}}
\put(8,2){\line(1, 1){1}} \put(9,3){\line(0,1){1}}
 \put(8,3){\line(1,0){2}}\put(7,4){\line(1,-1){1}}

\put(8,2){\circle*{0.1}} \put(6.7 ,1.6){\SMALL{(2,1,1)}}
\put(8,3){\circle*{0.1}} \put(6.7 ,2.7){\SMALL{ (2,2,1)}}
\put(9,3){\circle*{0.1}} \put(8.9,3.2){\SMALL{ (3,2,1)}}

 \put(
11,2 ){$+$}

 \put(12,3){\line(1,-1){1}}
 \put(12,2){\line( 1,0){4}} \put(13,2){\line(0,-1 ){2}}\put(14,4){\line(0,-1 ){4}}
 \put(15,1){\line(1,0){1}}\put(15,1){\line(0,-1){1}}
\put(14,2){\line(1,-1){1}}

\put(13,2){\circle*{0.1}} \put(11.7 ,1.6){\SMALL{ (0,1,0) }}
\put(14,2){\circle*{0.1}} \put(13.9,2.2){\SMALL{ (1,1,0)}}

\put(15,1){\circle*{0.1}} \put(14.9,1.2){\SMALL{ (2,0,0)}}
\end{picture}
$$
A Minkowski sum decomposition of a polyhedral complex can be
tracked by the decomposition of the compact faces and the addition
of the recession cones for unbounded polyhedra.

\end{ex}

Suppose that $\tilde{\Sigma}_{\cal Q}$ is the fan as in
Lemma~\ref{l:tilfan} corresponding to a $\Sigma$-separated
Minkowski sum decomposition
$Q^\sigma=Q^\sigma_0+Q^\sigma_1+\cdots+Q^\sigma_k$ of a
$\Sigma$-compatible   polyhedral complex ${\cal  Q}$ with $0\notin
|{\cal Q}|$. Assume in addition that $Q^\sigma_i $ are rational
polyhedra for all $i$ and $\sigma$ and the Minkowski sums satisfy
the condition:\\    $(*)$ the induced decomposition of a vertex of
$Q^\sigma$
 $${\rm vert}(Q^\sigma )={\rm vert}(Q^\sigma_0)+{\rm vert}(Q^\sigma_1)+\cdots+{\rm
 vert}(Q^\sigma_k)$$
has at most one of the summands not in the lattice $N$.

   The natural inclusion map
${N}_\RR\subset\tilde{N}_\RR$ is compatible with    the fans
$\Sigma$ and $\tilde{\Sigma}_{\cal Q}$ and we have

\begin{thm}\label{t:emb1} Let $\Sigma$ and  $\tilde{\Sigma}_{\cal Q}$ be as above. Associated to the map of
the fan $\Sigma$ to $\tilde{\Sigma}_{\cal Q}$, the morphism
between toric varieties $X_\Sigma\rightarrow
X_{\tilde{\Sigma}_{\cal Q}}$ is an embedding, whose image  is a
complete intersection given by the regular sequence
$$\prod_{ \langle e_i^*,v_\xi\rangle>0 } x_\xi^{ \langle
e_i^*,v_\xi\rangle}-\prod_{ \langle e_i^*,v_\xi\rangle<0 } x_\xi^{
-\langle e_i^*,v_\xi\rangle}=0,$$ for $i=1,\dots,k$, where $x_\xi$
are the homogeneous coordinates of the toric variety
$X_{\tilde{\Sigma}_{\cal Q}}$ corresponding to
  the rays $\xi\in\tilde{\Sigma}_{\cal Q}(1)$ with  primitive
lattice generators $v_\xi$.
\end{thm}

\begin{pf}
 Since the embedding property of
varieties can be checked locally, from Lemma~\ref{l:1} we get the
embedding part of the statement.

To prove the second part it is sufficient to show that locally for
each affine toric chart $X_{\tilde{\sigma}}\subset
X_{\tilde{\Sigma}_{\cal Q}}$  the given equations induce the same
equations as in Theorem~\ref{c:eq}, which determine $X_\sigma$ in
$X_{\tilde{\sigma}}$. We assume that the one-dimensional cones
$\Sigma(1)$ span $N_\RR$, but as noted in Section~\ref{s:toric}
this assumption is not essential if we add homogeneous coordinates
corresponding to ``non-existent'' one-dimensional cones of
$\Sigma$. We also assume  $\dim\tilde\sigma=\dim \tilde N_\RR$.
   Denote by  $\tilde I$
the homogeneous ideal in  $S(\tilde\Sigma_{\cal
Q})=\CC[x_\xi\mid\xi\in\tilde\Sigma_{\cal Q}(1)]$ generated by
$$\prod_{ \langle e_i^*,v_\xi\rangle>0 } x_\xi^{ \langle
e_i^*,v_\xi\rangle}-\prod_{ \langle e_i^*,v_\xi\rangle<0 } x_\xi^{
-\langle e_i^*,v_\xi\rangle}\quad, i=1,\dots,k.$$ By
Propositions~\ref{p:local}, the closed subvariety ${\bf V}(\tilde
I)$ in the affine open toric chart $X_{\tilde{\sigma}}$
corresponds to the ring
$$\bigl( \CC[x_\xi\mid\xi\in\tilde\Sigma_{\cal Q}(1)]_{x^{\hat{\tilde\sigma}}}/\tilde I_{x^{\hat {\tilde\sigma}}}\bigr)^{G(\tilde\Sigma_{\cal Q})}\simeq(\CC[x_\xi\mid\xi\in\tilde\sigma(1)]/\psi(\tilde I))^{G(\tilde\sigma)},$$
where $\psi: \CC[x_\xi\mid\xi\in\tilde\Sigma_{\cal Q}(1)]
\rightarrow\CC[x_\xi\mid\xi\in\tilde\sigma(1)]$
  is evaluation   $x_{\xi}=1$, for $\xi\not\subseteq \tilde\sigma$. It is clear
  that $\psi(\tilde I)$ is the ideal $I$ in
  Proposition~\ref{p:iso} and the result follows.
\end{pf}

\section{Deformations of     toric
varieties by complete intersections.}\label{s:darb}

 In this section we generalize our construction from Section~\ref{s:1}
to the case of an arbitrary toric variety $X_\Sigma$ associated to
a fan $\Sigma$ in $N_\RR$ using realization of the toric variety
as a complete intersection in a higher dimensional toric variety
from Theorem~\ref{t:emb1}.

Let $\tilde\Sigma_{\cal Q}$ be the fan as in Theorem~\ref{t:emb1},
corresponding to a $\Sigma$-separated Minkowski sum decomposition
$Q=Q_0+Q_1+\cdots+Q_k$, for $Q\in\cal Q$, of polyhedral complex
${\cal Q}$. By Definition~\ref{d:comp}, we get polyhedral
complexes ${\cal Q}_i=\{Q_i\mid Q\in{\cal Q}\}$.  For   $u\in M$,
denote by $\min \langle u, {\cal Q}\rangle $ the minimal value of
$u$ on $|{\cal Q}|$ if it exists, and similarly
 define  $\min \langle u, {\cal Q}_i\rangle$. By Theorem~\ref{t:emb1},
the toric variety $X_\Sigma$ is embedded into the toric variety
$X_{\tilde\Sigma_{\cal Q}}$ as a complete intersection in
homogeneous coordinates. For any $u\in M$, such that $\min \langle
u, {\cal Q}\rangle\ge-1$, $\min \langle u, {\cal Q}_i\rangle=\min
\langle u,   Q_i \rangle$ for some $Q=\sum_{i=0}^k Q_i \in{\cal
Q}$ and all $i$,
  and $\sigma(u)^c\subset \RR_{>0} \cdot |{\cal
Q} |$ for all $\sigma\in\Sigma$, we construct a $k$-parameter
embedded deformation of $X_\Sigma$ in $X_{\tilde\Sigma_{\cal Q}}$
by deforming the defining equations:
 \begin{equation*}
\begin{array}{ccccc}
X_\Sigma & \subset &   {\mathcal X}_{u{\cal Q}{\cal Q}_0\dots{\cal Q}_k } & \subset & X_{\tilde{\Sigma}_{\cal Q}}\times\CC^k\\
\downarrow & & \downarrow& \swarrow  &\\
\{0\}&\subset &\CC^k,&&
\end{array}
\end{equation*}
where the total space ${\mathcal X}_{u{\cal Q}{\cal Q}_0\dots{\cal
Q}_k}$ of the deformation family  in $X_{\tilde{\Sigma}_{\cal
Q}}\times\CC^k$ is a complete intersection    in homogeneous
 coordinates: \begin{equation}\label{e:deformh} \prod_{ \langle e_i^*,v_\xi\rangle>0 }
x_\xi^{ \langle e_i^*,v_\xi\rangle}-\prod_{ \langle
e_i^*,v_\xi\rangle<0 } x_\xi^{ -\langle
e_i^*,v_\xi\rangle}-\lambda_i x^{\tilde u}\prod_{ \langle
e_i^*,v_\xi\rangle<0 } x_\xi^{ -\langle
e_i^*,v_\xi\rangle}=0\end{equation} for $i=1,\dots,k$, with
$\tilde u=u-\sum_{i=1}^k[[\min\langle u,{\cal Q}_i\rangle]]
e_i^*$.

\begin{rem} To see that  the equations are given by polynomials in $\CC[x_\xi\mid \xi\in\tilde{\Sigma}_{\cal Q}(1)]$
 follow the arguments in  Remark~\ref{r:mon}.  Complete intersection and
compatibility   in the affine charts of $X_{\tilde{\Sigma}_{\cal
Q}}$ guarantee the flatness of this family.
\end{rem}

 Consider the polyhedral complex $\Sigma(u)=\{ \sigma(u) \mid \sigma\in\Sigma\}$,
as in Lemma~\ref{l:slice}. Let $\sigma(u)=Q_0+Q_1+\cdots+Q_k$, for
 $\sigma\in \Sigma$,  be a $\Sigma$-separated  Minkowski sum decomposition of this
 complex satisfying condition $(*)$.
 Like in the case of affine toric varieties, for a toric variety $\xs$, we expect that its unobstructed
 infinitesimal deformations are generated by Kodaira-Spencer map
 from the  families ${ \mathcal X}_{u\Sigma(u){\cal Q}_0\dots {\cal
 Q}_k}\rightarrow\CC$ and   will   verify this for compact simplicial toric
 varieties with at worst Gorenstein terminal singularities in the next section.

\begin{ex}
Let $\xs$ be  the weighted projective space $\PP(1,1,2)$ with the
primitive lattice generators   $n_1=(1,0)$, $n_2=(-1,-2)$,
$n_3=(0,1)$  of the rays of the fan:
$$\setlength{\unitlength}{1cm}
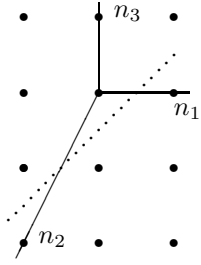
\begin{picture}(6,4.1)
\put( 3,2.7){$n_{1}$} \put( 1.2,1){$n_{2}$} \put( 2.2,4){$n_{3}$}
\put(1,1){\circle*{0.1}} \put(3,3){\circle*{0.1}}
\put(2,2){\circle*{0.1}} \put(2,4){\circle*{0.1}}
\put(3,4){\circle*{0.1}} \put(1,2){\circle*{0.1}}
\put(1,2){\circle*{0.1}}
 \put(1,3){\circle*{0.1}}
\put(1,4){\circle*{0.1}} \put(2,1){\circle*{0.1}}
\put(3,1){\circle*{0.1}} \put(3,2){\circle*{0.1}}
\put(2,3){\circle*{0.1}} \put(2,3){\line(0,1){1.2}}
\put(2,3){\line(1,0){1.2}} \put(2,3){\line(-1,-2){1.1}}
 \multiput(3,3.5)(-0.1,-0.1){23}{\circle*{0.03}}
\end{picture}\label{e:pic}
$$
Here, the dotted line is where the values of  $u=(-2,2)$ are $-1$.
The compact part of the polyhedral complex ${\cal Q}=\Sigma(u)$
obtained by intersecting the cones of the fan of $\xs$ with the
hyperplane $\{n\in N_\RR|\langle
 R,n\rangle=1\}$ consists of just one line segment
 $[\frac{n_1}{2},\frac{n_2}{2}]$ with its vertices.   It admits
 a Minkowski sum decomposition satisfying the condition $(*)$:
$[\frac{n_1}{2},\frac{n_2}{2}]=\frac{n_1}{2}+[0,\frac{n_2-n_1}{2}]$.
The intersections of the cones of the fan adjacent to this line
segment are the vertices $\frac{n_1}{2},\frac{n_2}{2}$ and they
admit the induced decompositions.   According to our construction
the rays of the fan ${\tilde{\Sigma}}={\tilde{\Sigma}}_{\cal Q}$
are spanned by $\frac{n_1}{2}-e_1$, $e_1$,
$\frac{n_2-n_1}{2}+e_1$, and $n_3$ in $ \tilde{N}_\RR$.
Identifying  naturally $\tilde{N}_\RR=N_\RR\oplus\RR^1\simeq
\RR^3$, we get that these points correspond to
$(\frac{1}{2},0,-1)$, $(0,0,1)$, $(-1,-1,1)$ and $(0,1,0)$ in
$\RR^3$. And replacing the first point with $(1,0,-2)$ we get the
primitive lattice generators of the fan ${\tilde{\Sigma}} $. The
maximal cones of the fan ${\tilde{\Sigma}} $ are spanned by all
possible triplets among these lattice points except for the
combination  $(0,1,0)$, $(0,0,1)$,$(-1,-1,1)$. This simply means
that $\xst$ is not compact. If  the homogeneous coordinates
$x_1,x_2,x_3,x_4$ of $X_{\tilde{\Sigma}}$ correspond to
$(1,0,-2)$, $(0,0,1)$, $(-1,-1,1)$ and $(0,1,0)$, then the
equation  defining $\xs$ in $\xst$ is  $x_1^2=x_2x_3$ (exponents
in the monomials correspond to  the ``distances'' of lattice
generators of $\tilde\Sigma$ from $N$). By adding the missing
maximal cone to $\tilde\Sigma$ we get the fan of $\PP^3$ and the
hypersurface $x_1^2=x_2x_3$ simply misses the point
$x_2=x_3=x_4=0$ in $\PP^3$. Using formulas for toric morphisms in
Section~\ref{s:toric}, it is easy to see that the morphism
$\PP(1,1,2)=\xs\rightarrow\xst\subset\PP^3$ is
$(y_1,y_2,y_3)\mapsto(y_1y_2,y_1^2,y_2^2,y_3)$ in homogeneous
coordinates. For $u=(-2,2)$, we get $\tilde u= (-2,2,0)\in
M\oplus\ZZ$ and the corresponding  deformation  is given by
$x_2x_3-x_1^2-\lambda x^2_4=0$ which smoothes out the singularity.
\end{ex}

\begin{ex}
Let $\xs$ be  a toric blow up of the weighted projective space
$\PP(1,1,2)$ with the  lattice generators   $n_1=(1,0)$,
$n_2=(-1,-2)$, $n_3=(0,1)$, $n_4=(0,-1)$  of the rays of the fan:
$$\setlength{\unitlength}{1cm}
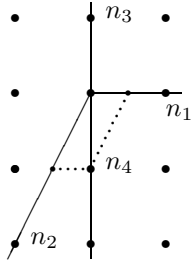
\begin{picture}(8,4.1 )
\put( 3,2.7){$n_{1}$} \put( 1.2,1){$n_{2}$} \put(
2.2,4){$n_{3}$}\put( 2.2,2){$n_{4}$} \put(1,1){\circle*{0.1}}
\put(3,3){\circle*{0.1}} \put(2,2){\circle*{0.1}}
\put(2,4){\circle*{0.1}} \put(3,4){\circle*{0.1}}
\put(1,2){\circle*{0.1}} \put(1,2){\circle*{0.1}}
 \put(1,3){\circle*{0.1}}
\put(1,4){\circle*{0.1}} \put(2,1){\circle*{0.1}}
\put(3,1){\circle*{0.1}} \put(3,2){\circle*{0.1}}
\put(2,3){\circle*{0.1}} \put(2,3){\line(0,1){1.2}}
\put(2,3){\line(1,0){1.2}} \put(2,3){\line(-1,-2){1.1}}
\put(2,3){\line(0,-1){2.2}}
\put(1.5,2){\circle*{0.08}}\put(2.5,3){\circle*{0.08}}
 \multiput(2,2)(-0.1,0){5 }{\circle*{0.03}}\multiput(2,2)( 0.05,0.1){11}{\circle*{0.03}}
\end{picture}
$$
Let ${\cal Q}$ be the polyhedral complex consisting of line
segments $[\frac{n_2}{2}, {n_4} ]$, $ [n_4,\frac{ n_1}{2}]$ with
their vertices. Take a Minkowski sum decomposition of this
complex: $[\frac{n_2}{2}, {n_4} ]=[\frac{n_2}{2}, {n_4} ]+0$,
$[n_4,\frac{ n_1}{2}]=n_4+[0,\frac{n_1}{2}-{n_4} ] $, which
satisfies condition $(*)$ and is $\Sigma$-separated. Then the fan
${\tilde{\Sigma}}_{\cal Q}$ has four 3-dimensional cones,
corresponding to the maximal cones of $\Sigma$, with the lattice
generators    $v_1=n_4-e_1=(0,-1,-1)$, $v_2=e_1=(0,0,1)$,
$v_3=n_3=(0,1,0)$, $v_4=n_1-2n_4+2e_1=(1,2,2)$, $v_5=n_2-2e_1=
(-1,-2,-2)$ in $N_\RR\oplus\RR\cong \RR^3$. By adding two more
cones to ${\tilde{\Sigma}}_{\cal Q}$, the toric variety
$X_{\tilde{\Sigma}_{\cal Q}}$ can be compactified  to
$\PP^2\times\PP^1$ and the embedding $\xs\rightarrow
X_{\tilde{\Sigma}_{\cal Q}}\subset\PP^2\times\PP^1$ is given by
$(y_1,y_2,y_3,y_4)\mapsto(y_1^2 y_4,y_2^2y_4,y_3)\times(y_1,y_2)$,
where $y_1,\dots,y_4$ correspond to $n_1,\dots,n_4$ for $\xs$. If
the homogeneous coordinates $x_1 ,\dots,x_5$ of $\PP^2\times\PP^1$
correspond to the lattice generators of the rays of
$\tilde{\Sigma}_{\cal Q}$ appearing in the same order as above,
then the image of $\xs$ is a hypersurface given by
$x_2x_4^2-x_1x_5^2=0$. For $u=(-1,1)$, we have $\min \langle u,
{\cal Q}\rangle=-1$, $\tilde u=(-1,1,0)$,  and the corresponding
deformation is $x_2x_4^2-x_1x_5^2-\lambda x_3x_4x_5=0$.
\end{ex}

\begin{ex} Our general construction also generalizes the results
of \cite{m1,m2}. To see this let $\xs$ be a toric blow up of a
Fano toric variety corresponding to a reflexive polytope $\Delta$,
such that the rays of $\Sigma$ pass through all lattice points
$n_1,\dots,n_k$ of an edge $\Gamma^*$ of the dual polytope
$\Delta^*$.

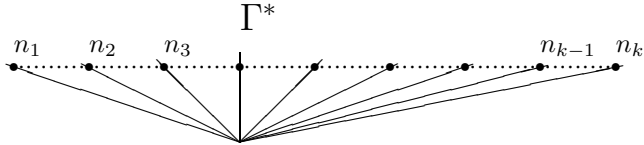
$$ \setlength{\unitlength}{1cm}
\begin{picture}(8,2)
\put(-1,1.2){$n_{1}$} \put(2,0){\line(-3,1){3.1}}
\put(-1,1){\circle*{0.1}} \put(0,1){\circle*{0.1}}
\put(1,1){\circle*{0.1}} \put(2,1){\circle*{0.1}}
\put(3,1){\circle*{0.1}}
\put(4,1){\circle*{0.1}}\put(5,1){\circle*{0.1}}\put(6,1){\circle*{0.1}}\put(7,1){\circle*{0.1}}
 \put(0,1.2){$n_{2}$}
\put(2,0){\line(-2,1){2.1}} \put(2,0){\line(-1,1){1.1}}
\put(1,1.2){$n_3 $} \put(2,0){\line(0,1){1.2}}
\put(2,0){\line(1,1){1.1}} \put(2,0){\line(2,1){2.1}}
\put(2,1.5){\LARGE$\Gamma^*$}
  \put(6,1.2){$n_{ {k-1}}$}
\put(2,0){\line(3,1){3.1}} \put(2,0){\line(4,1){4.1}}
\put(2,0){\line(5,1){5.1}} \put(7,1.2){$n_{k}$}
 \multiput(-1,1)( 0.1,0){80}{\circle*{0.03}}
\end{picture}\label{e:pic}
$$
Let $u\in M$ be in the relative interior of $\Gamma=\{u\in
\Delta\cap M\mid \langle u, \Gamma^*\rangle=-1\}$, the codimension
2 face of $\Delta$ dual to the edge $\Gamma^*$. Then the compact
part of the polyhedral complex ${\cal Q}=\Sigma(u)$   is just the
collection of line segments $[n_i,n_{i+1}]$ together with their
vertices. If we use the Minkowski sum decomposition of the
polyhedral complex as in Example~\ref{ex:line} and apply our
construction of $\tilde \Sigma_{\cal Q}$ we obtain the same fan as
$\Sigma(\Gamma^*)$ in \cite[Section~4]{m2}. The minimal lattice
generators of $\tilde \Sigma_{\cal Q}$  according to our
construction are $n_1-\sum_{i=1}^k e_i $, $n_2-\sum_{i=1}^k e_i$,
$e_j$, $n_2-n_1+e_j$ for $j=1,\dots,k$ which are  up to notation
the same as in that paper.
\end{ex}

\section{Kodaira-Spencer map for deformations of compact simplicial toric varieties.}

The   infinitesimal space $H^1(\xs,{\cal T}_\xs)$  of deformations
of compact smooth toric varieties was computed by N. Ilten in
\cite{i} using Euler exact sequence and Demazure's description of
cohomology of line bundles. Alternatively and more generally, this
can be done for compact simplicial toric varieties $\xs$ using \v
Cech cohomology presentation for cohomology of sheaves ${\cal
O}_\xs(D)$ of torus-invariant Weil divisors $D$ on $\xs$ as in
\cite{torvar}. The last approach is more suitable for us to
describe the Kodaira-Spencer map for deformations of compact
simplicial toric varieties in terms of homogeneous coordinates.

Denote by  ${  \widehat\Omega}^1_X$   the sheaf of {\it Zariski
$1$-forms} on an orbifold $X$ (see Appendix A.3 in \cite{ck}). By
Proposition~A.4.1 in \cite{ck}, the sheaf  ${\cal
T}_X:={\Hom}_{{\cal O}_X}(\widehat\Omega^1_X,{\cal O}_X)$ is
isomorphic to the tangent sheaf on $X$. From Theorem~12.1 in
\cite{bc}, we have the generalized Euler exact sequence
$$0@>>>\widehat\Omega^1_\xs@>>>\bigoplus_{\rho\in\Sigma(1)}{\cal
O}_\xs(-D_\rho)@>>>{\rm Cl}(\xs)\otimes_\ZZ{\cal O}_\xs@>>>0$$ for
compact simplicial toric varieties $\xs$, and its dual
$$0@>>>\Hom_\ZZ({\rm Cl}(\xs),\ZZ)\otimes_\ZZ{\cal O}_\xs@>>>\bigoplus_{\rho\in\Sigma(1)} {\cal
O}_\xs(D_\rho)@>>>{\cal T}_\xs@>>>0,$$ whose  associated long
exact sequence in cohomology gives isomorphism
 \begin{equation}\label{e:isomor}
 H^1(\xs,{\cal T}_\xs)\cong \bigoplus_{\rho\in\Sigma(1)} H^1(\xs,{\cal
 O}_\xs(D_\rho)).\end{equation}

We will use an exact sequence of complexes from the proof of
Theorem~9.1.2 in \cite{torvar} to represent $ H^1(\xs,{\cal
O}_\xs(D_{\rho }))$ by   \v Cech cocycles. Let ${\cal
U}=\{X_{\sigma_j}: j\in J\}$ be the affine open cover of a compact
 toric variety $\xs$, where $ X_{\sigma_j}
=\spe(\CC[\sigma_j^\ve\cap M])$ and $\{\sigma_j:j\in J\}$ is an
ordered set of maximal cones in $\Sigma$.
 The \v Cech cohomology
$\check{H}^1({\cal U},{\cal O}_\xs(D_{\rho}))$
 is calculated  by the \v Cech complex
$$\check{C}^p({\cal U}, {\cal
O}_\xs(D_{\rho}))=\bigoplus_{j_0<\cdots<j_p} H^0(X_{\sigma_{j_0}
\cap\cdots\cap  \sigma_{j_{p}}}, {\cal O}_\xs(D_{\rho}))$$ with
differential
$d^p(\alpha)_{j_0,\dots,j_{p+1}}=\sum_{k=0}^{p+1}(-1)^k\alpha_{j_0,\dots,\widehat{j_k},\dots,j_{p+1}}
|_{X_{\sigma_{j_0} \cap\cdots\cap  \sigma_{j_{p+1}}}}$. The
$M$-grading on
$$H^0(X_{\sigma}, {\cal O}_\xs(D_{\rho}))=\bigoplus_{u\in P_{\rho}^{\sigma}\cap M}\CC\cdot x^u,$$ where
$P_{\rho}^{\sigma}=\{u\in M_\RR|\,\langle u,v_{\rho}\rangle\ge-1,
\langle u,v_{\rho'}\rangle\ge0 \text{ for
}\rho'\in\sigma(1)\setminus\{\rho\}\}$ and
$x^u=\prod_{\rho\in\Sigma(1)} x_\rho^{\langle u,v_\rho\rangle}$,
induces an $M$-grading on the \v Cech complex and cohomlogy:
$$H^p(\xs, {\cal
O}_\xs(D_\rho))\cong\check{H}^p({\cal U}, {\cal
O}_\xs(D_\rho))=\bigoplus_{u\in M}\check{H}^p({\cal U}, {\cal
O}_\xs(D_\rho))_u.$$

Generalizing Proposition~1.1 in \cite{i}, we have

\begin{pr}\label{p:sections} Let $\xs$ be a compact   toric variety.  Then,
 for $\rho \in\Sigma(1)$, we have:\\
 {\rm (i)} $   H^1(\xs,{\cal
O}_\xs(D_{\rho}))_u=0$, if  $\langle u,v_{\rho}\rangle \ne -1 $,\\
{\rm (ii)} $  \dim H^1(\xs,{\cal O}_\xs(D_{\rho}))_u=\max\{0, \dim
H^0(|\Sigma(u)^c|\setminus\{v_\rho\},\CC)-1\}$, if  $\langle
u,v_{\rho}\rangle = -1$, where
$\Sigma(u)^c=\{\sigma(u)^c\mid\sigma\in\Sigma\}$  and
$\sigma(u)^c$ is the convex hull of vertices of
$\sigma(u)=\sigma\cap\{n\in
N_\RR\mid\langle u,n\rangle=-1\}$,\\
 {\rm (iii)} For $u\in M$,
 such that $\langle u,v_{\rho}\rangle=-1 $, let $C_0,\dots,C_k$ be the connected components of
$|\Sigma(u)^c|\setminus\{v_{\rho}\}$, then the space
 $ H^1(\xs,{\cal O}_\xs(D_{\rho}))_u\cong
 \check{H}^1({\cal U},{\cal O}_\xs(D_{\rho}))_u$  is spanned by
the \v Cech cocycles $ s^i_{u,\rho}:=
 \bigl\{(\delta_{j_1}^i-\delta_{j_0}^i)x^u \bigr\}_{j_0j_1}$, for
$i=0,\dots,k$, where $\delta_{j}^i=1$
 if $\sigma_j\cap C_i\ne\emptyset $ and $\delta_{j}^i=0$ if  $\sigma_j\cap C_i=\emptyset
 $. Moreover,   $ \sum_{i=0}^ks^i_{u,\rho}=0$ and
 $\{s^i_{u,\rho}\mid i=1,\dots,k\}$ is a basis of $\check{H}^1({\cal U},{\cal
 O}_\xs(D_{\rho}))_u$.
\end{pr}

\begin{pf}
Adapting the proof of Theorem 9.1.2 in \cite{torvar}, note:
$$x^u\in H^0(X_\sigma, {\cal
O}_\xs(D_{\rho}))\Longleftrightarrow u\in P_{\rho}^{\sigma}\cap
M\Longleftrightarrow\sigma\cap
|\Sigma(u)^c|\setminus\{v_{\rho}\}=\emptyset.$$ Since
$\sigma\cap|\Sigma(u)^c|\setminus\{v_{\rho}\}=\sigma(u)^c\setminus\{v_{\rho}\}$
is convex, this defines a canonical exact sequence
$$0@>>>H^0(X_\sigma, {\cal
O}_\xs(D_{\rho}))_u@>>>\CC@>>>H^0(\sigma\cap|\Sigma(u)^c|\setminus\{v_{\rho}\},\CC)@>>>0$$
and  an exact sequence of complexes
$$0@>>>\check{C}^\bullet({\cal U}, {\cal O}_\xs(D_\rho))_u@>>> B^\bullet@>>> C^\bullet@>>>0,$$ where
$$B^p=\bigoplus_{j_0<\cdots<j_p} \CC,\quad C^p=\bigoplus_{j_0<\cdots<j_p} H^0(\sigma_{j_0}\cap\cdots\cap\sigma_{j_p}\cap|\Sigma(u)^c|\setminus\{v_{\rho}\},\CC)$$
with differentials induced from the differential of the \v Cech
complex. By (9.1.10)  in \cite{torvar}, the vanishing
$H^1(B^\bullet)=0$ implies the exact sequence of
 cohomology:
$$H^0(\xs, {\cal
O}_\xs(D_{\rho}))_u@>>>\CC@>>>H^0(
|\Sigma(u)^c|\setminus\{v_{\rho}\},\CC)@>>>H^1(\xs, {\cal
O}_\xs(D_{\rho}))_u@>>>0.$$

 If  $\langle u,v_{\rho}\rangle \ne -1
$, then $v_{\rho}\not\in|\Sigma(u)^c|$ and
$|\Sigma(u)^c|\setminus\{v_{\rho}\}=|\Sigma(u)^c|$ is connected by
the completeness of the fan $\Sigma$. Also, note $H^0(\xs, {\cal
O}_\xs(D_{\rho}))_u\ne0$ if and only if
$|\Sigma(u)^c|\setminus\{v_{\rho}\}=\emptyset$. Hence, the above
exact sequence implies the vanishing $H^1(\xs, {\cal
O}_\xs(D_{\rho}))_u=0$, if $\langle u,v_{\rho}\rangle \ne -1 $.
 If
$\langle u,v_{\rho}\rangle = -1 $ and
$|\Sigma(u)^c|\setminus\{v_{\rho}\}\ne\emptyset$, then the
dimension formula for  $H^1(\xs, {\cal O}_\xs(D_{\rho}))_u$  in
part ({\rm ii}) follows as well.

To compute the \v Cech cocyles representing $H^1(\xs, {\cal
O}_\xs(D_{\rho}))_u$ from the connecting homomorphism
$H^0(C^\bullet)\rightarrow \check{H}^1({\cal U},{\cal
O}_\xs(D_{\rho}))_u$, consider the commutative diagram:
$$ \begin{CD}
\check{C}^1({\cal U}, {\cal
O}_\xs(D_{\rho}))_u@>>>  B^1@>>> C^1@>>>0  \\
    @AAA@AAA @AAA\\
\check{C}^0({\cal U}, {\cal O}_\xs(D_{\rho}))_u @>>>B^0@>>> C^0
@>>> 0
\end{CD}   $$
The basis of  $H^0( |\Sigma(u)^c|\setminus\{v_{\rho}\},\CC)$ is
represented  by the  cocycles  $\{c_{j_0}^i\}_{j_0}\in C^0$, where
$c_{j_0}^i\in
H^0(\sigma_{j_0}\cap|\Sigma(u)^c|\setminus\{v_{\rho}\},\CC)$ is
defined by $c_{j_0}^i(\sigma_{j_0}\cap C_i)=1$ and
$c_{j_0}^i(\sigma_{j_0}\cap C_j)=0$ for $j\ne i$. This cocycle
lifts to
 $\{\delta_{j_0}^i\}_{j_0}\in B^0=\bigoplus_{j_0}\CC$. Applying the
 differential $B^0\rightarrow B^1$ we get
 $\{\delta_{j_1}^i-\delta_{j_1}^i\}_{j_0j_1}\in B^1$ and by
 commutativity of the diagram this comes from the \v Cech cocycle  $ s^i_{u,\rho}=
 \bigl\{(\delta_{j_1}^i-\delta_{j_0}^i)x^u \bigr\}_{j_0j_1}$.

The relation  $ \sum_{i=0}^ks^i_{u,\rho}
\bigl\{\bigl(\sum_{i=0}^k\delta_{j_1}^i-\sum_{i=0}^k\delta_{j_0}^i\bigr)x^u
\bigr\}_{j_0j_1}= \{ (  1-1 )x^u  \}_{j_0j_1}=0$ follows from the
property that every cone $\sigma\in\Sigma$ intersects precisely
one connected component of $|\Sigma(u)^c|\setminus\{v_{\rho}\}$.
Hence, $\{s^i_{u,\rho}\mid i=1,\dots,k\}$ is a basis by the
dimension count in part (ii).
\end{pf}

Similar to Proposition 1.4 in \cite{i}, we have

\begin{pr}\label{p:dimH1}  Let $\xs$ be a compact simplicial toric variety.
Then, for $u\in M$,\\
{\rm(i)}  $\ds\dim H^1(\xs, {\cal
T}_\xs)_u=\sum_{\begin{Sb}\rho\in\Sigma(1), \\ \langle
u,v_\rho\rangle=-1\end{Sb}} \max\{0, \dim
H^0(|\Sigma(u)^c|\setminus\{v_\rho\},\CC)-1\}$,\\
 {\rm(ii)}   The space
 $ H^1(\xs, {\cal
T}_\xs)_u\cong
 \check{H}^1({\cal U}, {\cal
T}_\xs)_u$ is spanned by   the \v Cech cocycles $
\gamma^i_{u,\rho}:=
 \bigl\{(\delta_{j_1}^i-\delta_{j_0}^i)x^u  x_\rho\partial/\partial
 x_\rho\bigr\}_{j_0j_1}$,
for $\rho\in\Sigma(1)$,    such that  $ \langle
u,v_\rho\rangle=-1$, and connected components $C_i$ of
$|\Sigma(u)^c|\setminus\{v_\rho\}$. Moreover,   $
\sum_{i=0}^{k_{u,\rho}} \gamma^i_{u,\rho}=0$, where
$k_{u,\rho}=\dim H^0(|\Sigma(u)^c|\setminus\{v_\rho\},\CC)$ and
 $$\{\gamma^i_{u,\rho}\mid\rho\in\Sigma(1), \langle
u,v_\rho\rangle=-1; i=1,\dots,k_{u,\rho} \}$$ is a basis of
$\check{H}^1({\cal U},{\cal
 T}_\xs )_u$.
\end{pr}

\begin{pf} This is an immediate consequence of Proposition~\ref{p:sections}
since the isomorphism  $(\ref{e:isomor})$ is induced from the
homomorphisms $H^0(X_\sigma, {\cal O}_\xs(D_{\rho}))\rightarrow
H^0(X_\sigma, {\cal T}_\xs )$ given by $s\mapsto s\cdot
x_\rho\partial/\partial x_\rho$.
\end{pf}

For partial crepant   resolutions of Gorenstein Fano toric
varieties we get the following simple formula.

\begin{thm}\label{t:infdef} Let $\xs$ be a compact simplicial toric variety such
that its anticanonical divisor $-K_\xs=\sum_{\rho\in\Sigma(1)}
D_\rho$ is semiample, and let $$\Delta =\{m\in M_\RR\mid\langle
m,v_\rho\rangle\ge-1\,\forall\, \rho\in\Sigma(1)\}$$ be the
associated reflexive polytope with its dual
$\Delta^*=\co(v_\rho|\,\rho\in\Sigma(1))$.  Then $$\ds\dim
H^1(\xs, {\cal T}_\xs)=\sum_{\begin{Sb} \Gamma\prec\Delta \\{\rm
codim} \Gamma=2
\end{Sb}}l^*(\Gamma)l^*_\Sigma(\Gamma^*),$$
where   $\Gamma^*$ is the 1-dimensional face of $\Delta^*$ dual to
$\Gamma$, $l^* (\Gamma)$ is the number of interior lattice points
in the face $\Gamma$,   and $l^*_\Sigma(\Gamma^*)$ is the number
of $v_\rho$, for $\rho\in\Sigma(1)$, in the relative interior of
$\Gamma^*$. In particular, if $(\Delta^*\setminus\{0\})\cap
N=\{v_\rho\mid\rho\in\Sigma(1)\}$, then $\xs$ has at worst
Gorenstein  terminal singularities and the dimension of its
infinitesimal space of deformations is $$\dim H^1(\xs, {\cal
T}_\xs)=\sum_{\begin{Sb} \Gamma\prec\Delta \\{\rm codim} \Gamma=2
\end{Sb}}l^*(\Gamma)l^*(\Gamma^*).$$
\end{thm}

\begin{pf} Take $v_{\rho_0}\in\Delta^*$ corresponding to
$\rho_0\in\Sigma(1)$ and $u\in M$ such that $\langle
u,v_{\rho_0}\rangle=-1$. Then note that $ \Sigma(u)^c  $ is a
deformation retract of the simplical complex on the boundary of
$\Delta^*$ consisting of $\co(v_\rho\mid\rho\in\sigma(1),\langle
u,v_\rho\rangle\le-1)$, for $\sigma\in\Sigma$. Points in such
simplices are connected by pathes along the boundary of $\Delta^*$
to the face of $\Delta^*$ where $u$ takes minimal value. Hence,
the removal of the point $v_{\rho_0}$, where $u$ takes the maximal
value  on the underlying set of the simplicial complex, will not
have effect on the connectedness of $|\Sigma(u)^c|$, unless
$|\Sigma(u)^c|=\{n\in \Delta^*\mid\langle u,n\rangle=-1\}$ is a
face $\Gamma^*\prec\Delta^*$, such that $\dim\Gamma^*=1$ and
$v_{\rho_0}$ is in the relative interior of $\Gamma^*$. In this
case $u$ has to be in the relative interior of the dual face
$\Gamma=\{m\in \Delta\mid\langle
m,n\rangle=-1\,\forall\,n\in\Gamma^*\}$. Now, the  formula follows
from Proposition~\ref{p:dimH1}. The other part follows from
Lemma~4.1.2 and Proposition~A.4.2 in \cite{ck}.
\end{pf}

\begin{rem} By Corollary~2.2 in \cite{m2},  the above formula shows that the dimension of
the infinitesimal space of ``non-polynomial'' deformations of a
semiample Calabi-Yau hypersurface in a compact simplicial toric
variety $\xs$ with at worst terminal singularities coincides with
the dimension of $H^1(\xs,  {\cal T}_\xs)$. In particular,
deformations of such toric varieties in \cite{m1,m2}, which induce
deformations of Calabi-Yau hypersurfaces,  span $H^1(\xs, {\cal
T}_\xs)$ by the Kodaira-Spencer map.
\end{rem}

Let $\xs$ be a compact simplicial toric variety, and
$\rho_0\in\Sigma(1)$,  $u\in M$, such that $\langle
u,v_{\rho_0}\rangle=-1$. Let $k=\dim
H^0(|\Sigma(u)^c|\setminus\{v_{\rho_0}\},\CC)-1$, and let
$C_0,C_1,\dots,C_k$ be the connected components of
$|\Sigma(u)^c|\setminus\{v_{\rho_0}\}$. Consider  the following
Minkowski sum decomposition of the polyhedral complex
$\Sigma(u)^c$: $\sigma(u)^c=Q_0+Q_1+\cdots+Q_k$, where
\begin{equation}\label{e:simpdec} Q_0=\Biggl\{\begin{array}{ll }  \sigma(u)^c  &\text{if } \sigma(u)^c\cap C_0 \ne\emptyset, \\
 v_{\rho_0}&  \text{if } \sigma(u)^c\cap C_0 =\emptyset,\end{array}\,
  Q_i=\Biggl\{\begin{array}{cc } 0  &\text{if } \sigma(u)^c\cap C_i=\emptyset, \\
\sigma(u)^c- v_{\rho_0}&  \text{if } \sigma(u)^c\cap C_i
\ne\emptyset,\end{array}\end{equation}   for $i\ne0$. Since every
cone $\sigma\in \Sigma$  intersects at most one connected
component of $|\Sigma(u)^c|\setminus\{v_{\rho_0}\}$, this
construction is well defined and   determines a Minkowski sum
decomposition of $\Sigma(u)$  by adding the recession cone  of
$\sigma(u)$ to each summand. By Lemma~\ref{l:slice}, this
Minkowski sum decomposition is $\Sigma$-separated and we get the
fan $\tilde\Sigma_{u,\rho_0}:=\tilde\Sigma_{\cal Q}$ for ${\cal
Q}=\Sigma(u)$.

\begin{ex} A slice of a simplicial fan  consists of polyhedra
whose compact parts are faces of the polyhedra. In the following
picture we have a slice of  a complete simplicial fan $\Sigma$ in
$\RR^3$ in the hyperplane, where $u=(0,0,-1)$ has value $-1$, and
its Minkowski sum decomposition. In this case
$\Sigma(u)^c\setminus\{(3,2,1)\}$ consists of three connected
components corresponding to two triangles and one line segment.

$$ \setlength{\unitlength}{0.6cm}
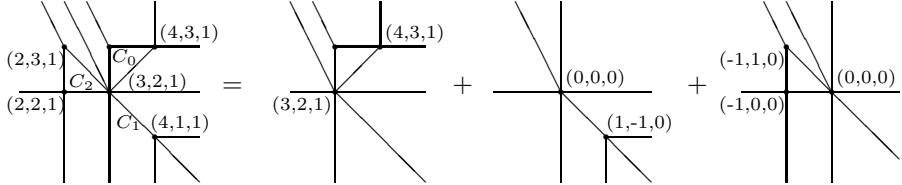
\begin{picture}(30,4)
 \put(2,2){\line(-1,0){1}} \put(2,2){\line(0,-1 ){2}}\put(3,3){\line(0,-1 ){3}}
\put(2,2){\line(0,1){1}}\put(2,2){\line(1,0){3}}
\put(4,3){\line(0,1){1}}\put(4,1){\line(0,-1){1}}
\put(3,3){\line(1,0){2}}\put(4,1){\line(1,0){1}}
\put(2,3){\line(1,-1){3}}\put(3,2){\line(1, 1){1}}
 \put(3,3){\line(-1,2){0.5}} \put(2,3){\line(-1,2){0.5}}\put(2,4){\line(1,-2){1}}

\put(3.05,2.65){\SMALL{$C_0$}} \put(2.1,2.1){\SMALL{$C_2$}}
\put(3.15,1.2){\SMALL{$C_1$}}

 \put(2,3){\circle*{0.1}}
\put(0.7,2.6){\SMALL{(2,3,1)}}\put(3,3){\circle*{0.1}}
\put(4,3){\circle*{0.1}}\put(4.1,3.2){\SMALL{(4,3,1)}}
 \put(2,2){\circle*{0.1}}\put(0.7,1.6){\SMALL{(2,2,1)}}
   \put(3,2){\circle*{0.1}}\put(3.4,2.1){\SMALL{(3,2,1)}}
\put(4,1){\circle*{0.1}}\put(3.9,1.2){\SMALL{(4,1,1)}}

\put( 5.5,2 ){$=$}

\put(8,2){\line(-1,0){1}} \put(8,2){\line(0,-1 ){2}}
\put(8,2){\line(0,1){1}}\put(8,2){\line(1,0){2}}
\put(8,2){\line(1, 1){1}} \put(9,3){\line(0,1){1}}
 \put(8,3){\line(1,0){2}}\put(8,3){\line(-1,2){0.5}}\put(8,2){\line(-1,2){1}}\put(8,2){\line( 1,-1){2}}

\put(8,2){\circle*{0.1}} \put(6.6 ,1.6){\SMALL{(3,2,1)}}
\put(8,3){\circle*{0.1}} \put(9,3){\circle*{0.1}}
\put(9.1,3.2){\SMALL{(4,3,1)}}

 \put(
10.6,2 ){$+$}

 \put(11.5,2){\line( 1,0){3.5}}  \put(13,4){\line(0,-1 ){4}}
 \put(14,1){\line(1,0){1}}\put(14,1){\line(0,-1){1}}
\put(13,2){\line(1,-1){2}}\put(13,2){\line(-1,2){1}}

\put(13,2){\circle*{0.1}} \put(13.1,2.2){\SMALL{(0,0,0)}}
\put(14,1){\circle*{0.1}} \put(14,1.2){\SMALL{(1,-1,0)}}

\put( 15.8,2 ){$+$}

\put(18,2){\line(-1,0){1}} \put(18,2){\line(0,-1
){2}}\put(19,4){\line(0,-1 ){4}}
\put(18,2){\line(0,1){1}}\put(18,2){\line(1,0){2.5}}
\put(18,3){\line(1,-1){2.5}}
  \put(18,3){\line(-1,2){0.5}}\put(19,2){\line(-1,2){1}}

\put(18,3){\circle*{0.1}} \put(16.5,2.6){\SMALL{(-1,1,0)}}

\put(19,2){\circle*{0.1}} \put(19.1,2.2){\SMALL{(0,0,0)}}

\put(18,2){\circle*{0.1}}  \put(16.5,1.6){\SMALL{(-1,0,0)}}
\end{picture}
$$
\end{ex}

The 1-dimensional cones of the fan $\tilde\Sigma_{u,\rho_0}$ are
$$\xi_0=\RR_{\ge0}\Biggl(v_{\rho_0}-\sum_{i=1}^k e_i\Biggr), \,\, \xi_i=\RR_{\ge0}e_i,\text{ for }
i=1,\dots,k, \, \text{ and }$$
 $$\xi_\rho= \Biggl\{\begin{array}{ll } \RR_{\ge0} ( (v_{\rho}/ \langle -u,v_\rho \rangle)-\sum_{i=1}^k
 e_i),
&\text{for }\rho\cap C_0 \ne\emptyset,\\
\RR_{\ge0}((v_{\rho}/\langle -u,v_\rho \rangle)-v_{\rho_0}+
e_i),&\text{for }\rho\cap C_i \ne\emptyset, i\ne0,\\
\RR_{\ge0} v_{\rho}, & \text{for }\rho\cap |\Sigma(u)|=\emptyset.
\end{array}$$
The corresponding primitive lattice generators are
$v_{\xi_0}=v_{\rho_0}-\sum_{i=1}^k e_i$, $v_{\xi_i}=e_i$, for
$i=1,\dots,k$, and $$v_{\xi_\rho}=\Biggl\{\begin{array}{ll }
 v_{\rho} + \langle u,v_\rho \rangle \sum_{i=1}^k
 e_i,
&\text{for }\rho\cap C_0 \ne\emptyset,\\
 v_{\rho}+\langle u,v_\rho \rangle(v_{\rho_0}-
e_i),&\text{for }\rho\cap C_i \ne\emptyset, i\ne0,\\
 v_{\rho}, & \text{for }\rho\cap
|\Sigma(u)|=\emptyset.\end{array}$$

Following the formulas in Section~\ref{s:toric}, we will describe
the toric morphism $\phi:\xs\rightarrow
X_{\tilde\Sigma_{u,\rho_0}}$, associated to the map of  $\Sigma$
to $\tilde{\Sigma}_u$, in terms of homogeneous coordinates. The
lattice generators of $\Sigma$ are linear combinations of the
generators of minimal cones of $\tilde\Sigma_{u,\rho_0}$
containing them:
$v_{\rho_0}=v_{\xi_0}+v_{\xi_1}+\cdots+v_{\xi_k}$,
$$v_{ \rho}=\Biggl\{\begin{array}{ll }
 v_{\xi_\rho}- \langle u,v_\rho \rangle \sum_{j\ne i}
 v_{\xi_j},
&\text{for }\rho\cap C_i \ne\emptyset,\\
v_{\xi_\rho}, & \text{for }\rho\cap
|\Sigma(u)|=\emptyset.\end{array}$$
 Then by (\ref{e:morphism}), the morphism  $\phi:\xs\rightarrow
 X_{\tilde\Sigma_{u,\rho_0}}$ is given by $(x_\rho)_{\rho\in\Sigma(1)}\mapsto (x_\xi)_{\xi\in\tilde\Sigma_{u,\rho_0}(1)} $
with $x_{\xi}= x_{\rho}$, for  $\xi=\xi_\rho$, $\rho\ne\rho_0$,
and
$$x_{\xi}=
 \prod_{\rho\cap
|\Sigma(u)|\setminus C_i\ne\emptyset}  x_{\rho}^{-\langle u,v_\rho
\rangle},
 \text{ for } \xi=\xi_i, \,i=0,\dots,k. $$

By Theorem~\ref{t:emb1}, the image of $\phi$ is a complete
intersection in $X_{\tilde\Sigma_{u,\rho_0}}$, given by
\begin{equation*}  x_{\xi_i}\prod_{\rho\cap C_i \ne\emptyset }
x_{\xi_\rho}^{ -\langle u,v_\rho\rangle }-x_{\xi_0}\prod_{\rho\cap
C_0 \ne\emptyset } x_{\xi_\rho}^{ -\langle u,v_\rho\rangle
}=0,\end{equation*} for $i=1,\dots,k$. The $k$-parameter
deformation ${ \mathcal X}_{u\Sigma(u){\cal Q}_0\dots {\cal
 Q}_k}$ is given by the complete intersection
\begin{equation}\label{e:deformc}  x_{\xi_i}\prod_{\rho\cap C_i \ne\emptyset }
x_{\xi_\rho}^{ -\langle u,v_\rho\rangle }-x_{\xi_0}\prod_{\rho\cap
C_0 \ne\emptyset } x_{\xi_\rho}^{ -\langle u,v_\rho\rangle
}-\lambda_i\prod_{\rho\cap |\Sigma(u)|=\emptyset} x_{\xi_\rho}^{
\langle u,v_\rho\rangle }=0,\end{equation} for $i=1,\dots,k$, in
$X_{\tilde\Sigma_{u,\rho_0}}\times\CC^k$ by the construction in
Section~\ref{s:darb}, where $\tilde u=u$
  since $ \langle u,Q_i\rangle=0$, for $i=1,\dots,k$.

\begin{thm}\label{t:koda} Let $\xs$ be a compact simplicial toric variety, $\rho_0\in\Sigma(1)$, $u\in M$, such that $\langle
u,v_{\rho_0}\rangle=-1$,  and let $\sigma(u)=Q_0+Q_1+\cdots+Q_k$,
for $\sigma\in\Sigma$, be a Minkowski sum decomposition of
$\Sigma(u)$ corresponding to $(\ref{e:simpdec})$. For the family
${\mathcal X}_{u\Sigma(u){\cal Q}_0\dots {\cal
Q}_k}\rightarrow\CC^k$,
 the Kodaira-Spencer map $\kappa_{u,\rho_0} :T_{\CC^k,0}\rightarrow
 T^1_{X_\Sigma}$
 sends
the basis vector $\partial/\partial \lambda_i$ to the element in
$H^1(\xs, {\cal T}_\xs)\subseteq T^1_{X_\Sigma}$ represented by
the \v Cech cocycle
$$\Biggl\{ \Biggl(\prod_{\rho\in\Sigma(1)} x_\rho^{\langle
u,v_\rho\rangle}\Biggr)x_{\rho_0}\frac{\partial}{\partial
x_{\rho_0}}\Biggr\}_{U^i_0U^i_1},$$ where
$$U^i_0=\Biggl\{(x_\rho)\in\xs:\prod_{\rho\cap C_i \ne\emptyset }
x_{\rho}\ne0\Biggl\},\quad
 U^i_1=\Biggl\{(x_\rho)\in\xs: \prod_{\rho\cap
 (\cup_{j\ne i} C_j)\ne\emptyset } x_{\rho}\ne0\Biggl\}.$$
\end{thm}

\begin{pf}  We will use the techniques of \cite{m1,m2}. Note that the two sets $U^i_0$ and $U^i_1$ form an open cover of $\xs$, since
the homogeneous coordinates $x_{\rho_1}$,  $x_{\rho_2}$ do not
vanish simultaneously if  $\rho_1\cap C_i \ne\emptyset$,
$\rho_2\cap C_j \ne\emptyset$ for $i\ne j$. Similar to the proof
of Theorem 4.5 in \cite{m2}, consider a one parameter subfamily of
${\mathcal X}_{u\Sigma(u){\cal Q}_0\dots {\cal Q}_k}$ consisting
of the complete intersection
$${\mathcal X}_i=\{((x_\xi),\lambda)\in X_{\tilde\Sigma_{u,\rho_0}}\times
\AA^1:  ((x_\xi),\lambda)\text{ satisfies } (\ref{e:deformc} )
\text{ with } \lambda_i=\lambda,\lambda_j=0 \text{ for }j\ne i\}
$$  in $X_{\tilde\Sigma_{u,\rho_0}}\times
\AA^1$.
 Then the   sets
$\ds \tilde U^i_0=\biggl\{(x_\xi)\in X_{\tilde\Sigma_{u,\rho_0}}:
\prod_{\rho\cap C_i \ne\emptyset } x_{\xi_\rho}\ne0\biggr\}$ and
$\ds\tilde U^i_1=\biggl\{(x_\xi)\in X_{\tilde\Sigma_{u,\rho_0}} :
\prod_{\rho\cap (\cup_{j\ne i} C_j)\ne\emptyset }
x_{\xi_\rho}\ne0\biggr\} $ form an open cover of
$X_{\tilde\Sigma_{u,\rho_0}}$ since the 1-dimensional cones $
\xi_{\rho_1},\xi_{\rho_2}\in\tilde\Sigma_{u,\rho_0}(1)$  with
$\rho_1\cap C_i \ne\emptyset$, $\rho_2\cap C_j \ne\emptyset$,
$i\ne j$, cannot lie in the same cone of the fan
$\tilde\Sigma_{u,\rho_0}$ by its construction.

We have the following isomorphisms:
$$\varphi_j:U^i_j\times\AA^1\stackrel{\phi\times{\rm
id}}\cong(\phi(\xs)\cap \tilde U^i_j)\times\AA^1\stackrel{
\psi_j}\cong  {\mathcal X}_i\cap (\tilde U^i_j\times\AA^1), \text{
for } j=0,1,$$ where  $\psi_0$ is induced by the isomorphism
$\tilde U^i_0\times\AA^1 \cong \tilde U^i_0\times\AA^1 $ which
changes only the coordinate
$$x_{\xi_i}\mapsto x_{\xi_i}-\lambda\prod_{\rho\cap C_i
\ne\emptyset } x_{\xi_\rho}^{ \langle u,v_\rho\rangle
}\prod_{\rho\cap |\Sigma(u)|=\emptyset} x_{\xi_\rho}^{ \langle
u,v_\rho\rangle },$$ and $\psi_1$ is induced by the isomorphism
$\tilde U^i_1\times\AA^1 \cong \tilde U^i_1\times\AA^1 $ which
changes only the coordinates
$$x_{\xi_j}\mapsto x_{\xi_j}+\lambda\prod_{\rho\cap C_j
\ne\emptyset } x_{\xi_\rho}^{ \langle u,v_\rho\rangle
}\prod_{\rho\cap |\Sigma(u)|=\emptyset} x_{\xi_\rho}^{ \langle
u,v_\rho\rangle }, \text{ for } j\ne i.$$

 It is not difficult to
compute  that
$$\varphi_0^{-1}\varphi_1:(U_0^i\cap U_1^i)\times\AA^1\cong(U_0^i\cap
U_1^i)\times\AA^1$$ changes only the coordinate $x_{\rho_0}\mapsto
x_{\rho_0}+\lambda x_{\rho_0}\prod_{\rho\in\Sigma(1)}
x_\rho^{\langle u,v_\rho\rangle}$. Hence, like in \cite{m1}, the
one parameter deformation   ${\mathcal X}_i\rightarrow \AA^1$ of
$\xs$ is locally trivial and corresponds to the above regluing of
open toric  subvarieties $U_0^i$ and $U_1^i$ of $\xs$  by the
automorphism on their intersection.

 To compute the corresponding infinitesimal deformation change the base $\AA^1$ to
 $\spe(\CC[\varepsilon]/(\varepsilon^2))$ and restrict to the affine cover $X_\sigma=\spe(\CC[\sigma^\ve\cap M])$. Then
 $\varphi_0^{-1}\varphi_1$ on the level of rings corresponds to
 the ring  isomorphisms
 $$\varphi_{\sigma}^*:S(\Sigma)_{x^{\hat\sigma}g}^{G(\Sigma)}\otimes\CC[\varepsilon]/(\varepsilon^2)\rightarrow
S(\Sigma)_{x^{\hat\sigma}g}^{G(\Sigma)}\otimes\CC[\varepsilon]/(\varepsilon^2),
$$ where $S(\Sigma)=\CC[x_\rho\mid\rho\in\Sigma(1)]$ is the homogeneous
coordinate ring of $\xs$,   $$x^{\hat\sigma}= \prod_{\rho
 \not\subseteq\sigma} x_\rho, \quad g=\prod_{\rho\cap
|\Sigma(u)|\setminus \rho_0\ne\emptyset}  x_{\rho}.$$ By the
coordinate change  $x_{\rho_0}\mapsto x_{\rho_0}+\varepsilon
x_{\rho_0}\prod_{\rho\in\Sigma(1)} x_\rho^{\langle
u,v_\rho\rangle}$,
\begin{multline*}\varphi_{\sigma}^*\Biggl(\prod_{\rho\in\Sigma(1)}
x_\rho^{\langle
a,v_\rho\rangle}+\varepsilon\prod_{\rho\in\Sigma(1)}
x_\rho^{\langle b,v_\rho\rangle}\Biggr)=
\\ =\prod_{\rho\in\Sigma(1)} x_\rho^{\langle
a,v_\rho\rangle}+\varepsilon\prod_{\rho\in\Sigma(1)}
x_\rho^{\langle b,v_\rho\rangle} +\varepsilon\langle
a,v_{\rho_0}\rangle\prod_{\rho\in\Sigma(1)} x_\rho^{\langle
a,v_\rho\rangle}\prod_{\rho\in\Sigma(1)} x_\rho^{\langle
u,v_\rho\rangle}=\\ =\prod_{\rho\in\Sigma(1)} x_\rho^{\langle
a,v_\rho\rangle}+\varepsilon\prod_{\rho\in\Sigma(1)}
x_\rho^{\langle b,v_\rho\rangle}+
\varepsilon\Biggl(\prod_{\rho\in\Sigma(1)} x_\rho^{\langle
u,v_\rho\rangle}\Biggr)x_{\rho_0}\frac{\partial}{\partial
x_{\rho_0}}\Biggl(\prod_{\rho\in\Sigma(1)} x_\rho^{\langle
a,v_\rho\rangle}\Biggr).\end{multline*}
 Hence, we get the \v Cech cocycle representing our infinitesimal
 deformation.
\end{pf}

\begin{rem} The affine open cover ${\cal
U}=\{X_{\sigma_j}: j\in J\}$ is a refinement of the open cover
$\{U^i_0,U^i_1\}$ of $\xs$, since
$X_{\sigma_j}=\{(x_\rho)\in\xs:\prod_{\rho\not\subseteq\sigma_j }
x_{\rho}\ne0\}\subseteq U^i_0$ if $\sigma_j\cap C_i=\emptyset$,
and $X_{\sigma_j} \subseteq U^i_1$ if $\sigma_j\cap
C_i\ne\emptyset$. One can easily check that the \v Cech cocycle in
Theorem~\ref{t:koda} and the  \v Cech cocycle $
\gamma^i_{u,{\rho_0}} =
 \bigl\{(\delta_{j_1}^i-\delta_{j_0}^i)x^u  x_{\rho_0}\partial/\partial
 x_{\rho_0}\bigr\}_{j_0j_1}$ in Proposition~\ref{p:dimH1} represent the same element in $H^1(\xs, {\cal
 T}_\xs)$. Since these cocycles form a basis, the Kodaira-Spencer
 map $\kappa_{u,\rho_0}$ is injective.
\end{rem}

Combining Theorem~\ref{t:koda} with Proposition~\ref{p:dimH1} and
the last remark, we have one of our main results:

\begin{thm} Let $\xs$ be a compact simplicial toric variety. Then
$$H^1(\xs, {\cal T}_\xs)=\bigoplus_{\begin{Sb}\rho\in\Sigma(1),u\in M \\ \langle
u,v_\rho\rangle=-1\end{Sb}} {\rm Im} ( \kappa_{u,\rho} ),$$ where
$\kappa_{u,\rho}$   are the Kodaira-Spencer maps for the families
of complete intersections in
$X_{\tilde\Sigma_{u,\rho}}\times\CC^k$ given by (\ref{e:deformc}).
\end{thm}

In particular,  our deformation families span the infinitesimal
space of deformations $T^1_{X_\Sigma}=H^1(\xs, {\cal T}_\xs)$ for
compact simplicial toric varieties $\xs$ with at worst Gorenstein
terminal singularities.

 \begin{rem} Using a  combinatorial construction similar to our Minkowski sum decompositions of polyhedral complexes,
  locally trivial deformations of   compact smooth toric varieties
were constructed in the language of T-varieties in \cite{iv}.
 \end{rem}

\section{Deformations of Gorenstein Fano toric varieties.}

In this section we will apply our general construction in
Section~\ref{s:darb} to combine  deformation  families in the case
of a Gorenstein Fano toric variety.

A Gorenstein Fano toric variety corresponds to a reflexive
polytope $\Delta$ and has a (normal) fan $\Sigma_\Delta$ whose
cones are generated by the proper faces of the dual polytope
$$\Delta^*=\{n\in N_\RR\mid \langle m,n\rangle\ge-1\, \forall\, m\in \Delta\}.$$

 Let $\cal
Q=\{\Gamma^*\mid\Gamma^*\prec\Delta^*\}$ be the polyhedral complex
consisting of the proper faces of the polytope $\Delta^*$.
 To obtain a Minkowski sum
decomposition of this polyhedral complex we can simply use a
Minkowski sum decomposition of the polytope
\begin{equation}\label{e:deldec}\Delta^*=\Delta_0^*+\Delta_1^*+\cdots+\Delta_k^*
\end{equation}
 by lattice polytopes. This induces Minkowski sum
decompositions for each face
$$\Gamma^*=\Gamma_0^*+\Gamma_1^*+\cdots+\Gamma_k^*.$$

Let $\tilde{N}=N\oplus\ZZ^k$ and construct the cones
$$\tilde{\sigma}_{\Gamma^*}=\langle
  \Gamma_0^*-e_1-\dots-e_k, \Gamma_1^*+e_1,\dots,\Gamma_k^* +e_k\rangle$$
in $\tilde{N}_\RR$
 for each $\sigma_{\Gamma^*}\in\Sigma_\Delta$, where $\sigma_{\Gamma^*}$
 is the cone generated by $\Gamma^*$. By Corollary~\ref{c:subd} and Example~\ref{e:subdiv}, these cones
  form a fan which we denote
 $\tilde\Sigma_\Delta$.

Applying  Theorem~\ref{t:emb1}, we get

\begin{thm} Associated to the map of fan
$\Sigma_\Delta$ to $\tilde{\Sigma}_\Delta$, the map between toric
varieties $X_{\Sigma_\Delta}\rightarrow X_{\tilde{\Sigma}_\Delta}$
is an embedding, whose image  is a complete intersection given by
the equations $$\prod_{ \langle e_i^*,v_\xi\rangle>0 } x_\xi
-\prod_{ \langle e_i^*,v_\xi\rangle<0 } x_\xi =0,$$ for
$i=1,\dots,k$, where $x_\xi$ are the homogeneous coordinates of
the toric variety $X_{\tilde{\Sigma}_\Delta}$ and $v_\xi$ are the
primitive lattice generators of 1-dimensional cones
$\xi\in\tilde{\Sigma}_\Delta(1)$.
\end{thm}

\begin{rem}\label{r:refl} The exciting part here is that the fan
$\tilde{\Sigma}_\Delta$ consists of the cones over the faces of
the  reflexive polytope
$$ {\rm Conv}(
  \Delta_0^*-e_1-\dots-e_k\cup \Delta_1^*+e_1\cup\cdots\cup\Delta_k^* +e_k),$$
whence
   the noncompact toric variety
$X_{\tilde{\Sigma}_\Delta}$ can be compactified to a Gorenstein
Fano toric variety. Relation of this construction to  the Cayley
trick and Mirror Symmetry will be explored in  \cite{m3}.
\end{rem}

If $l(\Delta)$ denotes the number of lattice points in the
reflexive polytope $\Delta$, then we construct
$(kl(\Delta)-k)$-parameter embedded deformations of $\xsd$
corresponding to the diagram
\begin{equation*}
\begin{array}{ccccc}
\xsd & \subset & { {\mathcal    X}_\Delta} & \subset & X_{{\tilde{\Sigma}_\Delta}}\times\CC^{kl(\Delta)-k}\\
\downarrow & & \downarrow& \swarrow  &\\
\{0\}&\subset &\CC^{kl(\Delta)-k},&&
\end{array}
\end{equation*}
with the total space $ {  \mathcal   X}_\Delta$ in
$X_{\tilde{\Sigma}_\Delta}\times\CC^{kl(\Delta)-k}$   given by the
equations in homogeneous
 coordinates: \begin{equation}\label{e:defFano} \Biggl(x^{ e_i^*}-1+\sum_{\Gamma\prec\Delta}\sum_{u\in\inte({\Gamma})}\lambda_{i,u}
x^{u-\sum_{l=1}^k \langle u,\Gamma_l^*\rangle e_l^*
}\Biggr)\prod_{ \langle e_i^*,v_\xi\rangle<0 }
x_\xi=0\end{equation} for $i=1,\dots,k$, where we combine
deformations   from $(\ref{e:deformh})$ with  $$\tilde
u=u-\sum_{l=1}^k \min \langle u,\Delta_l^*\rangle
e_l^*=u-\sum_{l=1}^k \langle u,\Gamma_l^*\rangle e_l^*$$ for $  u$
in the relative interior of $ {\Gamma}=\{m\in \Delta \mid\langle
m,n\rangle=-1\, \forall\, n\in \Gamma^*\}$, the face of $\Delta$
dual  to $\Gamma^*$.

\begin{rem} The above construction of deformations of Fano toric
varieties can be described without homogeneous coordinates.
Indeed, let $$\sigma=\{( t\Delta^*,t)\mid t\in\RR_{\ge0}\}\subset
N_\RR\oplus\RR$$ be the reflexive Gorenstein cone  and let
$R=(0,1)\in M\oplus\ZZ$. Then the slice of $\sigma$ in the
hyperplane where $R$ has value 1 is the polytope $( \Delta^*,1)$.
From (\ref{e:deldec}) we get the induced decomposition
$$( \Delta^*,1)=(\Delta^*_0,1)+(\Delta_1^*,0)+\cdots+( \Delta_k^*,0).$$
Applying (\ref{e:cons}) we get the cone $$ \tilde\sigma=\langle
(\Delta^*_0,1)-e_1-\dots-e_k,(\Delta_1^*,0)+e_1,\dots, (
 \Delta_k^*,0)+e_k\rangle $$
in $N_\RR\oplus\RR^{k+1}$, which is reflexive Gorenstein of index
$k+1$ by \cite{bb}. By Proposition~\ref{p:comp}  we have embedding
$X_\sigma\hookrightarrow X_{\tilde\sigma} $ given by
$\chi^{R+e_1^*}-\chi^R,\dots,\chi^{R+e_k^*}-\chi^R$. Its
deformation then corresponds to the quotient ring $\CC{
[{\tilde{\sigma}}^{\ve}}\cap
\tilde{M}]\otimes\CC[\underline{\lambda}]/\tilde I$, where
$\tilde{M}=M\oplus\ZZ^{k+1}$ and $\tilde I$ is generated by
$$ \chi^{R+e_i^*}-\chi^R+\sum_{\Gamma\prec\Delta}\sum_{u\in\inte({\Gamma})}\lambda_{i,u}\chi^{    (u,1)-\sum_{l=1}^k \langle
 u,\Gamma_l^*\rangle e_l^*  },$$
 for $i=1,\dots,k$.
For the Gorenstein Fano toric variety $\xsd={\rm Proj}(\CC{ [{
{\sigma}}^{\ve}}\cap  {M}]) $, we get its deformations inside
  $X_{\tilde{\Sigma}_\Delta}={\rm Proj}(\CC{
[{\tilde{\sigma}}^{\ve}}\cap \tilde{M}])$ as  complete
intersections $ {\rm Proj}(\CC{ [{\tilde{\sigma}}^{\ve}}\cap
\tilde{M}]/\tilde I)$ for different values of $\lambda$'s. An
important point of this  remark  is that the Gorenstein cone
construction used in Mirror Symmetry (see \cite{bb}) is  a
particular case of Altmann's construction in \cite{al2}!
\end{rem}

 Finally, we relate the above deformations of Gorenstein Fano toric varieties $X_\Delta$ to deformations of
  their partial crepant resolutions arising from Minkowski sum decompositions of a subdivision of the boundary of $\Delta^*$.  Notice that  that a refinement $\tilde\Sigma'_\Delta$ of the fan $\tilde{\Sigma}_\Delta$
whose $1$-dimensional cones are generated by the lattice points on
the boundary of the reflexive polytope in Remark~\ref{r:refl}
corresponds to a partial crepant resolution
$X_{\tilde\Sigma'_\Delta}\rightarrow X_{\tilde{\Sigma}_\Delta}$.
Intersection of the fan $\tilde\Sigma'_\Delta$ with the linear
subspace $N_\RR\subset \tilde N_\RR$ gives a subdivision
$\Sigma'_\Delta$ of the fan $\Sigma_\Delta$ and a commutative
diagram:
$$ \begin{CD}
 &X_{ \Sigma'_\Delta}& \,\,\hookrightarrow\,\, &  X_{\tilde\Sigma'_\Delta}  &  \\
    & @VVV     @VVV  & \\
 &X_{\Sigma_\Delta}& \,\,\,\hookrightarrow \,\,\,&
 X_{\tilde\Sigma_\Delta}.&
\end{CD}   $$
This induces a morphism between the deformation families of
complete intersections:
$$
\begin{array}{ccccc}
 { {\mathcal    X }'_\Delta} &  \rightarrow &{ {\mathcal    X}_\Delta}\\
\downarrow & & \downarrow&  \\
\CC^{kl(\Delta)-k}&=&\CC^{kl(\Delta)-k},
\end{array}
$$ where ${ {\mathcal    X }'_\Delta}\subset X_{\tilde\Sigma'_\Delta}\times\CC^{kl(\Delta)-k}$ is given by
the equations as in $(\ref{e:defFano})$ but with  homogeneous
coordinates $x_{\xi'}$, for $\xi'\in \tilde\Sigma'_\Delta(1)$, of
the toric variety  $X_{\tilde\Sigma'_\Delta}$.

\begin{ex} In some cases, as in this example, the family ${ {\mathcal    X }'_\Delta}$
is a complete family of deformations of $X_{ \Sigma'_\Delta}$.
  Consider the  fan $\Sigma'_\Delta$  obtained by subdividing   the
reflexive polytope $\Delta^*$ in Example~\ref{ex:2nd} using 6
lattice points in the relative interiors of the facets. This toric
variety $X_{\Sigma'_\Delta}$ is the crepant resolution of the
Gorenstein Fano toric variety which is mirror of $\PP^2$. One can
show that the Minkowski sum decomposition of the polyhedral
complex in Example~\ref{ex:2nd} is $\Sigma'_\Delta$-separated
either directly or by  using Corollary~\ref{c:subd} and
Example~\ref{e:subdiv}. This gives the corresponding fan
$\tilde\Sigma'_\Delta$. Note that $k=2$ and $l(\Delta)=4$ in this
case, whence $kl(\Delta)-k=6$. On the other hand, the dimension
formula in Theorem~\ref{t:infdef} tells us that $\dim
H^1(X_{\Sigma'_\Delta},{\cal T}_{X_{\Sigma'_\Delta}})=6$. One can
show that the Kodaira-Spencer map is an isomorphism for the family
${ {\mathcal    X }'_\Delta}\rightarrow\CC^6$.
\end{ex}

Applications of deformations of Gorenstein Fano toric varieties to
deformations of Calabi-Yau complete intersections and Mirror
Symmetry  will follow in \cite{m3}.

\section{Appendix. Complete intersection toric
ideals.}\label{s:ap}

Complete intersection binomial prime ideals in a polynomial ring
have been characterized in \cite{fs} in terms of mixed dominating
matrices. Here, we will review these results, which we used in the
proof of Proposition~\ref{p:iso}.

  Let $L\subseteq\ZZ^n$ be a sublattice. Denote ${\bf x}^u=\prod
  x_i^{u_i}$ for $u=(u_1,\dots,u_n)$.
Then the {\it lattice ideal} $I_L\subseteq\CC[x_1,\dots,x_n]$
associated to $L$ is the ideal
$$I_L=\langle {\bf x}^u-{\bf x}^v\mid u,v\in
\NN^n \text{ and } u-v\in L\rangle.$$ Equivalently, $I_L=\langle
{\bf x}^{l_+}-{\bf x}^{l_-}\mid l\in L \rangle$, where
  $l=l_+-l_-$ for $l_+,l_-\in\NN^n$ which  have
  disjoint support (i.e., no coordinate is positive in both vectors).
A prime lattice ideal is called a {\it toric ideal}. Every prime
binomial ideal in $\CC[x_1,\dots,x_n]$ is toric by Proposition
1.1.11 in \cite{torvar}.

 By \cite[Theorem~7.4]{ms}, the lattice ideal $I_L$ is prime if and only if $L$ is
saturated in $\ZZ^n$ (i.e., if $s\cdot u\in L$ for $u\in\ZZ^n$ and
$s\in\ZZ_{>0}$, then $u\in L$). And, by Proposition~7.5 in
\cite{ms}, the Krull dimension of $S/I_L$ equals $n-{\rm
rank}(L)$.

Let $l^1,\dots,l^k$ be a lattice basis of $L$. Then the ideal
$I_{\bf L}=\langle {\bf x}^{l^i_+}-{\bf x}^{l^i_-}\mid i=1,\dots,k
\rangle$ is called a {\it lattice  basis ideal}. It is related to
the lattice ideal by the formula
 $I_L= I_{\bf L}:\langle x_1\cdots x_n\rangle^\infty$ (see
 \cite[Lemma~7.6]{ms}). The equality  $ I_L= I_{\bf L} $ holds when the entries in the matrix made of the vectors
 $l^1,\dots,l^k$ satisfy a certain sign pattern.

\begin{defn} \cite{fs}
A matrix $M$ is called {\it mixed} if every row of $M$ contains
both positive and negative entries. A matrix $M$ is called {\it
dominating} if it does not contain a square mixed submatrix.
\end{defn}

\begin{thm} \cite[Theorem~2.9]{fs} Let $L\subseteq\ZZ^n$ be a sublattice such
that $L\cap\NN^n=\{0\}$ and let $M$ be the $k\times n$ matrix
whose rows consist of a lattice basis $l^1,\dots,l^k$   of $L$.
Then $ I_L= I_{\bf L} $ if and only if $M$ is mixed dominating.
\end{thm}

As a consequence (see Corollary~2.4 and 2.10 in \cite{fs}):

\begin{cor}\label{c:toricci} If  $L=\bigoplus_{i=1}^k\ZZ l^i$ is a saturated sublattice in $\ZZ^n$
such that $L\cap\NN^n=\{0\}$, then
  $ I_L= I_{\bf L}=\langle {\bf x}^{l^i_+}-{\bf
x}^{l^i_-}\mid i=1,\dots,k \rangle$ if and only if the matrix $M$
with rows $l^1,\dots,l^k$ is mixed dominating. In this case, the
binomial sequence ${\bf x}^{l^i_+}-{\bf x}^{l^i_-}$,
$i=1,\dots,k$, is regular and   the toric ideal  $ I_L=  I_{\bf
L}$ is a complete intersection.
\end{cor}

\begin{rem} \label{r:satur} If $L=\bigoplus_{i=1}^k\ZZ l^i$ is a  sublattice in $\ZZ^n$ and the matrix $M$
with rows $l^1,\dots,l^k$ is mixed dominating, then
$L\cap\NN^n=\{0\}$ by Corollary~2.7 in \cite{fs}.

\end{rem}

\end{document}